\documentclass[11pt,reqno]{amsart}

\numberwithin{equation}{section}
\newtheorem{theorem}{\sc Theorem}[section]
\newtheorem{lemma}{\sc Lemma}[section]

\newtheorem{remark}{\sc Remark}
\newtheorem{definition}{\sc Definition}[section]

\title[Motion of elastic solids inside of a viscous fluid]
      {On the motion of an elastic solid inside of an incompressible
       viscous fluid}

\author[D. Coutand and S. Shkoller]{}

\email{coutand@math.ucdavis.edu}
\email{shkoller@math.ucdavis.edu}

\begin{document}

\maketitle

\centerline{\scshape   Daniel Coutand and Steve Shkoller}
 \medskip

  {\footnotesize \centerline{ Department of Mathematics }
  \centerline{ University of California at Davis } \centerline{Davis, CA
   95616 } }

\begin{abstract}
The motion of an elastic solid inside of an incompressible viscous fluid is
ubiquitous in nature.  Mathematically, such motion is described by
a PDE system that couples the  parabolic and hyperbolic
phases, the latter inducing a loss of regularity which has left the 
basic question of existence open until now. 

In this paper, we prove the existence and uniqueness of such motions (locally in
time),
when the elastic solid is the linear Kirchhoff elastic material.
The solution is found using a topological fixed-point theorem that 
requires the analysis of a linear problem consisting of the coupling
between the time-dependent Navier-Stokes equations set in Lagrangian variables and the linear equations of 
elastodynamics, for which we prove the existence of a unique weak solution. We then
establish the regularity of the weak solution; this regularity is obtained in
function spaces that scale in a hyperbolic fashion in both the fluid and solid
phases.  Our functional framework is optimal, and provides the a priori
estimates necessary for us to employ our fixed-point procedure.

\end{abstract}

%\tableofcontents

\section{Introduction}
We are concerned with establishing the existence (and uniqueness) of solutions
for the equations of
motion of linearly elastic solids moving and interacting with an incompressible 
viscous fluid, with the natural conditions of continuity of the velocity fields 
and normal components of the stress tensors along the moving interface between
the two materials.  

The analysis of interacting fluid-structure problems has been the subject of 
active research since the late nineties. As of now, only the question of the
possible  motion of a solid
inside of a viscous flow, in which the solid is either {\it rigid}  or 
consists of a {\it finite number of modes}, has been settled. 
In \cite{GrMa}, existence and uniqueness (locally in time) of smooth solutions 
has 
been obtained using a Lagrangian framework, for the rigid body case, provided 
that 
the rigid disk is sufficiently heavy. In \cite{DeEs1999}, for the same problem, but
with an arbitrary number of rigid solid bodies, existence of at least one weak
solution has been established in an Eulerian formulation by a
global variational approach; their  result holds for all time in two space
dimensions as
long as no collisions occur between solids or with the boundary, and is local 
in 
time for the three dimensional case. In \cite{DeEsGrLe}, by generalizing the methods
of \cite{DeEs1999}, the case of an elastic body following the linear Kirchhoff 
law, 
with the important restrictions of allowing only a {\it finite number of modes}, 
and a relaxation of the continuity of the normal stress along the boundary of 
the 
solids, has been considered. The above list of references for contributions to 
this
area is by no means exhaustive; see for instance \cite{Conca}, \cite{FlOr}, 
\cite{Gunz}.  Note also that the related problem of the free fall of a rigid
body in a Stokes flow in the full space has been considered in \cite{Wein}, 
for the stationary case, and in \cite{Serre} for the stationary as well as the 
time-dependent case.

More recently, the interaction of a viscous incompressible flow with an 
elastic plate (without the restriction of a finite number of modes), whose 
constitutive law comprises a {\it parabolic} hyperviscosity term in
the plate, has been studied in \cite{ChDeEsGr}.
We remark that this additional hyperviscosity term  is of crucial
importance in that study.  (Note also that two dimensional plate models that
approximate thin three dimensional structures usually contain fourth-order
operators arising from bending stresses, whereas models of elastic solids 
have only second-order operators;  as such, plate models can provide better
a priori control for the motion of the material interface.)

In the {\it steady-state} situation in which both phases are governed by
elliptic operators, \cite{Gr} has
obtained an existence result (for the case that the solid follows the 
nonlinear Saint Venant-Kirchhoff law) by the use of a fixed-point method that 
iterates between fluid and solid phases.  
This approach is
indeed natural for the steady-state problem since the analysis can make use of
elliptic regularity theory.  For the dynamic problem, however, such an iteration
procedure appears to 
{\it fails} because of a consequent loss of regularity induced by either a 
fluid-solid-fluid iteration or a solid-fluid-solid iteration.  This loss of
regularity is due to the fact that
hyperbolic and parabolic systems do not have the same regularity requirements 
and 
properties, which is in fact the heart of the difficulty in the coupling  of
the two phases.

Whereas the coupling between the Navier-Stokes equations and the linear 
Kirchhoff law is perhaps the most fundamental problem to consider in regards 
to the motion and interaction of an elastic body in a viscous incompressible 
fluid, none of the methods that have been developed to date can handle this 
system, mostly because of 
the differences between parabolic and hyperbolic regularity, i.e. in both the
requirements on the function spaces for the prescribed data, as well as the 
functional framework of the solution space.

We now come to the formulation of the problem.  The motion of the fluid is described
by the time-dependent incompressible Navier-Stokes equations, while the deformation 
of the solid body is governed by the linear Kirchhoff equations.  The two models are 
coupled along the moving material interface by imposing the continuity of the 
normal component of the stress tensors as well as the particle displacement fields.
From the point-of-view of mathematical analysis, the Navier-Stokes equations
are traditionally studied in the Eulerian (or spatial) description, while the
elastic body is studied in the Lagrangian (or material) frame.  Because the
material interface is fixed in the Lagrangian representation, we shall study
this problem entirely in material coordinates. This Lagrangian framework also has 
the advantage of keeping the hyperbolic problem (where the loss of regularity 
occurs) linear, which is of paramount importance here. Note, however, that a
semi-linear elastic system, as for some plate or shell models (see for 
instance \cite{PGC1}), 
can be handled without any difficulty by our methodology. The question of 
existence for the case of a quasilinear elasticity law can also be obtained 
(and shall be addressed in later work), requiring a smoother functional 
framework leading to more compatibility conditions at the origin. 

Let us now set the equations. Let $\Omega \subset{\mathbb R}^3$ denote an open, 
bounded, connected and smooth domain with smooth boundary $\partial\Omega$ which 
represents the fluid container in which both the solid and fluid move.
Let $ \overline{\Omega ^s(t)}\subset  \Omega $ denote the closure of an open
and bounded
subset representing the solid body at each instant of time $t\in [0,T]$ with
$\Omega^f(t):= \Omega / \overline{\Omega^s(t)}$ denoting the fluid domain at 
each $t\in [0,T]$. 
Note that in our analysis $\Omega^s (t)$ is not necessarily connected,
which allows us to handle the case of several elastic bodies moving in the 
fluid.

\begin{remark}
If a function $u$ is defined on all of $\Omega$, we will denote 
$u^f=u\ 1_{\overline \Omega_0^f}$ and
$u^s=u\ 1_{\overline\Omega_0^s}$. This allows us to indicate from which phase
the traces on 
$$\Gamma(0):= \overline{\Omega^f(0)}\cap \overline{\Omega^s(0)} $$ 
of various discontinuous terms arise, and also to
specify functions that are associated with the fluid and solid phases.
\end{remark}

For each $t\in (0,T]$, we wish to find the location of these domains
inside $\Omega$, the
divergence-free velocity field $u^f(t, \cdot)$ of the fluid, the fluid pressure
function $p(t, \cdot)$ on $\Omega^f(t)$,  the fluid Lagrangian 
volume-preserving configuration
$ \eta ^f(t, \cdot ):\Omega ^f(0)=\Omega_0^f \rightarrow \Omega ^f(t)$, 
and the elastic  Lagrangian configuration field 
$\eta^s(t, \cdot):\Omega^s(0)=\Omega_0^s \rightarrow \Omega ^s(t)$ such that
\begin{subequations}
  \label{ns-elastic}
\begin{align}
\Omega&=  \eta^s(t,\overline{\Omega^s_0}) \cup \eta^f(t,\Omega^f_0) \,,
         \label{ns-elastic.a}\\
\intertext{where}
\eta^f_t(t,x) &= u^f(t, \eta^f(t,x)) \,,
         \label{ns-elastic.b}\\
\intertext{and $u^f$ solves the Navier-Stokes equations in $\Omega^f(t)$}
u^f_t + (u^f\cdot \nabla)u^f  &=  \operatorname{div} T^f + f_f \,,
         \label{ns-elastic.c}\\
   \operatorname{div} u^f &= 0     \,,
         \label{ns-elastic.d}\\
\intertext{with} T^f&=\nu\ \operatorname{Def} u^f -p\ \text{I} \,,\label{fluidlaw}\\
\intertext{and $\eta^s$ solves the elasticity equations on $\Omega^s(0)$}
\ddot{\eta}^s  &= \operatorname{div} T^s + f_s \,,
         \label{ns-elastic.e}\\
\intertext{with} T^s&= \lambda\ \text{Trace}(\nabla\eta^s-I)I\ +\mu\ (\nabla\eta^s+{\nabla\eta^s}^T-2\ I) \,,\\
\intertext{and where the equations are coupled together by the continuity of
the normal component of stress along the material interface 
$\Gamma(t):= \overline{\Omega^s}(t) \cap \overline{\Omega^f}(t)$ expressed 
in the Lagrangian representation on}
\Gamma_0&:=\Gamma(0)\nonumber\\
\intertext{as}
T^s \ N
&=  
[T^f \circ \eta^f] \ [( \nabla \eta^f)^{-1} \ N ]
\,,\\
\intertext{and the continuity of particle displacement fields along
$\Gamma_0$}
         \label{ns-elastic.f}
\eta^f &=  \eta^s
\,,\\
\intertext{together with the initial conditions}
   u(0,x) &= u_0(x)     \,,
         \label{ns-elastic.h}\\
   \eta(0,x) &= x     \,,
         \label{ns-elastic.i}
\intertext{and the Dirichlet (no-slip) condition on the boundary 
$\partial\Omega$ of the container} 
u^f&=0\,,
\end{align}
\end{subequations}
where $\nu>0$ is the kinematic viscosity of the fluid, $\lambda>0$ and $\mu>0$ denote the Lam\'e constants of the elastic material, 
$N$ is the outward unit normal to $\Gamma_0$ and $\operatorname{Def}u$ is twice the rate of deformation tensor 
of $u$, given in coordinates by $u^i,_j + u^j,_i$.  All Latin indices run 
through $1,2,3$, the Einstein  summation convention is employed, and indices 
after commas denote partial derivatives.

We now briefly outline the proof. As the solid and fluid phases are naturally 
expressed in the Lagrangian and Eulerian framework, respectively, we begin 
by transforming the fluid phase into Lagrangian coordinates, leading us to the 
system of equations (\ref{nsl}) of Section \ref{2}. 
This system of PDE is both 
parabolic (in the fluid) and hyperbolic (in the solid) in character;
hence, one of the fundamental difficulties that must be overcome is 
an appropriate functional framework accommodating both features. 
Sections \ref{3} and \ref{4} are devoted to the setting of our functional 
framework, which appears to be of hyperbolic-type  in both solid and fluid 
phases, and is necessitated by the estimate of the elastic energy. 
This hyperbolic scaling in turn requires the initial data to possess more
regularity, and thus produces more compatibility conditions in the fluid phase 
than if a parabolic scaling were used  (as seen in the statement of the 
existence theorem in Section \ref{5}). Whereas the choice of working in 
Eulerian or Lagrangian variables may seem arbitrary, at the level of the 
functional framework, it appears that the problem truly requires this
hyperbolic functional framework for both phases, regardless of the choice
of spatial or material coordinates.

In order to solve (\ref{nsl}), we use a fixed-point approach, where we solve
the linear system
(\ref{linear}) for the Lagrangian velocity $w$, the coefficients $a^i_j(\eta)$ 
coming from the flow map $\eta$ of a {\it given} velocity $v$. 
The study of the regularity of the solutions to 
this problem, which constitutes the main part of this paper, is given in 
Sections \ref{8} and \ref{9}. 
It appears that the regularity theory for (\ref{linear}) cannot be obtained 
directly by solving the problem with the actual coefficients $a_i^j (\eta)$.
In Section \ref{7} we explain the smoothing process for the problem:  we
introduce smoothed velocity fields $v_n$ which provide us with smoothed
coefficients $a_i^f(\eta_n)$ (which we denote generically by $\tilde v$ and
$\tilde a$). We also present two versions of what we term the {\it Lagrange 
multiplier lemma} (which associates a pressure function to the weak solution) 
that will be of basic use throughout this paper.

We study in Section \ref{8} the existence of weak solutions $\tilde w$ to 
(\ref{linear}) (with regularized coefficients), as the limit of penalized 
problems. Whereas these penalized problems are not necessary merely to obtain
existence of weak solutions, they are of {\it paramount} importance in getting 
the appropriate regularity results for $\tilde w_t$ and $\tilde w_{tt}$, 
the primary reason being that the pressure associated to
(\ref{linear}) with the Dirichlet boundary condition cannot be obtained simply 
from the variational form of the problem, and requires the study of the time 
differentiated problem in order to get more information on $\tilde w_t$ 
(which would need to be
in $L^2(0,T;H^{-1}(\Omega;{\mathbb R}^3))$ for the Lagrange multiplier lemma). 
Unfortunately, this time differentiated problem contains $p \circ \eta$ in its 
formulation, which leads to a  circular argument, and thus explains the need 
for the penalized problem. We then obtain the regularity for the problem by
the energy inequality for $\tilde w_{tt}$ and some difference-quotient 
inequalities for $\tilde w_t$ and $\tilde w$ carried-out in Lagrangian variables
in a neighborhood of the interface $\Gamma_0$. This, in turn, provides us 
with an estimate for the {\it trace} of $\tilde w$ and $\tilde w_t$ on 
$\Gamma_0$, which after a return to the
Eulerian variables for the fluid phase, immediately provides the regularity in 
the fluid domain. The regularity in the solid phase is then obtained in
a straightforward manner from 
elliptic regularity and the already-obtained trace estimate.  We note 
that the estimates proved at this stage blow-up as the regularized 
coefficients tend to the true coefficients, i.e., as the regularization
parameter tends to zero.

For this reason,
in Section \ref{9}, we obtain a different set of estimates (founded 
upon interpolation inequalities) for the solutions of the regularized problems, 
and conclude that the norms of the 
regularized solutions are actually uniformly bounded in the appropriate 
spaces, which thus provides by weak convergence, a solution to 
(\ref{linear}) with the appropriate  a priori estimates.

Finally, we conclude the proof of the existence theorem in Sections \ref{10} and
\ref{11} by means of the Tychonoff fixed-point theorem. Although it might be
possible to employ the Schauder theorem instead, it appears that the strong 
convergence requirements of the Schauder theorem are not very 
convenient to write and are, in particular, unnecessary for the use of
the Tychonoff theorem. 

Uniqueness is obtained with more regularity on the  forcing functions and 
initial data in Section
\ref{12}, in order to get information on the second-time derivative of the
pressure $q_{tt}$ of the solution built, that we do not have with the assumptions of Theorem \ref{main}. 

\section{Notational simplification}
\label{1bis}
Although a fluid with a Neumann (free-slip) boundary condition indeed obeys the 
constitutive law (\ref{fluidlaw}), it turns out that the notation
is substantially simplified (particularly in Section \ref{8} wherein we 
analyze the twice differentiated-in-time  problem in Lagrangian coordinates)
if we replace (\ref{fluidlaw}) with
\begin{equation}
\label{newfluid}
T^f=\nu \nabla u^f-p \text{I};
\end{equation}
this amounts to replacing the energy 
$\displaystyle\int_{\Omega_0^f}\operatorname{Def}u^f : \operatorname{Def} v$ 
by $\displaystyle\int_{\Omega_0^f}\nabla u^f:\nabla v$, 
which is an equivalent form when $u^f=0$ on $\partial \Omega$ due to the
well-known Korn inequality. Henceforth, we shall take (\ref{newfluid}) as the
fluid constitutive law.

\section{Lagrangian formulation of the problem}
\label{2}
In regards to the forcing functions, we shall use the convention of denoting 
both the fluid forcing $f_f$ and the solid forcing $f_s$ by the same letter $f$. 
Since $f_f$ has to be defined in 
$\Omega$ (because of the composition with $\eta$), and $f_s$ must be defined 
in $\Omega_0^s$, we will assume that the forcing $f$ is defined over the entire
domain $\Omega$.
 
Let 
\begin{equation}\label{a}
a(x) = [\nabla\eta^f(x)]^{-1}, 
\end{equation}
where $(\nabla \eta^f(x))^i_j = \partial (\eta^f)^i/\partial x^j(x)$
denotes the matrix of partial derivatives of $\eta^f$.  
Clearly, the matrix $a$ depends on $\eta$ and we shall sometimes use the 
notation $a^i_j(\eta)$ to denote the formula (\ref{a}).

Let $v=u\circ \eta$ denote
the Lagrangian or material velocity field, $q=p \circ \eta$ is the Lagrangian 
pressure function (in the fluid), 
and $F= f^f \circ \eta^f$ is the fluid forcing function in
the material frame.  
Then, as long as no collisions occur between the solids (if there are initially
more than one) or between a solid and $\partial\Omega$, the system
(\ref{ns-elastic}) can be reformulated as 
\begin{subequations}
  \label{nsl}
\begin{alignat}{2}
\eta_t &=v&\ &\text{in} \ \ (0,T)\times \Omega \,, 
         \label{nsl.a}\\
v^i_t - \nu (a^j_l a^k_l v^i,_k),_j + (a^k_i q),_k &= F^i 
&&\text{in} \ \ (0,T)\times \Omega_0^f \,, 
         \label{nsl.b}\\
   a^k_i v^i,_k &= 0     &&\text{in} \ \ (0,T)\times \Omega_0^f \,, 
         \label{nsl.c}\\
v^i_t - [c^{ijkl}\int_0^t v^k,_l],_j &= f^i
&&\text{in} \ \ (0,T)\times \Omega_0^s \,, 
         \label{nsl.d}\\
\nu\ v^i,_k a^k_l a^j_l N_j - q a^j_i N_j &= 
c^{ijkl}\int_0^t v^k,_l\ N_j
&&\text{on} \ \ (0,T)\times \Gamma_0 \,, 
         \label{nsl.e}\\
v(t,\cdot)& \in H^1_0 (\Omega;{\mathbb R}^3) &&\text{a.e. in } \ \ (0,T)\,, \label{nsl.f}\\
   v &= u_0  
 &&\text{on} \ \ \Omega_0\times \{ t=0\} \,, 
         \label{nsl.g}\\
   \eta &= \text{Id}     
 &&\text{on} \ \ \Omega_0\times \{ t=0\} \,, 
         \label{nsl.h}
\end{alignat}
\end{subequations}
where $N$ denotes the outward-pointing unit normal to $\Gamma_0$ (pointing
into the solid phase), and
$$c^{ijkl}=\lambda \delta^{ij}\delta^{kl}+\mu (\delta^{ik}\delta^{jl}
+\delta^{il}\delta^{jk})\ . $$
Throughout the paper, all Greek indices run through $1,2$ and 
all Latin indices run through $1,2,3$.  Note that the continuity of the
velocity (\ref{ns-elastic.f}) along the interface is satisfied in the sense of
traces on $\Gamma_0$ by condition (\ref{nsl.f}), whereas the continuity of 
the normal stress along the interface is represented by (\ref{nsl.e}).

\begin{remark}
The case in which  the viscosity or Lam\'e coefficients are variable functions
depending on $x \in \Omega$ and satisfying the usual 
assumptions, can be handled by our methodology without any supplementary 
mathematical difficulties.
\end{remark}

%\begin{remark}
%Remember that we use (\ref{newfluid}) as the constitutive law for the fluid.
%\end{remark} 

\section{Notation and conventions}
\label{3}
We begin by specifying our notation for certain vector and matrix operations.
\begin{itemize}
\item[] We write the Euclidean inner-product between two vectors $x$ and $y$ 
as $x\cdot y$, so that $x\cdot y=x^i\ y^i$.
\item[] The transpose of a matrix $A$ will be denoted by $A^T$, {\it i.e.}, 
$(A^T)^i_j=A^j_i$.
\item[] We write the product of a matrix $A$ and a vector $b$ as $A\ b$, 
{\it i.e}, $(A\ b)^i=A^i_j b^j$.
\item[] The product of two matrices $A$ and $S$ will be denoted by 
$A\cdot S$, {\it i.e.}, $(A\cdot S)^i_j=A^i_k\ S^k_j$.
\item[] The trace of the product of two matrices $A$ and $S$ will be denoted by 
$A: S$, {\it i.e.}, $ A:S=\operatorname{Trace}(A\cdot S)=A^i_j\ S^j_i$.
\end{itemize}

For $s\ge 0$ and a Hilbert space $(X,\|\cdot\|_X)$, $H^s(\Omega; {\mathbb R}^3)$ denotes the
Sobolev space of $\mathbb R^3$-valued functions with $s$ distributional derivatives in 
$L^2(\Omega;{\mathbb R}^3)$, while $L^2 (0,T; X)$ denotes the equivalence class of functions 
which are measurable and have finite
$\|\cdot\|_{L^2}$-norm, where 
$\|f\|^2_{L^2(0,T;X)}= \int_0^T\|f(t)\|^2_X dt$.

We also set $H^1_{\partial \Omega} (\Omega_0^f;{\mathbb R}^3)=\{u\in 
H^1(\Omega_0^f;{\mathbb R}^3)|\ u=0\ \text{on}\ \partial \Omega\}\ .$

For $T>0$, we set
\begin{align*}
%V^1(T)&= L^2(0,T;  H^1_0(\Omega;{\mathbb R}^3)), \\
V^2_f(T)&=\{ w \in L^2(0,T;  H^2(\Omega_0^f;{\mathbb R}^3) )\ | \ 
w_t \in L^2(0,T;  H^1(\Omega_0^f;{\mathbb R}^3))\\
&\qquad w_{tt} \in L^2(0,T;  L^2(\Omega_0^f;{\mathbb R}^3)) \},  \\
V^3_f(T)&=\{ w \in L^2(0,T;  H^3(\Omega_0^f;{\mathbb R}^3)) \ | \ 
w_t \in L^2(0,T;  H^2(\Omega_0^f;{\mathbb R}^3))\ |\ \\
&\qquad w_{tt}\in  L^2(0,T;  H^1(\Omega_0^f;{\mathbb R}^3))\},  \\
%F^1_f(T)&=\{ w \in L^2(0,T;  H^1(\Omega_0^s;{\mathbb R}^3)) \ | \ 
%w_t \in L^2(0,T;  L^2(\Omega_0^s;{\mathbb R}^3))\ |\ \\
%& \qquad w_{tt} \in L^2(0,T;  H^1_{\partial\Omega}(\Omega_0^f;{\mathbb R}^3)') \},  \\
V^2_s(T)&=\{ w \in L^2(0,T;  H^2(\Omega_0^s;{\mathbb R}^3)) \ | \ 
w_{t} \in L^2(0,T;  H^1(\Omega_0^s;{\mathbb R}^3))\ |\ \\
&\qquad w_{tt} \in L^2(0,T;  L^2(\Omega_0^s;{\mathbb R}^3))  \}\ ,  \\
V^3_s(T)&=\{ w \in L^2(0,T;  H^3(\Omega_0^s;{\mathbb R}^3) )\ | \ 
w_t \in L^2(0,T;  H^2(\Omega_0^s;{\mathbb R}^3))|\\
&\qquad w_{tt} \in L^2(0,T;  H^1(\Omega_0^s;{\mathbb R}^3))
 \}.
\end{align*}

We will solve (\ref{nsl}) by a fixed point method, set in an appropriate subset
of $V^3_f (T)\times V^3_s (T)$. We assume in what follows that $v\in V^3_f (T)$
is given in such a way that the matrix $a_i^j(\eta)$ associated to the flow
$\eta$ of this velocity field $v$ is well-defined.  

We then introduce the space (of weak solutions)
\begin{align*}
{\mathcal V}_v([0,T]) = &
\{ w \in L^2(0,T;  L^2(\Omega;{\mathbb R}^3)) \ | \int_0^\cdot w \in L^2(0,T;  H^1 (\Omega;{\mathbb R}^3)),\\
&w\in L^2(0,T;  H^1(\Omega_0^f;{\mathbb R}^3)),\ \  
 a_i^j w^i,_j=0\ \text{in}\ [0,T]\times \Omega_0^f,\ \ w=0\ \text{on}\ 
\partial \Omega \}.
\end{align*}

Note that we impose the condition 
$\displaystyle\int_0^\cdot w \in L^2(0,T;  H^1 (\Omega;{\mathbb R}^3))$ to ensure continuity
of the displacement field, in the sense of traces, between the solid and fluid
phases along the interface $\Gamma_0$.  
We will also denote  for $t\in [0,T]$
$${\mathcal V}_v(t)=\{ \psi \in H^1_0(\Omega; {\mathbb R}^3) \ | \ a_i^j (t) \psi^i,_j=0 \ 
\text{in}\ \Omega_0^f \}\ .$$

Furthermore, we will need the space
\begin{align*}
{\mathcal W}([0,T]) = &
\{ w \in L^2(0,T;  L^2(\Omega;{\mathbb R}^3)) \ | \int_0^\cdot w \in L^2(0,T;  H^1 (\Omega;{\mathbb R}^3)),\\
&w\in L^2(0,T;  H^1(\Omega_0^f;{\mathbb R}^3))\ \ w=0\ \text{on}\ 
\partial \Omega \},
\end{align*}
with the ``divergence-free'' constraint removed.

In order to specify the initial data for the weak formulation, we 
introduce the space
$$L^2_{div,f}
=\{ \psi \in L^2(\Omega; {\mathbb R}^3) \ | \ \operatorname{div}\psi=0 \ 
\text{in}\ \Omega_0^f, \ \  \psi\cdot N=0\ \text{on}\ \partial \Omega \}\ ,$$
which is endowed with the $L^2(\Omega; {\mathbb R}^3)$ scalar product.

The space of velocities, $X_T$, is defined as the following separable Hilbert space:
\begin{align}
X_T= \{ u\in L^2(0,T; H^1_0(\Omega;{\mathbb R}^3))\ 
|\  \ (u^f, \int_0^\cdot u^s)\in V^3_f (T)
\times V^3_s (T)\}\ ,
\label{XT}
\end{align}
 endowed with its natural Hilbert norm
\begin{align*}
\|u\|^2_{X_T}&=\|  u \|^2_ {L^2(0,T;H^1(\Omega;{\mathbb R}^3))}
+\|  u \|^2_ {L^2(0,T;H^3(\Omega_0^f;{\mathbb R}^3))}+ 
\|u_t\|^2_ {L^2(0,T;H^2(\Omega_0^f;{\mathbb R}^3))}\nonumber \\
& \displaystyle 
+ \|  \int_0^\cdot u \|^2_ {L^2(0,T;H^3(\Omega_0^s;{\mathbb R}^3))}
+ \|  u \|^2_ {L^2(0,T;H^2(\Omega_0^s;{\mathbb R}^3))}
+ 
\|u_t\|^2_ {L^2(0,T;H^1(\Omega_0^s;{\mathbb R}^3))}\nonumber \\
& + \|u_{tt}\|^2_ {L^2(0,T;H^1(\Omega_0^f;{\mathbb R}^3))}
\,.
\end{align*}

The existence of solutions to (\ref{nsl}) will be obtained in the following 
separable Hilbert space 
\begin{align*}
Y_T= \{& (u,p)\in X_T\times L^2(0,T;H^2(\Omega_0^f;{\mathbb R}))|\ p_t\in L^2(0,T;H^1(\Omega_0^f;{\mathbb R}))\}\ , 
\end{align*}
 endowed with its natural Hilbert norm
\begin{align*}
\|(u,p)\|^2_{Y_T}=\|  u \|^2_ {X_T}
+ \|  p \|^2_ {L^2(0,T;H^2(\Omega_0^f;{\mathbb R}))}+ 
\| p_t\|^2_ {L^2(0,T;H^1(\Omega_0^f;{\mathbb R}))}
\,.
\end{align*}

\begin{remark}
Note well that our method does not require any a priori knowledge of the 
regularity of the second time derivative of the pressure function $p_{tt}$; 
this is due to the Dirichlet boundary condition on 
$\partial\Omega$ as well as the Lagrangian representation of the problem
that we employ.
\end{remark}

We shall also need $L^{\infty}$-in-time control of certain norms of the
velocity, which necessitates the use of the  following closed subspace of $X_T$:
\begin{align*}
W_T=\{u\in X_T|\ &\ u_{tt}\in L^{\infty}(0,T;L^2(\Omega;{\mathbb R}^3)),\ u_t\in L^{\infty}(0,T;H^1(\Omega_0^s;{\mathbb R}^3)),\\
& u\in L^{\infty}(0,T;H^2(\Omega_0^s;{\mathbb R}^3)),\  \int_0^{\cdot} u\in L^{\infty}(0,T;H^3(\Omega_0^s;{\mathbb R}^3))\}\ ,
\end{align*}
endowed with the following norm
\begin{align*}
\|u\|^2_{W_T}&=\|  u \|^2_ {X_T}
+\|  u_{tt} \|^2_ {L^{\infty}(0,T;L^2(\Omega;{\mathbb R}^3))}+ 
\|\int_0^{\cdot} u\|^2_ {L^{\infty}(0,T;H^3(\Omega_0^s;{\mathbb R}^3))}\nonumber \\
& \displaystyle 
+ \|  u \|^2_ {L^{\infty}(0,T;H^2(\Omega_0^s;{\mathbb R}^3))}
+ \|  u_t \|^2_ {L^{\infty}(0,T;H^1(\Omega_0^s;{\mathbb R}^3))}
\,.
\end{align*}
For some of our estimates, we will also make use of the space
\begin{align*}
Z_T= \{& (u,p)\in W_T\times L^2(0,T;H^2(\Omega_0^f;{\mathbb R}))|\ p_t\in L^2(0,T;H^1(\Omega_0^f;{\mathbb R}))\}\ , 
\end{align*}
endowed with its natural norm
\begin{align*}
\|(u,p)\|^2_{Z_T}=\|  u \|^2_ {W_T}
+ \|  p \|^2_ {L^2(0,T;H^2(\Omega_0^f;{\mathbb R}))}+ 
\| p_t\|^2_ {L^2(0,T;H^1(\Omega_0^f;{\mathbb R}))}
\,.
\end{align*}

Throughout the paper, we shall use
$C$ to denote a generic constant, which may possibly depend on the coefficients
$\nu$, $\lambda$, $\mu$,  or on the initial geometry given by $\Omega$ and
$\Omega_0^f$ (such as a Sobolev constant or an elliptic constant). Similarly,
we will denote by $C(M)$ a generic constant which depends on the same variables
as $C$ as well as on $M$ (which is a variable defined in the next section) 
and $\|u_0\|_{H^5(\Omega_0^f;{\mathbb R}^3)}$, 
$\|f(0)\|_{H^3(\Omega;{\mathbb R}^3)}$ and the fixed time 
$\bar{T}$ for which the forcing functions are defined. We note that
these constants do not blow-up whenever the quantities they depend 
upon remain finite.

For the sake of notational convenience, we will also write
$u(t)$ for $u(t,\cdot)$.

\section{The main theorem} 
\label{5}
\begin{theorem}\label{main}

Let $\Omega\subset {\mathbb R}^3$ be a bounded domain of class $H^3$, and let 
$\Omega_0^s$ be an open set (with a finite number $\ge 1$ of 
connected components) of class $H^4$ such that $\overline{\Omega_0^s}\subset \Omega$. 
Let us denote $\Omega_0^f=\Omega\cap (\overline{\Omega_0^s})^c$.
 Let $\nu>0$, $\lambda> 0$, $\mu>0$ be given. Let 
\begin{subequations}
\label{f_regularity}
\begin{align} 
f \in L^2(0,\bar T; H^2(\Omega;{\mathbb R}^3)), &\ 
f_t \in L^2(0,\bar T; H^1(\Omega;{\mathbb R}^3)),\  f_{tt} \in L^2(0,\bar T; L^2(\Omega;{\mathbb R}^3)),\\
& f(0)\in H^3(\Omega;{\mathbb R}^3)\ .
\end{align}
\end{subequations}
Assume that the initial data  satisfies
$$u_0 \in H^5(\Omega_0^f;{\mathbb R}^3)\cap H^2(\Omega_0^s;{\mathbb R}^3) 
\cap H^1_0(\Omega;{\mathbb R}^3)\cap L^2_{{div},f}$$
as well as  the compatibility conditions 
\begin{subequations}
\label{compatibility}
\begin{align}
& [\nabla u_0^f\  N]_{\operatorname{tan}}=0\ \text{ on }\ \Gamma_0,\ \  w_1=0\ \text{on}\ \partial\Omega,\ \ \nu\triangle u_0^f-\nabla q_0=0\ \text{on}\ \Gamma_0,\\
&[(\nu [\nabla w_1^f\ N]^i + \nu [\ {u_0^f},_k^i {(a^k_la^j_l)}_t(0)\ ]\ N_j)_{i=1}^3]_{\operatorname {tan}} \nonumber\\
&\qquad\qquad\qquad\qquad\qquad\qquad\qquad\qquad = [(c^{ijkl} {u_0^f},_l^k N_j)_{i=1}^3]_{\operatorname{tan}}\ \text{on}\ \Gamma_0,
\end{align}
\end{subequations}
with $q_0\in H^4(\Omega_0^f;{\mathbb R})$ defined by 
\begin{subequations}
\label{defq0}
\begin{align}
\triangle q_0&=\operatorname{div} f(0)+(a_i^j)_t (0) u_0^i,_j\ \text{in}\ \Omega_0^f,\\
q_0&=\nu [\nabla u_0^f\ N]\cdot N\ \text{on}\ \Gamma_0,\label{defq0.b}\\
\frac{\partial q_0}{\partial N} &=f(0)\cdot N+\nu \triangle u_0\cdot N\ \text{on}\ \partial\Omega ,
\end{align}
\end{subequations}
and $w_1\in H^1_0(\Omega;{\mathbb R}^3)\cap H^3(\Omega_0^s;{\mathbb R}^3)\cap H^3(\Omega_0^f;{\mathbb R}^3)$ defined by 
\begin{subequations}
\label{def1}
\begin{align}
 w_1&=\nu\triangle u_0-\nabla q_0 +f(0)\ \text{in}\ \Omega_0^f\\
w_1&=f(0)\ \text{in}\ \Omega_0^s\ .
\end{align}
\end{subequations}
(Note that $({a_i^j})_t|_{t=0}$ depends only on $u_0$ and not on the values 
taken by $u$ at times $t>0$.)

Then there 
exists $T\in (0,\bar T)$ depending on $u_0$, $f$, and $\Omega_0^f$, such that 
there exists a solution $(v,q) \in Z_T$  
of the problem (\ref{nsl}).  Furthermore, 
$\eta \in  C^0([0, T]; H^3(\Omega_0^f; {\mathbb R}^3)
\cap  H^3(\Omega_0^s; {\mathbb R}^3)\cap H^1(\Omega; {\mathbb R}^3))$.
\end{theorem}
\begin{remark}
In Theorem \ref{unique}, assumptions ensuring uniqueness of the solutions are
also given.
\end{remark}

\begin{remark}
If we had not made the notational simplification of Section \ref{1bis}, we would have to modify (\ref{compatibility}) by 
\begin{align*}
& [\operatorname{Def} u_0^f\  N]_{\operatorname{tan}}=0\ \text{ on }\ \Gamma_0,\ \  w_1=0\ \text{on}\ \partial\Omega,\ \ \nu\triangle u_0^f-\nabla q_0=0\ \text{on}\ \Gamma_0,\\
&[(\nu [\operatorname{Def}w_1^f\ N]^i + \nu [\ {u_0^f},_k^i {(a^k_la^j_l)}_t(0) + {u_0^f},_k^l {(a^k_i a^j_l)}_t(0)\ ]\ N_j)_{i=1}^3]_{\operatorname {tan}} \nonumber\\
&\qquad\qquad\qquad\qquad\qquad\qquad\qquad\qquad = [(c^{ijkl} {u_0^f},_l^k N_j)_{i=1}^3]_{\operatorname{tan}}\ \text{on}\ \Gamma_0,
\end{align*}
and (\ref{defq0.b}) would be replaced by $q_0=\nu\ [\operatorname{Def} u_0^f\ N]\cdot N$ on $\Gamma_0$.
\end{remark}

\begin{remark}
The regularity of our solution $v \in W_T$ implies that for
each $t\in (0, T]$ the solid domain $\Omega_s(t)$ is of class $H^3$.  Also, 
although we have stated our results for three-dimensional motion, all
of our results hold in the two-dimensional case as well.
\end{remark}

\begin{remark}
We have stated our results using the convention of Section \ref{2}, wherein
the forcing function $f$ is taken to be defined over the entire domain
$\Omega$.  It is certainly possible to define separate forcing functions for
the solid and fluid phase, in which case we would need the following
regularity:
\begin{align*} 
f_f \in L^2(0,\bar T; H^1(\Omega;{\mathbb R}^3)), &\ 
{f_f}_t \in L^2(0,\bar T; L^2(\Omega;{\mathbb R}^3)),\  {f_f}_{tt} \in L^2(0,\bar T; H^1_{\partial\Omega}(\Omega;{\mathbb R}^3)'),\\
f_s \in L^2(0,\bar T; H^2(\Omega_0^s;{\mathbb R}^3)), &\ 
{f_s}_t \in L^2(0,\bar T; H^1(\Omega_0^s;{\mathbb R}^3)),\  {f_s}_{tt} \in L^2(0,\bar T; L^2(\Omega_0^s;{\mathbb R}^3)),\\
& f_f(0)\in H^3(\Omega_0^f;{\mathbb R}^3)\ .
\end{align*}
The compatibility condition $\nu\triangle u_0^f-\nabla q_0=0$ on $\Gamma_0$ in
(\ref{compatibility}) would be replaced by $\nu\triangle u_0^f-\nabla q_0 
+f_f(0)=f_s(0)$ on $\Gamma_0$. In the definition of $q_0$, 
$f_f(0)$ replaces $f(0)$ and in the definition of $w_1$, $f(0)$ is replaced 
by $f_f(0)$ and $f_s(0)$ respectively in $\Omega_0^f$ and $\Omega_0^s$.
\end{remark} 

\begin{remark}
Note that the supplementary regularity condition for $u(0)$ and $f(0)$ is due 
to the hyperbolic scaling of the velocity and forcing in the fluid. A parabolic
scaling in the fluid, which may appear to be more appropriate, would not,
however, lead to the necessary estimates, except for the case in which  
the initial solid-fluid interface is flat (which is  not the case considered 
herein). This is due to an elastic energy integral (which we shall shortly
identify) that requires the hyperbolic  scaling in order to be estimated.
\end{remark}

\begin{remark}  
Note also the presence of two compatibility conditions for the stresses on 
$\Gamma_0$, which is also a consequence of the hyperbolic scaling. A 
fluid-fluid interface problem would require only one compatibility condition.
\end{remark}

\begin{remark}
Note that the proof of {\it existence} of solutions  
requires only the ``minimal'' regularity assumptions 
(\ref{f_regularity}) on the forcing function $f$;
this is due to our method
of proof which employs (just as in \cite{CoSh2003} for the 
Navier-Stokes equations with surface tension) the Tychonoff fixed-point theorem
instead of a Banach-type contraction mapping. Note also that unlike the case 
of a free-surface fluid problem, a Banach contraction method does not work for 
the problem that we study herein.  We  will see later that some additional 
Lipschitz assumptions (\ref{Lip}) are necessary for uniqueness.
\end{remark}

\begin{remark}
We also remark that our technique is restricted to the case where the elastic 
constitutive law in the solid is either linear or semi-linear. Whereas the 
paper is written with a linear elasticity law, we can handle in the 
same fashion and with the exact same methods, the case where an extra 
contribution of the type $F(\nabla\eta,\eta)$ is added, with $F$ satisfying 
the usual regularity and growth assumptions. In that case the linear problem 
(\ref{linear}), defined hereafter, which is used in the fixed-point approach 
would be replaced by a similar problem, with (\ref{linear.d}) replaced by 
$$ w^i_t - [c^{ijkl}\int_0^t w^k,_l],_j + F(\nabla\int_0^t w,\int_0^t w)=  f^i
\ \text{in} \ \ (0,T)\times \Omega_0^s \ ,$$
which does not create any additional difficulties with respect to the analysis
of the linear case. 

The consideration of a quasilinear elastic law such as the nonlinear 
Saint-Venant Kirchhoff material, involves a smoother functional framework 
and will be developed in a future article. 
\end{remark}

\section{A bounded convex closed set of $W_T$}
\label{4}
\begin{definition}
Let $M>0$ be given. We let $C_T(M)$ denote the subset of $W_T$ 
consisting of elements $u\in W_T$ such that
\begin{align}
\|u\|_{W_T}^2\le M ,
\end{align}
and such that
\begin{equation}
\label{convex2}
u(0)|_{\Omega_0^f}=u_0|_{\Omega_0^f},\ \ \text{and}\ u_t(0)|_{\Omega_0^f}=w_1|_{\Omega_0^f},
\end{equation}
with $w_1$ defined in Theorem \ref{main}, and where we continue to assume that
the conditions stated in Theorem \ref{main} for the forcing function $f$ and 
the initial data $u_0$ are satisfied.
\end{definition}

\begin{lemma}
\label{convexnonempty}
There exists $M_0>0$ such that $C_T(M)$ is non-empty for $M>M_0$.  Furthermore,
$C_T (M)$ is a convex, bounded and closed subset of $X_T$. 
\end{lemma}
\begin{proof}
We note that if $\check v(t)=u_0+t w_1$, then  $\check v\in C_M(T)$ 
for $M\ge M_0= \|\check v\|^2_{W_T}$.
The fact that $C_T(M)$ is closed follows from Mazur's lemma.
\end{proof}

\begin{remark}
Note also that if $0<T'\le T$, then $C_{T'}(M)$ is non-empty. Henceforth,  
$M$ is assumed to be  larger that $M_0$. 
\end{remark}

In the remainder of the paper, we will assume that
$$0< T < \bar T$$ 
where the forcing $f$
is defined on the time interval $[0,\bar T]$; we  will have to 
choose $T$ sufficiently small to ensure existence of solutions to our problem.

We will need the following series of simple lemmas on the set $C_T (M)$.
\begin{lemma}
There exists $T_0\in (0,\bar{T})$ such that for all $T\in (0,T_0)$ and 
for all $v\in C_T(M)$,
$$\text{the matrix } a \text{ is well defined},$$
and satisfies the estimate (which is independent of $v\in C_T (M)$)
\begin{align}
\label{c0}
&\|a\|_{L^{\infty}(0,T;H^2(\Omega_0^f;{\mathbb R}^9))}
+\|a_t\|_{L^{\infty}(0,T;H^1(\Omega_0^f;{\mathbb R}^9))}
+\|a_{tt}\|_{L^{\infty}(0,T;L^2(\Omega_0^f;{\mathbb R}^9))}\nonumber\\
& +  \|a_{t}\|_{L^{2}(0,T;H^2(\Omega_0^f;{\mathbb R}^9))}
+ \|a_{tt}\|_{L^{2}(0,T;H^1(\Omega_0^f;{\mathbb R}^9))}
+ \|a_{ttt}\|_{L^{2}(0,T;L^2(\Omega_0^f;{\mathbb R}^9))}\le C(M)\ .
\end{align}
\end{lemma}
\begin{proof}
Notice that in the separable Hilbert space $H^3(\Omega_0^f;{\mathbb R}^3)$ 
(for which the Bochner integral is well defined),
$$\eta(t)=\text{Id}+\int_0^t v(s)ds\ ;$$
this together with the
Jensen and Cauchy-Schwarz inequalities shows that
$$\|\eta-\text{Id}\|_{L^{\infty}(0,T;H^3(\Omega_0^f;{\mathbb R}^3))}\le
C\ \sqrt{ T}\ \| v\|_{L^2(0,T;H^3(\Omega_0^f;{\mathbb R}^3))} ,$$
and thus
\begin{equation}
\label{c1}
\|\nabla\eta-\text{I}\|_{L^{\infty}(0,T;H^2(\Omega_0^f;{\mathbb R}^9))}
\le C\ \sqrt{T}\ \| v\|_{X_T}\le C\ \sqrt{ T}\ \sqrt{M}  .
\end{equation}
Next, choose  $R>0$ be such that 
for any $3\times 3$ matrix $b$ satisfying $\|b-\text{I}\|_{{\mathbb R}^9}\le R$,
we have $\operatorname{det}b\ge\frac{1}{2}$. 

We then see from (\ref{c1}) and the Sobolev inequalities that for 
$T\le T_0=\frac{CR^2}{{ M}}$, $\nabla\eta(t)$ is invertible for 
$t\in [0,T]$ in $\Omega_0^f$ for any $v\in C_T (M)$. From now on, 
$T$ is assumed 
to be in $(0,T_0)$. Since $$a(t)=\frac{1}{\operatorname{det}\nabla\eta(t)}\ 
\text{Cof}\nabla\eta(t)\ \text{in}\ \Omega_0^f\ ,$$
we then see from (\ref{c1}) that
\begin{equation*}
\|a\|_{L^{\infty}(0,T;H^2(\Omega_0^f;{\mathbb R}^9))}\le C(1+\sqrt{T M})^5\ .
\end{equation*}
Similarly,
\begin{equation}
\label{c2}
\|v-u_0\|_{L^{\infty}(0,T;H^2(\Omega_0^f;{\mathbb R}^3))}
\le C\ \sqrt{ T M}  ,
\end{equation}
providing
\begin{equation*}
\|a_t\|_{L^{\infty}(0,T;H^1(\Omega_0^f;{\mathbb R}^9))}\le C(1+\|u_0\|_{H^3(\Omega_0^f;{\mathbb R}^3)}+\sqrt{TM})^5\ .
\end{equation*}
In the same fashion,
\begin{equation}
\label{c3}
\|v_t-w_1\|_{L^{\infty}(0,T;H^1(\Omega_0^f;{\mathbb R}^3))}
\le C\ \sqrt{ T}\ M  ,
\end{equation}
providing
\begin{equation*}
\|a_{tt}\|_{L^{\infty}(0,T;H^1(\Omega_0^f;{\mathbb R}^9))}\le C(1+\|w_1\|_{H^1(\Omega_0^f;{\mathbb R}^3)}+\sqrt{T}M)^5\ .
\end{equation*}
The $L^2$ in time estimates are established in a straightforward manner from the
definition of $C_T(M)$, which concludes the proof of the lemma.
\end{proof}

\begin{remark}
Note that $T_0$ also depends on $M$.
\end{remark}

In the following, $T$ is taken in $(0,T_0)$ (and $M$ is still taken in 
$(0,M_0)$).
By the same arguments as above, we can easily prove the following results:
\begin{lemma}
For all $v\in C_M(T)$, we have
\begin{equation}
\label{c6}
\|a-a (0)\|^2_{L^{\infty}(0,T;H^2(\Omega_0^f;{\mathbb R}^9))}+\|a_{t}-a_t (0)\|^2_{L^{\infty}(0,T;H^1(\Omega_0^f;{\mathbb R}^9))}\le C(M)\ T\ .
\end{equation}
\end{lemma}

\begin{lemma}
There exists $T_1\in (0,T_0)$ which depends on $M$, and a constant $C>0$ which 
depends on
$u_0$ but does not depend on $M$, such that we have for all $v\in C_M(T)$,
\begin{equation}
\label{etah3}
\|\eta\|^2_{L^{\infty}(0,T;H^3(\Omega_0^f;{\mathbb R}^3))}
+ \|v\|^2_{L^{\infty}(0,T;H^2(\Omega_0^f;{\mathbb R}^3))}
+ \|v_t\|^2_{L^{\infty}(0,T;H^1(\Omega_0^f;{\mathbb R}^3))}
\le C \ .
\end{equation}
\end{lemma}

The next result concerns potential solid-solid or solid-container collisions
for a short time.

\begin{lemma} 
\label{collision}
Let $d>0$ denote the infimum of the distances between two distinct
 connected components of $\Omega_0^s$ (if we have more than one solid in the problem) and of the distance between $\Omega_0^s$ and $\partial\Omega$. Then, there exists $T_2\in (0,T_0)$ such that for all $v\in C_M(T_2)$, 
\begin{equation}
\label{collisionbis}
\displaystyle\int_0^{T_2} \|v\|_{L^{\infty}(\Omega_0^f; {\mathbb R}^3)}\le \frac{d}{2}. 
\end{equation}
\end{lemma}
\begin{proof}
The inequality $\displaystyle \int_0^T \|v\|_{L^{\infty}
(\Omega_0^f; {\mathbb R}^3)}\le 
C \sqrt{T}\ [\int_0^T \|v\|^2_{H^2(\Omega_0^f; {\mathbb R}^3)}]^{\frac{1}{2}}$ 
proves the result.
\end{proof}

Henceforth, we shall require
$$T\in (0,T_M), \ \ \  T_M=\min (T_1,T_2).$$ 
The series of estimates in Sections \ref{8} and \ref{9} will show that $M$ 
must first be chosen sufficiently large, and then $T$ 
must be chosen sufficiently small.

The next result is crucial for the derivation of appropriate estimates; while it
appears that we should require an
estimate of $q_{tt}$ in $L^2(0,T;L^2(\Omega_0^f; {\mathbb R}))$, we are not
able to obtain such an estimate, and effectively replace it with an estimate
of $q_t$ in $L^\infty(0,T;L^2(\Omega_0^f; {\mathbb R}))$.

\begin{lemma}
For all $v\in C_M(T)$,
\begin{equation}
\label{c4}
\|a_{tt}(t)\|_{L^{\infty}(0,T;L^3(\Omega_0^f;{\mathbb R}^9))}\le C(M)\ .
\end{equation}
\end{lemma}
\begin{proof}
%Let $\rho\in\mathcal D({\mathbb R}^3)$ be such that $\rho\ge 0$ and $\int_{{\mathbb R}^3} \rho=1$. Let then $\rho_n (x)=\frac{1}{n^3} \rho(\frac{x}{n})$ be the classical mollifying sequence associated and Let
Let $\displaystyle\psi(t)=\int_{\Omega_0^f}|a_{tt}(t)|^3+1\ge 1$. 
We then have in the distributional sense 
$$ \displaystyle\psi'(t)=3\int_{\Omega_0^f}| a_{tt}(t)|^2  a_{ttt}(t)\ .$$ 
Thus,
$$\psi'(t)\le C\ \|a_{tt}(t)\|^2_{L^4(\Omega_0^f;{\mathbb R}^9)} 
\|a_{ttt}(t)\|_{L^2(\Omega_0^f;{\mathbb R}^9)}\ ,$$
which by interpolation yields
$$\psi'(t)\le C\ \| a_{tt}(t)\|_{L^3(\Omega_0^f;{\mathbb R}^9)} 
\|a_{tt}(t)\|_{H^1 (\Omega_0^f;{\mathbb R}^9)} \| 
a_{ttt}(t)\|_{L^2(\Omega_0^f;{\mathbb R}^9)}\ ,$$
{\it i.e}
$$\psi'(t)\le C\ [\psi(t)]^{\frac{1}{3}}\  \|a_{tt}(t)\|_{H^1 (\Omega_0^f;{\mathbb R}^9)} \|a_{ttt}(t)\|_{L^2(\Omega_0^f;{\mathbb R}^9)}\ .$$
Thus, since $\psi(t)\ge 1$,
$$\psi(t)\le [\psi(0)^{\frac{2}{3}}+C \int_0^t \|a_{tt}(t)\|_{H^1 (\Omega_0^f;{\mathbb R}^9)} \|a_{ttt}(t)\|_{L^2(\Omega_0^f;{\mathbb R}^9)}]^{\frac{3}{2}}\ ,$$
which by (\ref{c0}) provides
$$\psi(t)\le [\psi(0)^{\frac{2}{3}}+C(M)]^{\frac{3}{2}}\ ,$$
which establishes (\ref{c4}).
\end{proof}

\begin{remark}
Note that in the above $L^3$ estimate, the exponent $3$ is the limiting case 
for this lemma.
\end{remark}

\begin{remark}
Had we  not made the notational (constitutive)  simplification of Section 
\ref{1bis}, we would require the following Korn-type lemma in the 
Lagrangian setting (this is the only mathematical issue that the actual
constitutive law (\ref{fluidlaw}) requires):
\end{remark}
\begin{lemma}
There exists $T_3\in (0,T)$ such that for any $T\in (0,T_3)$ and 
$v\in C_T(M)$,
for all $\phi\in H^1_0(\Omega_0^f;{\mathbb R}^3)$ and $t \in [0,T]$ 
$$
\int_{\Omega_0^f} (a_j^k(t) \phi,_k^i+a_i^k(t)\phi,_k^j)(a_j^k(t) 
\phi,_k^i+a_i^k(t)\phi,_k^j)\ge C\ 
\|\phi\|^2_{H^1_0(\Omega_0^f;{\mathbb R}^3)}\ .$$
\end{lemma}
\begin{proof}
To prove this result, we let $a(t)=I + [a(t)-a(0)]$ and apply
(\ref{c6}) followed by the Korn inequality.
\end{proof}

\section{The basic linear problem}
\label{6}
Suppose that $M\ge M_0$, $T\in (0,T_M)$ and $ v\in C_T (M)$ are given. 
Let $\displaystyle\eta=Id+\int_0^{\cdot}v$ and let $a_i^j$ be the quantity
associated with $\eta$ through (\ref{a}).

We are concerned with the following time-dependent linear problem, whose 
fixed-point $w=v$ provides a solution to (\ref{nsl}):
\begin{subequations}
  \label{linear}
\begin{alignat}{2}
 w^i_t - \nu (a^j_l a^k_l  w^i,_k),_j + (a^k_i q),_k &=  f\circ{\eta} 
&&\text{in} \ \ (0,T)\times \Omega_0^f \,, 
         \label{linear.b}\\
   a^k_i  w^i,_k &= 0    &&\text{in} \ \ (0,T)\times \Omega_0^f \,, 
         \label{linear.c}\\
 w^i_t - [c^{ijkl}\int_0^t w^k,_l],_j &=  f^i
&&\text{in} \ \ (0,T)\times \Omega_0^s \,, 
         \label{linear.d}\\
\nu w^i,_k a^k_l a^j_l N_j - q a^j_i N_j &= 
c^{ijkl}\int_0^t  w^k,_l\ N_j
&&\text{on} \ \ (0,T)\times \Gamma_0 \,, 
         \label{linear.e}\\
w(t,\cdot)& \in H^1_0 (\Omega;{\mathbb R}^3) &&\text{a.e. in } \ \ (0,T)\,, \label{linear.f}\\
   w &= u_0  
 &&\text{on} \ \ \Omega_0\times \{ t=0\} \,, 
         \label{linear.g}\\
   \eta &= \text{Id}     
 &&\text{on} \ \ \Omega_0\times \{ t=0\} \,, 
         \label{linear.h}
\end{alignat}
\end{subequations}
 
The following regularity result will be of paramount importance in our analysis:

\begin{theorem}\label{thm1}
Given $f$ and $u_0$ satisfying the assumptions of Theorem \ref{main}, there exists $M>0$, $T>0$, such that for any $v\in C_T (M)$, there exists a unique solution
$( w,p) \in Z_T$ 
of (\ref{linear}). Furthermore, $w\in C_T (M)$. 

\end{theorem}

Sections \ref{8} and \ref{9} are 
devoted to the proof of this theorem. In the following, we set 
\begin{align}
N(u_0,f)^2=&(1 + \ \|u_0\|^2_{H^5(\Omega_0^f; {\mathbb R}^3)}+ 
\|u_0\|^2_{H^2(\Omega_0^s; {\mathbb R}^3)}
+\|f(0)\|^2_{H^3(\Omega; {\mathbb R}^3)}
+\|f\|^2_{L^2(0,T; H^2(\Omega; {\mathbb R}^3))}\nonumber\\
&+ \|f_t\|^2_{L^2(0,T; H^1(\Omega; {\mathbb R}^3))}
+\|f_{tt}\|^2_{L^2(0,T; L^2(\Omega; {\mathbb R}^3))})^4\ .
\end{align}

\medskip

\section{Preliminary results}
\label{7}
\subsection{Divergence, extension and regularization type results.}  
We first state the following result, whose proof follows the same argument as 
for the case of a smooth boundary, with the exception that the regularity 
results for elliptic systems of \cite{Eben} are used instead of the more 
classical results wherein  the boundary is smooth.
\begin{lemma}
\label{lift}
Let $\Omega'$ be a domain of class $H^k$ ($k\ge 3$).\hfill\break
Then, for $0\le m\le k-2$, there exists a continuous linear operator 
$$L(\Omega'):\{(d,r)\in H^m(\Omega';{\mathbb R})\times 
H^{m+0.5}(\partial\Omega';{\mathbb R}^3)|\ \displaystyle \int_{\Omega} 
d=\int_{\partial\Omega} r  \cdot n\} 
\rightarrow H^{m+1}(\Omega';{\mathbb R}^3)$$
such that $u=L(\Omega')(d,r)$ satisfies
\begin{align*}
\operatorname{div}u &=d\ \text{in}\ \Omega',\\
u&=r\ \text{on}\ \partial\Omega'\ .
\end{align*}
Furthermore, the operator norm of $L(\Omega')$ remains bounded as
the norm of the charts defining $\Omega'$ stays in a bounded set of $H^k$. 
\end{lemma}

We will need the following extension
\begin{lemma}
\label{extension}
Recalling that $\Omega_0^f$ is of class $H^3$,
for each $1\le m\le 3$, there exists a continuous  linear operator 
$$E: H^m(\Omega_0^f;{\mathbb R}^3)\cap 
H^1_{\partial\Omega}(\Omega_0^f;{\mathbb R}^3) 
\rightarrow
H^m(\Omega;{\mathbb R}^3)\cap H^1_0(\Omega;{\mathbb R}^3)$$
such that $E(u)=u$ in $\Omega_0^f$. . 
\end{lemma}

\begin{proof}
The result is well-known in
the case that $\Omega_0^f= {\mathbb R}^3_+$;  let $\Pi: {\mathbb R}^3_+
\rightarrow {\mathbb R}^3$ denote this extension operator,
and let $\{ \Psi_i\}_{i=1}^N$ denote a collection of charts
in a neighborhood of $\Gamma_0$ (each $\Psi_i$ is a map of class $H^3$ 
from the unit ball in ${\mathbb R}^3$ into an open set containing a 
coordinate patch of $\Gamma_0$), and let
$(\theta_i)_{i=1}^N$ denote the associated partition of the unity.
We see that 
$$F(u)=\sum_{i=1}^N \Pi[(\theta_i u)\circ\Psi_i]\circ\Psi_i^{-1}$$
is an extension of $u$ into a neighborhood of $\Gamma_0$. 
By introducing a smooth cut-off function $\xi$, equal to $1$ in $\Omega^f_0$ 
and equal to $0$ in the complementary part of this neighborhood included in 
$\Omega_0^s$, we see that $E(u)=\xi\ F(u)$ satisfies the statement of the lemma.
\end{proof}

In a similar fashion, we can also extend from $\Omega$ to ${\mathbb R}^3$, with the same arguments.

\begin{lemma}
\label{extension2}
 There exists a linear and continuous operator $E_g$ from $H^m(\Omega;{\mathbb R}^3)$ into $H^m({\mathbb R}^3;{\mathbb R}^3)$ (for each $1\le m\le 3$) such that $E_g(u)=u$ in $\Omega$. 
\end{lemma}

We also need a regularization lemma for the coefficients $a_i^j$, which we shall
use to obtain estimates for the 
solutions of the regularized problems (whose coefficients by definition use 
these  regularized coefficients); we will then pass to the limit as the 
regularization parameter tends to zero.

\begin{lemma}
\label{regularize}
Let $ v\in C_T(M)$ and $ \displaystyle \eta=Id+\int_0^{\cdot}  v$. Then, there
exists a sequence $v_n\in V^{reg}_f(T)=\{u\in  
L^2(0,T; H^3(\Omega_0^f;{\mathbb R}^3))|\  u_t\in  
L^2(0,T; H^3(\Omega_0^f;{\mathbb R}^3)),\ u_{tt}\in  
L^2(0,T; H^3(\Omega_0^f;{\mathbb R}^3))\}$, 
such that $v_n(0)|_{\Omega_0^f}=u_0|_{\Omega_0^f}$, 
${v_n}_t(0)|_{\Omega_0^f}=w_1|_{\Omega_0^f}$, and 
$$\|v_n-v\|_{V^3_f(T)}+\|(v_n-v)_{tt}\|_{L^{\infty}
(0,T;L^2(\Omega_0^f;{\mathbb R}^3))}\rightarrow 0.$$
\end{lemma}
\begin{proof}

Let $\rho\in\mathcal D(B(0,1))$ be such that $\rho\ge 0$ and $\displaystyle\int_{B(0,1)} \rho=1$, and let $\rho_n(x)=n^3 \rho({x}{n})$ denote the usual mollifier.

From
Lemma \ref{extension}, for any $t\in [0,T]$, let $\bar v(t)=E( v(t))$, so that 
$\bar v\in V_3^f(T)\cap V^s_3(T)$ with $\|\bar v\|_{V^3_f(T)}+ 
\|\bar v\|_{V^3_s(T)} \le C\ \|v\|_{V^3_f(T)}\ $. 
We extend to ${\mathbb R}^3$ by  setting $v'=E_g(\bar v)$.

Then let $\tilde v_n$ be defined for any $t\in [0,T]$ by 
\begin{align*}
 \tilde v_n(t)&=\rho_n\star v'(t)  .
\end{align*}
From the properties of the space convolution, we have that 
$v_n\in V^{reg}_f (T)\cap V^3_s(T)$ and that 
$\|\tilde v_n- v'\|_{V^3_f(T)}+\|\tilde v_n- v'\|_{V^3_s(T)}+ 
\|(v_n-v')_{tt}\|_{L^{\infty}(0,T;L^2(\Omega_0^f;{\mathbb R}^3))}\rightarrow 0$ 
as $n\rightarrow\infty$. This in turn implies that 
$\|\tilde v_n-v\|_{V^3_f(T)}+ \|(v_n-v)_{tt}\|
_{L^{\infty}(0,T;L^2(\Omega_0^f;{\mathbb R}^3))}\rightarrow 0$ 
as $n\rightarrow\infty$. Now, for the initial conditions, 
let us define 
$$v_n(t)=u_0+t w_1-\int_0^t (t'-t)(\tilde v_n)_{tt}\ dt'\ .$$ 
(the Bochner integral being 
well-defined in the Hilbert space $H^3(\Omega_0^f;{\mathbb R}^3)$).
We then have that $v_n(0)=u_0$, ${v_n}_t(0)=w_1$.

Moreover, $$v_n=\tilde v_n +E_g(E(u_0))-\rho_n\star E_g (E(u_0))+ t\ [E_g (E(w_1))-\rho_n\star E_g(E(w_1))], $$ 
which yields 
 $\| v_n-v\|_{V^3_f(T)}+  \|(v_n-v)_{tt}\|_{L^{\infty}(0,T;L^2(\Omega_0^f;{\mathbb R}^3))}\rightarrow 0$, as $n\rightarrow\infty$. 
\end{proof}

\begin{remark}
Our construction does not necessarily yield
$v_n=0$ on $\partial\Omega$. Consequently, with
$\displaystyle\eta_n=\text{Id}+\int_0^{\cdot} v_n$ we do not have that 
$\eta_n (\Omega)=\Omega$. It, nevertheless, does not matter for the purpose of 
our analysis.
\end{remark}

\begin{remark}
\label{smoothcoefficients}
In the following, we will solve (\ref{linear}) as the limit as $n\rightarrow\infty$ of the solutions $w_n$ to the problems (\ref{linear}) associated to these regularized $v_n$. 
The interest of this regularizing process is that for a given $n$, $a(\eta_n)$
and its first and second time derivatives are in 
$L^{\infty}(0,T; H^2(\Omega_0^f;{\mathbb R}^9))\subset L^{\infty}
(0,T; L^{\infty}(\Omega_0^f;{\mathbb R}^9))$ and its third time derivative is 
in $L^{2}(0,T; H^2(\Omega_0^f;{\mathbb R}^9))\subset L^{2}
(0,T; L^{\infty}(\Omega_0^f;{\mathbb R}^9))$ which is necessary in order to 
get the existence of regular solutions to (\ref{linear}). These bounds in 
those spaces of course blow up as $n\rightarrow \infty$ (except for the 
estimate for $a(\eta_n)$). 

Nevertheless,  using the fact that $\|v_n-v\|_{V^3_f (T)}\rightarrow 0$ as
$n\rightarrow\infty$, 
$a(\eta_n)$, and its first, second and third time derivatives satisfy
the same type of estimates as (\ref{c0}), 
(\ref{c4}) and (\ref{c6}), respectively, with a constant $C(M)$ which does not 
depend on $n$. This fact will be used, together with interpolation inequalities 
(that hold since the solutions $w_n$ are regular) in order to get estimates in
$Y_T$ for $w_n$ which are independent of $n$. By weak convergence, this will 
provide our smooth solution to (\ref{linear}).

We will also use  the convention of denoting the regularized velocity
fields  $v_n$ by $\tilde v$, and
the corresponding regularized matrix $a(\eta_n)$ by $\tilde a$. 
\end{remark}

\begin{remark}
Since the fluid forcing in (\ref{linear}) is given $f\circ\eta$, we need to 
extend $f$ to ${\mathbb R}^3$. Hence, when we solve this problem with the
regularized coefficients arising from $\tilde v$, we in fact implicitly use 
the extension $E_g (f)$.  This extension has the same regularity as $f$ 
with ${\mathbb R}^3$ replacing $\Omega$; this follows from the fact that
$E_g$ commutes with the time derivative. For notational convenience, we shall
continue to denote the extended forcing function by the same letter $f$.
\end{remark}

\subsection{Pressure as a Lagrange multiplier}

\begin{lemma}\label{lemma_lagrange}
For all $p \in L^2(\Omega_0^f; {\mathbb R})$, $t\in [0,T]$, there exists a constant
$C>0$ and
$\phi \in H^1_0(\Omega, {\mathbb R}^3)$ such that $a_i^j (t) \phi^i,_j =p$ in $\Omega_0^f$ and
\begin{equation}\label{v-p} 
\|\phi\|^2_{H^1_0(\Omega; {\mathbb R}^3)} \le  C\|p\|^2_{L^2(\Omega_0^f; {\mathbb R})}. 
\end{equation}
\end{lemma}
\begin{proof}Let $p_1\in L^2 (\Omega_0^s;\mathbb R)$ be such that $p_1 >0$ in $\Omega_0^s$. Let $\bar p$ be defined by $\bar{p}=p$ in $\Omega_0^f$ and
$\displaystyle\bar{p}=-{\frac{\int_{\Omega_0^f} p\ \text{det}\nabla\eta(t)}{\int_{\Omega_0^s} p_1 \ \text{det}\nabla\eta(t)}} p_1$ in $\Omega_0^s$.
Since $$\int_{\eta(t,\Omega)}\bar p\circ{\eta(t)}^{-1}\ dx =\int_{\Omega}\bar p(x)\ \ \text{det}\nabla\eta(t)dx =0\ ,$$ we then see that 
$\phi=L(\eta(t,\Omega))(\bar p\circ{\eta(t)}^{-1},0)\circ{\eta(t)}
\in H^1_0(\Omega;{\mathbb R}^3)$ satisfies 
$$\operatorname{div} (\phi\circ{\eta(t)}^{-1})=\bar p\circ{\eta(t)}^{-1}\ \text{in}\ \eta(t,\Omega)\ ,$$
and thus
$$a_i^j(t) \phi^i,_j=\operatorname{div} (\phi\circ{\eta(t)}^{-1})\circ\eta(t)=\bar p=p\ \text{in}\ \Omega_0^f\ .$$

The inequality (\ref{v-p}) is then a simple consequence of the properties of
$L$ and of the condition $v\in C_T (M)$.
\end{proof}

We can now follow \cite{SolSca1973}.
Define the linear functional on $H^1_0(\Omega; {\mathbb R}^3)$
by $(p,a_i^j (t) \varphi^i,_j)_ {L^2(\Omega_0^f; {\mathbb R})}$, where $\varphi\in H^1_0(\Omega; {\mathbb R}^3)$.
By the Riesz representation theorem, there is a bounded linear operator
$Q (t): L^2(\Omega_0^f;{\mathbb R}^3)\rightarrow  H^1_0(\Omega; {\mathbb R}^3)$ such that
$$ 
\forall \varphi\in H^1_0(\Omega; {\mathbb R}^3),\ (p,\ a_i^j (t) \varphi^i,_j)_ {L^2(\Omega_0^f; {\mathbb R})}=
(Q(t)p,\ \varphi)_{H^1_0(\Omega; {\mathbb R}^3)}. 
$$ 
Letting $\varphi=Q(t)p$ shows that 
\begin{equation}\label{Qp0}
\|Q(t)p\|_{H^1_0(\Omega; {\mathbb R}^3)} \le C
\|p \|_ {L^2(\Omega_0^f; {\mathbb R})}
\end{equation}
for some constant $C>0$. 
Using Lemma \ref{lemma_lagrange},  we have the estimate
\begin{equation}\label{Qp}
\|p\|^2_{L^2(\Omega_0^f; {\mathbb R})}\le C\|Q(t)p\|_{H^1_0(\Omega; {\mathbb R}^3)}\|\phi\|_{H^1_0(\Omega; {\mathbb R}^3)} 
\le C\|Q(t)p\|_{H^1_0(\Omega; {\mathbb R}^3)}\|p\|_{L^2(\Omega_0^f; {\mathbb R})}, 
\end{equation}
which shows that $R(Q(t))$ is closed in $H^1_0(\Omega; {\mathbb R}^3)$.
Since ${\mathcal V}_v(t) \subset R(Q(t))^\perp$ and $R(Q(t))^\perp \subset {\mathcal V}_v (t)$, it follows that 
\begin{equation}\label{hodge}
H^1_0(\Omega; {\mathbb R}^3) = R(Q(t)) \oplus_ {H^1_0(\Omega; {\mathbb R}^3)} {\mathcal V}_v(t).
\end{equation}

We can now introduce our Lagrange multiplier 
\begin{lemma} \label{Lagrange}
Let ${\mathfrak L}(t)\in H^{-1} (\Omega;{\mathbb R}^3)$ be such that ${\mathfrak L}(t) \varphi=0$ for any
$\varphi\in {\mathcal V}_v(t)$. Then there exists a unique $q(t)\in L^2(\Omega_0^f; {\mathbb R})$, which is termed the pressure function, satisfying 
$$\forall \varphi\in H^{1}_0 (\Omega;{\mathbb R}^3),\ \ 
{\mathfrak L}(t) (\varphi)=(q(t),\ a_i^j \varphi^i,_j)_{L^2(\Omega_0^f; {\mathbb R})}.$$ 
Moreover, there is a $C>0$ (which does not depend on $t\in [0,T]$ and on the choice of $v\in C_M (T)$) such that$$\|q(t)\|_{L^2(\Omega_0^f; {\mathbb R})}\le C\ \|{\mathfrak L}(t)\|_{H^{-1}(\Omega; {\mathbb R}^3)}.$$
\end{lemma}

\begin{proof}
By the decomposition (\ref{hodge}), for $v\in{H^1_0(\Omega, {\mathbb R}^3)}$, we let
$\varphi=v_1+v_2$, where $v_1 \in {\mathcal V}_v (t)$ and $v_2 \in R(Q(t))$.  From our assumption, it follows that
$$ 
{\mathfrak L}(t)(\varphi) = {\mathfrak L}(t)(v_2) = ( \psi(t), v_2)_{H^1_0(\Omega, {\mathbb R}^3)} 
= ( \psi(t), \varphi)_{H^1_0(\Omega, {\mathbb R}^3)}, $$ \ for a unique  \ $\psi(t) \in R(Q(t))$.

From the definition of $Q(t)$ we then get the existence of a unique 
$q(t)\in L^2(\Omega_0; {\mathbb R})$ such that 
$$\forall \varphi\in H^{1}_0 (\Omega;{\mathbb R}^3),\ \ {\mathfrak L}(t) (\varphi)
=(q(t),\ a_i^j \varphi^i,_j)_{L^2(\Omega_0^f; {\mathbb R})}.$$
The estimate stated in the lemma is then a simple consequence of (\ref{Qp}).
\end{proof}

We will also need a version of the Lagrange multiplier lemma for the case 
where ${\mathfrak L}(t)\in H^{-1}(\Omega_0^f;{\mathbb R}^3)$, which implies an 
estimate on the pressure modulo a constant. We first have

\begin{lemma}\label{lemma_lagrangebis}
For all $p \in L^2(\Omega_0^f; {\mathbb R})$ such that $\displaystyle\int_{\Omega_0^f} p\ \text{det}\ \nabla\eta=0$, $t\in [0,T]$, there exists a constant
$C>0$ and
$\phi \in H^1_0(\Omega_0^f; {\mathbb R}^3)$ such that $a_i^j (t) \phi^i,_j =p$ in $\Omega_0^f$ and
\begin{equation}\label{v-pbis} 
\|\phi\|^2_{H^1_0(\Omega^f_0; {\mathbb R}^3)} \le  C\|p\|^2_{L^2(\Omega_0^f; {\mathbb R})}. 
\end{equation}
\end{lemma}
\begin{proof}
Since $$\int_{\eta(t,\Omega_0^f)}\bar p\circ{\eta(t)}^{-1}\ dx 
=\int_{\Omega_0^f}\bar p(x)\ \ \text{det}\nabla\eta(t)dx =0\ ,$$ 
we then define $\phi=L({\eta(t,\Omega_0^f)})(\bar p\circ{\eta(t)}^{-1},0)
\circ{\eta(t)}\in H^1_0(\Omega_0^f;{\mathbb R}^3)$. 

The inequality (\ref{v-p}) is then a simple consequence of the properties of
$L$ and of the condition $v\in C_T (M)$.
\end{proof}

In a completely similar fashion to the proof of Lemma \ref{Lagrange}, we can now
establish our second Lagrange multiplier 
\begin{lemma} \label{Lagrangebis}
Let ${\mathfrak L}(t)\in H^{-1} (\Omega_0^f;{\mathbb R}^3)$ be such that 
${\mathfrak L}(t) \varphi=0$ for any
$\varphi\in {\mathcal V}_v(t)\cap H^1_0(\Omega_0^f;{\mathbb R}^3)$. 
Then there exists a unique $q(t)\in L^2(\Omega_0^f; {\mathbb R})$, 
satisfying $\displaystyle \int_{\Omega_0^f} q(t) \text{det}\ \nabla\eta=0$, 
which is termed the pressure function, satisfying 
$$\forall \varphi\in H^{1}_0 (\Omega_0^f;{\mathbb R}^3),\ \ \mathfrak{L}(t) (\varphi)=(q(t),
\ a_i^j \varphi^i,_j)_{L^2(\Omega_0^f; {\mathbb R})}.$$ 
Moreover, there is a $C>0$ (which does not depend on $t\in [0,T]$ and on the 
choice of $v\in C_M (T)$) such that
$$\|q(t)\|_{L^2(\Omega_0^f; {\mathbb R})}\le 
C\ \|\mathfrak{L}(t)\|_{H^{-1}(\Omega_0^f; {\mathbb R}^3)}.$$
\end{lemma}

\begin{remark}
The four previous lemmas do not rely on the fact that $v=0$ on $\partial\Omega$. 
Therefore, they are also true for the case where the coefficients $\tilde a$ are
associated to $\tilde v$. The important point is that the estimates (\ref{c0}), 
(\ref{c4}) and (\ref{c6}) are also satisfied by the regularized matrix $\tilde a$ 
and velocity $\tilde v$.
\end{remark}

\section{Estimates for (\ref{linear}): the case of the regularized coefficients}
\label{8}

\subsection{Weak solutions}
\begin{definition} 
A vector $w \in {\mathcal V}_v([0,T])$ with $w_t \in
{\mathcal V}_v(t)'$ for a.e. $t\in (0,T)$ is a weak solution of (\ref{linear})
provided that a.e. $t\in (0,T)$, 
\begin{subequations}
\label{weakw}
\begin{alignat}{2}
\operatorname{(i)} \ & \langle w_t, \phi\rangle 
+\nu (a^r_k w^i,_r, a^s_k \phi^i,_s)_{L^2(\Omega_0^f;{\mathbb R}^9)} 
+ (c^{ijkl}\int_0^t  w^k,_l, \phi^i,_j)_{L^2(\Omega_0^s;{\mathbb R})} \label{weakw.a}\\ & \ =  ( F, \phi)_{L^2 (\Omega_0^f;{\mathbb R}^3)} 
+ (f,\phi)_{L^2 (\Omega_0^s;{\mathbb R}^3)}    , \ \forall \phi\in 
{\mathcal V}_v(t)\,,\ \operatorname{and} \nonumber\\
\operatorname{(ii)} \ & w(0,\cdot) = u_0 \label{weakw.b},
\end{alignat}
\end{subequations}
for a.e. $0\le t \le T$,
where $\langle\cdot,\cdot\rangle$ denotes the duality product between 
${\mathcal V}_v (t)$ and its dual ${\mathcal V}_v (t)'$.
\end{definition}

\subsection{Penalized problems}

Whereas the existence of a weak solution can be proved directly in the space
${\mathcal V}_v([0,T])$, with $w_t\in {\mathcal V}_v(t)'$, this framework
is not amenable to finding the pressure estimate required by our analysis.
(Even for the well-studied problem of the Navier-Stokes equations on a fixed
and smooth bounded domain, the weak solution only provides a pressure
estimate of the form $\int_0^\cdot p \in L^2(0,T; L^2(\Omega^f_0; {\mathbb R}))$.)
A penalized form of the problem, however, together with the penalized form for
the time-differentiated problem, provide the correct pressure estimate in the
limit as  the penalization parameter tends to zero.

As we noted following Lemma \ref{regularize}, we will work with a regularized
sequence of
velocities $v_n$, and we shall generically denote elements of this sequence
simply as $\tilde v$, and the associated regularized matrices
$a(\eta_n)$ as $\tilde a$.

Given the regularity assumptions in (\ref{f_regularity}), 
$f_t\in C([0,T];L^2 (\Omega;{\mathbb R}^3)), \ \
%g_t\in C([0,T];H^{\frac{1}{2}}(\Gamma_0;{\mathbb R}^3)),
$ 
so that
$f_t(0)\in L^2(\Omega_0^f;{\mathbb R}^3))$. 

Then, let $w_2\in L^2(\Omega;{\mathbb R}^3)$ be defined by
\begin{subequations}
\label{w2def}
\begin{align}
w_2^i= &\ \nu \triangle  w_1^i +\nu (( a_l^j a_l^k)_t (0) u_0^i,_k),_j + F_t(0) - ((a_i^j)_t(0) q_0),_j - q_1,_i
\ \text{in}\ \Omega_0^f\ ,\label{w2def.a}\\
{w}_2^i =&\  f_t^i(0) +[c^{ijkl} u_0^k,_l],_j\ \text{in}\ \Omega_0^s\ ,
\end{align}
\end{subequations}
where $q_1\in H^1(\Omega_0^f;{\mathbb R})$ is defined by
\begin{subequations}
\label{q1def}
\begin{align}
\triangle q_1= &\ \frac{\partial}{\partial x^i}[\nu \triangle (w^i)_1 + 
(F^i)_t(0)+ 
\nu (( a_l^j a_l^k)_t (0) u_0^i,_k),_j - ((a_i^j)_t(0) q_0),_j]\nonumber\\
& + 2 (a_i^j)_t (0) w_1^i,_j+(a_i^j)_{tt} (0) u_0^i,_j
\ \text{in}\ \Omega_0^f\ ,\\
q_1 =&\  \nu [\nabla w_1\ N\cdot N + (a_l^k a_l^j)_t(0) u_0^i,_k N_j N_i] 
-c^{ijkl} u_0^k,_l N_j N_i\ \text{on}\ \Gamma_0 \nonumber \\
& + q_0(a^j_i)_t(0) N_jN_i \, \\
\frac{\partial q_1}{\partial N}=& F_t(0)\cdot N - [(a_i^j)_t(0) q_0],_j +\nu \triangle w_1\cdot N+\nu ((a_l^j a_l^k)_t(0) u_0^i,_k),_j N_i\ \text{on}\ \partial\Omega\ . 
\end{align}
\end{subequations}

Once again, we remind the reader  that $(a_i^j)_t(0)$ and $(a_i^j)_{tt}(0)$ depend
only on $u_0$ and $w_1$, and we note that they are equal to $(\tilde a_i^j)_t(0)$ 
and $(\tilde a_i^j)_{tt}(0)$, respectively.

Letting $\epsilon>0$ denote the penalization parameter,  we define 
$w_{\epsilon}\in {\mathcal W}([0,T])$ to be  the unique weak solution of the 
problem (whose existence can be obtained via
a standard Galerkin method in a basis of $H^1_0 (\Omega;{\mathbb R}^3)$):
\begin{subequations}
\label{weakpenalty}
\begin{align}
\operatorname{(i)} \ & \langle {w_{\epsilon}}_t,\ \phi\rangle 
+\nu (\tilde a^r_k w_{\epsilon}^i,_r,\ \tilde a^s_k \phi^i,_s)_{L^2(\Omega_0^f;{\mathbb R})} 
+ (c^{ijkl}\int_0^t  w_{\epsilon}^k,_l,\ \phi^i,_j)_{L^2(\Omega_0^s;{\mathbb R})}\nonumber \\ & +(\frac{1}{\epsilon} \tilde a^i_j w_{\epsilon}^j,_i-q_0-t q_1,\ \tilde a^l_k \phi^k,_l)_{L^2(\Omega_0^f;{\mathbb R})}=  (F,\ \phi)_{L^2 (\Omega_0^f;{\mathbb R}^3)} 
+ (f,\ \phi)_{L^2 (\Omega_0^s;{\mathbb R}^3)} , \\
& \forall \phi\in 
H^1_0 (\Omega;{\mathbb R}^3)\,,\ 
 \operatorname{and} \nonumber\\ 
\operatorname{(ii)} \ & w(0,\cdot) = u_0 ,
\end{align}
\end{subequations}
where $\langle\cdot,\cdot\rangle$ denotes the duality product between
$H^1_0 (\Omega;{\mathbb R}^3)$ and its dual.

\subsection{Weak solutions for the penalized problem} The aim of this section is 
to establish the existence of $w_{\epsilon}$, as well as the energy equalities
satisfied by $w_{\epsilon}$ and ${w_{\epsilon}}_t$, and the energy inequality
satisfied by ${w_{\epsilon}}_{tt}$. It turns out that the exposition is simplified if we first study the twice differentiated-in-time problem, that we introduce now.

\noindent{\bf Step 1. Galerkin sequence.}

By introducing a basis $(e_l)_{l=1}^{\infty}$ of $H^1_0 (\Omega;{\mathbb R}^3)$ 
and $L^2(\Omega;{\mathbb R}^3)$, and taking the approximation at rank 
$l\ge 2$ under the form $\displaystyle w_l (t,x)=\sum_{k=1}^l y_k (t)\ e_k (x)\ ,$ 
and satisfying on $[0,T]$ 
\begin{align*}
&\text{(i)}\  ( {w_l}_{ttt}, \phi)_{L^2(\Omega;{\mathbb R}^3)}
+{\nu} (\tilde a_k^r {w_l}_{tt},_r,\ \tilde a_k^s \phi,_s)_{L^2(\Omega_0^f;{\mathbb R}^3)} 
+ (c^{ijkl}   {w_l}_{t}^k,_l, \phi^i,_j)_{L^2(\Omega_0^s;{\mathbb R})}  \nonumber \\ &\qquad + {\nu} ((\tilde a_k^r \tilde a_k^s)_{tt} {w_l},_r,\  \phi,_s)_{L^2(\Omega_0^f;{\mathbb R}^3)}+ 2 {\nu} ((\tilde a_k^r \tilde a_k^s)_{t} {w_l}_t,_r,\  \phi,_s)_{L^2(\Omega_0^f;{\mathbb R}^3)}\nonumber\\
&\qquad - ((\tilde a_i^j q_l)_{tt},\ \phi^i,_j)_{L^2(\Omega_0^f;{\mathbb R})}\nonumber\\
&\qquad\qquad\qquad\qquad =  (F_{tt}, \phi)_{L^2 (\Omega_0^f;{\mathbb R}^3)}
+ (f_{tt},\phi)_{L^2 (\Omega_0^s;{\mathbb R}^3)}, \ \forall \phi\in 
span(e_1,...,e_l)\,,\\
&\text{(ii)}\  {w_l}_{tt} (0)=({w}_2)_l,\ {w_l}_t (0)=(w_1)_l,\ 
w_l (0)=(u_0)_l,\ \text{in}\ \Omega\ ,
\end{align*}
where $\displaystyle q_l=q_0+t\ q_1-\frac{1}{\epsilon} \tilde a_i^j w_l^i,_j$ 
and $(w_2)_l$ denotes the
$L^2(\Omega;{\mathbb R}^3)$ projection of $w_2$ onto $span(e_1,...,e_l)$, and
$(w_1)_l$ and $(u_0)_l$ denote the respective $H^1_0(\Omega;{\mathbb R}^3)$ 
projections of $w_1$ and $u_0$ on $span(e_1,...,e_l)$, we see that the 
Cauchy-Lipschitz theorem gives us the local well-posedness for $w_l$. The use of
the test function $(w_l)_{tt}$ in this system of ODEs (which is allowed as
it belongs to $span(e_1,...,e_l)$) gives us in turn the energy law 
\begin{align*}
&  \frac{1}{2} \frac{d}{dt} \|{w_l}_{tt}\|^2_{L^2(\Omega;{\mathbb R}^3)} 
+{\nu} (\tilde a_k^r {w_l}_{tt},_r,\ \tilde a_k^s {w_l}_{tt},_s)_{L^2(\Omega_0^f;{\mathbb R}^3)}\nonumber\\
& + \frac{1}{2} \frac{d}{dt} (c^{ijkl}   {w_l}_t^k,_l, {w_l}_t^i,_j)_{L^2(\Omega_0^s;{\mathbb R})} + {\nu} ((\tilde a_k^r \tilde a_k^s)_{tt} {w_l},_r,\  {w_l}_{tt},_s)_{L^2(\Omega_0^f;{\mathbb R}^3)}\nonumber\\
& + 2 {\nu} ((\tilde a_k^r \tilde a_k^s)_{t} {w_l}_t,_r,\  {w_l}_{tt},_s)_{L^2(\Omega_0^f;{\mathbb R}^3)} - ((\tilde a_i^j q_l)_{tt},\ {w_l}_{tt}^i,_j)_{L^2(\Omega_0^f;{\mathbb R})}\nonumber\\
&\qquad\qquad\qquad\qquad =  (F_{tt}, {w_l}_{tt})_{L^2 (\Omega_0^f;{\mathbb R}^3)} 
+ (f_{tt},{w_l}_{tt})_{L^2 (\Omega_0^s;{\mathbb R}^3)} \ .
\end{align*}

After transforming the term with $(q_l)_{tt}$ (since it involves 
${\nabla w_l}_{tt}$)
 and integrating this relation from $0$ to $t\in(0,T)$,
\begin{align}
&  \frac{1}{2}  \|{w_l}_{tt}\|^2_{L^2(\Omega;{\mathbb R}^3)} 
+{\nu} \int_0^t (\tilde a_k^r {w_l}_{tt},_r,\  \tilde a_k^s {w_l}_{tt},_s)_{L^2(\Omega_0^f;{\mathbb R}^3)} 
+ \frac{1}{2}  (c^{ijkl}   {w_l}_{t}^k,_l, {w_l}_{t}^i,_j)_{L^2(\Omega_0^s;{\mathbb R})}  \nonumber \\
 & + \epsilon \int_0^t\int_{\Omega_0^f} {q_l}^2_{tt}+\int_0^t\int_{\Omega_0^f} {q_l}_{tt}\ [ 2(\tilde a_i^j)_t 
 {w_l}_t^i,_j + (\tilde a_i^j)_{tt} 
 {w_l}^i,_j] -2 \int_0^t\int_{\Omega_0^f} (\tilde a_i^j)_t {q_l}_t{w_l}_{tt}^i,_j\nonumber \\
& - \int_0^t\int_{\Omega_0^f} (\tilde a_i^j)_{tt} {q_l}{w_l}_{tt}^i,_j + {\nu} \int_0^t ((\tilde a_k^r \tilde a_k^s)_{tt} {w_l},_r,\  {w_l}_{tt},_s)_{L^2(\Omega_0^f;{\mathbb R}^3)}\nonumber \\
& +2\ {\nu} \int_0^t ((\tilde a_k^r \tilde a_k^s)_{t} {{w_l}_t},_r,\  {w_l}_{tt},_s)_{L^2(\Omega_0^f;{\mathbb R}^3)}\nonumber\\
&\qquad\le C\ N(u_0,f)^2 +\int_0^t (F_{tt}, {w_l}_{tt})_{L^2 (\Omega_0^f;{\mathbb R}^3)}+\int_0^t (f_{tt},{w_l}_{tt})_{L^2 (\Omega_0^s;{\mathbb R}^3)} .
\label{twice1}
\end{align}

By noticing that, $\epsilon$ being fixed, the fourth term of the left-hand side
of this inequality involving the square of $(q_l)_{tt}$ acts as a viscous
energy term, and taking into account the $L^{\infty}(0,T;L^{\infty}(\Omega_0^f;{\mathbb R}))$ bound of each one of
the regularized coefficients $\tilde a_i^j$ and their first and second time
derivatives, we then get 
\begin{align*}
& {\frac{1}{2}}   \|{w_l}_{tt} (t)\|^2_{L^2(\Omega;{\mathbb R}^3)} 
+ \frac{\nu}{2}\int_0^t \|\nabla {w_l}_{tt}\|^2_{L^2(\Omega_0^f;{\mathbb R}^9)} 
+ \frac{1}{2}  (c^{ijkl}   {w_l}_{t}^k,_l (t), {w_l}_{t}^i,_j (t))_{L^2(\Omega_0^s;{\mathbb R})}  \nonumber \\
 & + \frac{\epsilon}{4} \int_0^t\int_{\Omega_0^f} {q_l}^2_{tt}-\tilde C_{\epsilon} [\int_0^t \int_0^{t'} \|\nabla {w_l}_{tt}\|^2_{L^2(\Omega_0^f;{\mathbb R}^9)}+N(u_0,f)^2 + \int_0^t \int_0^{t'} \|{q_l}_{tt}\|^2_{L^2(\Omega_0^f;{\mathbb R})}]
\nonumber\\
&\qquad\le C\ N(u_0,f)^2 +\int_0^t (F_{tt}, {w_l}_{tt})_{L^2 (\Omega_0^f;{\mathbb R}^3)}+\int_0^t (f_{tt},{w_l}_{tt})_{L^2 (\Omega_0^s;{\mathbb R}^3)} ,
\end{align*}
where $\tilde C_{\epsilon}$ depends on the regularizing parameter of $\tilde a$ and on $\epsilon$, but not on $l$.
By Gronwall's inequality, we then get an estimate on each of the integral terms
multiplying $\tilde C_{\epsilon}$ which in turn implies
\begin{align*}
   \|{w_l}_{tt} (t)\|^2_{L^2(\Omega;{\mathbb R}^3)} 
+ \int_0^t \|\nabla {w_l}_{tt}\|^2_{L^2(\Omega_0^f;{\mathbb R}^9)} 
&+ \frac{1}{2}  (c^{ijkl}   {w_l}_{t}^k,_l (t), {w_l}_{t}^i,_j (t))_{L^2(\Omega_0^s;{\mathbb R})}  \nonumber \\
 & + \epsilon \int_0^t\int_{\Omega_0^f} {q_l}^2_{tt}\le \tilde C_{\epsilon} \ N(u_0,f)^2  .
\end{align*}

\noindent{\bf Step 2. Weak solution $w_{\epsilon}$ of the penalized problem and 
its time differentiated problem.}

We can then infer that $w_l$ is defined on $[0,T]$, and that there is a subsequence, still denoted with
the subscript $l$, satisfying 
\begin{subequations}
\label{twice2}
\begin{align}
w_{l}&\rightharpoonup   w_{\epsilon}\ \ \text{in}\ L^2(0,T; H^1_0(\Omega;{\mathbb R}^3))
\\
{w_l}_t &\rightharpoonup   {w_{\epsilon}}_t\ \ \text{in}\ L^2(0,T; H^1_0(\Omega;{\mathbb R}^3))
\\ 
{w_l}_{tt} &\rightharpoonup  {w_{\epsilon}}_{tt}\ \ \text{in}\ L^2(0,T; L^2(\Omega;{\mathbb R}^3))\ \text{and in}\ L^2(0,T; H^1(\Omega_0^f;{\mathbb R}^3))\\
{q_l}_{tt}&\rightharpoonup {q_{\epsilon}}_{tt}\ \ \text{in}\ L^2(0,T; L^2(\Omega_0^f;{\mathbb R}^3))\ ,
\end{align}
\end{subequations}
where 
\begin{equation}
q_{\epsilon}=q_0+t q_1-\frac{1}{\epsilon} \tilde a_i^j {w_{\epsilon}}^i,_j\ .
\end{equation}
From the standard procedure for weak solutions, we can now infer from these weak convergences and
the definition of $w_l$ that ${w_{\epsilon}}_{ttt}\in L^2(0,T;H^{-1}(\Omega;{\mathbb R}^3))$. In turn, ${w_{\epsilon}}_{tt}\in C^0([0,T]; H^{-1}(\Omega;{\mathbb R}^3))$, ${w_{\epsilon}}_{t}\in C^0([0,T]; L^2(\Omega;{\mathbb R}^3))$,
${w_{\epsilon}}\in C^0([0,T]; H_0^1(\Omega;{\mathbb R}^3))$, with
$w_{\epsilon}(0)=u_0$, ${w_{\epsilon}}_t(0)=w_1$, ${w_{\epsilon}}_{tt}(0)=w_2$.

We moreover have for ${w_l}_{t}$
\begin{align*}
&\text{(i)}\  ( {w_l}_{tt}, \phi)_{L^2(\Omega;{\mathbb R}^3)}
+{\nu} (\tilde a_k^r {w_l}_{t},_r,\ \tilde a_k^s \phi,_s)_{L^2(\Omega_0^f;{\mathbb R}^3)} 
+ (c^{ijkl}   {w_l}^k,_l, \phi^i,_j)_{L^2(\Omega_0^s;{\mathbb R})}  \nonumber \\ &\qquad + {\nu} ((\tilde a_k^r \tilde a_k^s)_{t} {w_l},_r,\  \phi,_s)_{L^2(\Omega_0^f;{\mathbb R}^3)}- ((\tilde a_i^j q_l)_{t},\ \phi^i,_j)_{L^2(\Omega_0^f;{\mathbb R})}\nonumber\\
&\qquad =  (F_{t}, \phi)_{L^2 (\Omega_0^f;{\mathbb R}^3)}
+ (f_{t},\phi)_{L^2 (\Omega_0^s;{\mathbb R}^3)} + c_l(\phi), \ \forall \phi\in 
span(e_1,...,e_l)\,,\\
&\text{(ii)}\ \ {w_l}_t (0)=(w_1)_l,\ 
w_l (0)=(u_0)_l,\ \text{in}\ \Omega\ ,
\end{align*}
where $c_l (\phi)\in \mathbb R$ is given by
\begin{align*}
c_l(\phi)=&\ ( (w_2)_{l}, \phi)_{L^2(\Omega;{\mathbb R}^3)}
+{\nu} (\tilde a_k^r(0) {w_1}_{l},_r,\ \tilde a_k^s(0) \phi,_s)_{L^2(\Omega_0^f;{\mathbb R}^3)} 
+ (c^{ijkl}   (w_1)_l ^k,_l, \phi^i,_j)_{L^2(\Omega_0^s;{\mathbb R})}  \nonumber \\ &+ {\nu} ((\tilde a_k^r \tilde a_k^s)_{t} (0) {(w_1)_l},_r,\  \phi,_s)_{L^2(\Omega_0^f;{\mathbb R}^3)} - ((\tilde a_i^j q_l)_{t} (0),\ \phi^i,_j)_{L^2(\Omega_0^f;{\mathbb R})}\nonumber\\
&- (F_{t}(0), \phi)_{L^2 (\Omega_0^f;{\mathbb R}^3)}
- (f_{t} (0),\phi)_{L^2 (\Omega_0^s;{\mathbb R}^3)}\ . 
\end{align*}
Thus, $c_l (\phi)$ converges to the same expression, where the approximate
initial data $(w_i)_l$ are replaced by the actual initial data $w_i$ ($i=0,1,2$).
From our compatibility conditions (\ref{compatibility}) together with
(\ref{w2def}), this leads us to
\begin{equation}
\label{cl}
\|c_l\|_{H^{-1} (\Omega;{\mathbb R}^3)} \rightarrow 0,\ \text{as}\ l\rightarrow\infty\ .
\end{equation}

Similarly, for ${w_l}$
\begin{subequations}
\label{wl}
\begin{align}
&\text{(i)}\  ( {w_l}_{t}, \phi)_{L^2(\Omega;{\mathbb R}^3)}
+{\nu} (\tilde a_k^r {w_l},_r,\ \tilde a_k^s \phi,_s)_{L^2(\Omega_0^f;{\mathbb R}^3)} 
+ (c^{ijkl}   \int_0^{\cdot}{w_l}^k,_l, \phi^i,_j)_{L^2(\Omega_0^s;{\mathbb R})}  \nonumber \\
&\qquad - (\tilde a_i^j q_l,\ \phi^i,_j)_{L^2(\Omega_0^f;{\mathbb R}^3)}\nonumber\\
 &\qquad  =  (F, \phi)_{L^2 (\Omega_0^f;{\mathbb R}^3)}
+ (f,\phi)_{L^2 (\Omega_0^s;{\mathbb R}^3)} + c_l(\phi) t+d_l (\phi), \ \forall \phi\in 
span(e_1,...,e_l)\,,\\
&\text{(ii)}\  
w_l (0)=(u_0)_l,\ \text{in}\ \Omega\ ,
\end{align}
\end{subequations}
where $d_l (\phi)\in \mathbb R$ is given by
\begin{align*}
d_l(\phi)=&\ ( (w_1)_{l}, \phi)_{L^2(\Omega;{\mathbb R}^3)}
+{\nu} (\tilde a_k^r(0) {u_0}_{l},_r,\ \tilde a_k^s(0) \phi,_s)_{L^2(\Omega_0^f;{\mathbb R}^3)} 
  \nonumber \\ & - ((\tilde a_i^j q_l) (0),\ \phi^i,_j)_{L^2(\Omega_0^f;{\mathbb R})}- (F(0), \phi)_{L^2 (\Omega_0^f;{\mathbb R}^3)}
- (f (0),\phi)_{L^2 (\Omega_0^s;{\mathbb R}^3)} \ . 
\end{align*}
 Similarly as for $c_l (\phi)$, from our compatibility conditions (\ref{compatibility})
\begin{equation}
\label{dl}
\|d_l\|_{H^{-1} (\Omega;{\mathbb R}^3)} \rightarrow 0,\ \text{as}\ l\rightarrow\infty\ .
\end{equation}

We can thus infer now that at the limit $w_{\epsilon}$ satisfies for all $\phi\in L^2(0,T;H^1_0(\Omega;{\mathbb R}^3))$,
 \begin{align}
& \int_0^T ( {w_{\epsilon}}_t, \phi)_{L^2(\Omega;{\mathbb R}^3)}\ dt  
+\nu \int_0^T (\tilde a_k^r   {w_{\epsilon}},_r,\  \tilde a_k^s \phi,_s)_{L^2(\Omega_0^f;{\mathbb R}^3)}\ dt \nonumber\\
&+ \int_0^T (c^{ijkl}\int_0^t {w_{\epsilon}}^k,_l, {\phi}^i,_j)_{L^2(\Omega_0^s;{\mathbb R})}\ dt  - \int_0^T ( q_{\epsilon},\ \tilde a_k^l \phi^k,_l)_{L^2(\Omega_0^f;{\mathbb R})} dt\nonumber\\
&\qquad\qquad = \int_0^T (F, \phi)_{L^2(\Omega_0^f;{\mathbb R}^3)}+ (f, \phi)_{L^2(\Omega_0^s;{\mathbb R}^3)}  \ dt \ ,
\label{weakwepsilon}
\end{align}
which, combined with $w_{\epsilon}(0)=u_0$, shows us that $w_{\epsilon}$ is a 
weak solution of (\ref{weakpenalty}).

Moreover, ${w_{\epsilon}}_t$ satisfies for all $\phi\in L^2(0,T;H^1_0(\Omega;{\mathbb R}^3))$,
 \begin{align}
& \int_0^T ( {w_{\epsilon}}_{tt}, \phi)_{L^2(\Omega;{\mathbb R}^3)}\ dt  
+\nu \int_0^T ((\tilde a_k^s \tilde a_k^r   {w_{\epsilon}},_r)_t,\  \phi,_s)_{L^2(\Omega_0^f;{\mathbb R}^3)}\ dt \nonumber\\
&+ \int_0^T (c^{ijkl}{w_{\epsilon}}^k,_l, {\phi}^i,_j)_{L^2(\Omega_0^s;{\mathbb R})}\ dt  - \int_0^T ( (\tilde a_k^l q_{\epsilon})_t,\  \phi^k,_l)_{L^2(\Omega_0^f;{\mathbb R})} dt\nonumber\\
&\qquad\qquad = \int_0^T (F_t, \phi)_{L^2(\Omega_0^f;{\mathbb R}^3)}+ (f_t, \phi)_{L^2(\Omega_0^s;{\mathbb R}^3)}  \ dt \ .
\label{weakwepsilont}
\end{align}

\noindent{\bf Step 3. Strong convergence for the Galerkin approximation.}

Since $w_{\epsilon}\in L^2(0,T;H^1_0(\Omega;{\mathbb R}^3))$, we can use it
as a test function in (\ref{weakwepsilon}), which provides us on $(0,T)$  
with the equality
\begin{align}
& \frac{1}{2} \|{w_{\epsilon}}(t)\|^2_{L^2(\Omega;{\mathbb R}^3)}  
+\nu \int_0^t (\tilde a_k^r {w_{\epsilon}},_r,\  \tilde a_k^s {w_{\epsilon}},_s)_{L^2(\Omega_0^f;{\mathbb R}^3)}\ dt \nonumber\\
&+ \frac{1}{2}(c^{ijkl}\int_0^t  {w_{\epsilon}}^k,_l, \int_0^t {w_{\epsilon}}^i,_j)_{L^2(\Omega_0^s;{\mathbb R})}  + \int_0^t \epsilon\|q_{\epsilon}\|^2_{L^2(\Omega_0^f;{\mathbb R})}-\epsilon (q_0+t q_1,\ q_{\epsilon})_{L^2(\Omega_0^f;{\mathbb R})} dt\nonumber\\
&\qquad\qquad = \frac{1}{2} \| u_0\|^2_{L^2(\Omega;{\mathbb R}^3)}+ \int_0^t (F, w_{\epsilon})_{L^2(\Omega_0^f;{\mathbb R}^3)}+ (f, w_{\epsilon})_{L^2(\Omega_0^s;{\mathbb R}^3)}  
 \ dt \ .
\label{twice3}
\end{align}
Similarly since $w_l(t)\in span(e_1,...,e_l)$ for all $t\in [0,T]$, we can use it as a test function in (\ref{wl}), which gives us 
\begin{align}
& \frac{1}{2} \|{w_l}(t)\|^2_{L^2(\Omega;{\mathbb R}^3)}  
+\nu \int_0^t (\tilde a_k^r {w_l},_r,\  \tilde a_k^s {w_l},_s)_{L^2(\Omega_0^f;{\mathbb R}^3)}\ dt \nonumber\\
&+\frac{1}{2} (c^{ijkl}\int_0^t  {w_l}^k,_l, \int_0^t {w_{\epsilon}}^i,_j)_{L^2(\Omega_0^s;{\mathbb R})}  + \int_0^t \epsilon \|q_l\|^2_{L^2(\Omega_0^f;{\mathbb R})}-\epsilon (q_0+tq_1,\ q_l)_{L^2(\Omega_0^f;{\mathbb R})} dt\nonumber\\
&= \frac{1}{2} \| (u_0)_l\|^2_{L^2(\Omega;{\mathbb R}^3)}+ \int_0^t (F, w_{\epsilon})_{L^2(\Omega_0^f;{\mathbb R}^3)}+ (f, w_{\epsilon})_{L^2(\Omega_0^s;{\mathbb R}^3)}  
 \ dt \nonumber\\
& \qquad+ \int_0^t t d_l(w_l)+c_l(w_l)\ dt .
\label{twice4}
\end{align}
By (\ref{twice2}), (\ref{cl}) and (\ref{dl}), we then infer by comparing (\ref{twice4}) and (\ref{twice3}), that as $l\rightarrow\infty$, for all $t\in [0,T]$,
\begin{align*}
& \frac{1}{2} \|{w_l}(t)\|^2_{L^2(\Omega;{\mathbb R}^3)}  
+\nu \int_0^t (\tilde a_k^r {w_l},_r,\  \tilde a_k^s {w_l},_s)_{L^2(\Omega_0^f;{\mathbb R}^3)}\ dt \nonumber\\
&+\frac{1}{2} (c^{ijkl}\int_0^t  {w_l}^k,_l, \int_0^t {w_{\epsilon}}^i,_j)_{L^2(\Omega_0^s;{\mathbb R})}\ dt  + \epsilon\int_0^t \|q_l\|^2_{L^2(\Omega_0^f;{\mathbb R})} dt\rightarrow \\
& \frac{1}{2} \|{w_{\epsilon}}(t)\|^2_{L^2(\Omega;{\mathbb R}^3)}  
+\nu \int_0^t (\tilde a_k^r {w_{\epsilon}},_r,\  \tilde a_k^s {w_{\epsilon}},_s)_{L^2(\Omega_0^f;{\mathbb R}^3)}\ dt \nonumber\\
&+\frac{1}{2} (c^{ijkl}\int_0^t  {w_{\epsilon}}^k,_l, \int_0^t {w_{\epsilon}}^i,_j)_{L^2(\Omega_0^s;{\mathbb R})}\ dt  + \epsilon\int_0^t \|q_{\epsilon}\|^2_{L^2(\Omega_0^f;{\mathbb R})} dt\nonumber\ ,
\end{align*}
which gives in turn the strong convergences
\begin{subequations}
\label{twice5}
\begin{align}
w_{l}&\rightarrow  w_{\epsilon}\ \ \text{in}\ L^2(0,T; H^1(\Omega_0^f;{\mathbb R}^3))
\\
{w_l} &\rightarrow   {w_{\epsilon}}\ \ \text{in}\ L^2(0,T; L^2(\Omega;{\mathbb R}^3))
\\ 
\int_0^{\cdot}{w_l} &\rightarrow \int_0^{\cdot}{w_{\epsilon}}\ \ \text{in}\ L^2(0,T; H^1(\Omega_0^s;{\mathbb R}^3))\\ 
{q_l}&\rightarrow{q_{\epsilon}}\ \ \text{in}\ L^2(0,T; L^2(\Omega_0^f;{\mathbb R}))\ .
\end{align}
\end{subequations}
Since ${w_{\epsilon}}_t\in L^2(0,T; H^1_0(\Omega;{\mathbb R}^3))$, we can
prove in a similar fashion the strong convergences
\begin{subequations}
\label{twice6}
\begin{align}
{w_{l}}_t&\rightarrow  {w_{\epsilon}}_t\ \ \text{in}\ L^2(0,T; H^1(\Omega_0^f;{\mathbb R}^3))
\\
{{w_l}}_t &\rightarrow   {w_{\epsilon}}_t\ \ \text{in}\ L^2(0,T; L^2(\Omega;{\mathbb R}^3))
\\ 
{w_l} &\rightarrow {w_{\epsilon}}\ \ \text{in}\ L^2(0,T; H^1(\Omega_0^s;{\mathbb R}^3))\\ 
{q_l}_t&\rightarrow{q_{\epsilon}}_t\ \ \text{in}\ L^2(0,T; L^2(\Omega_0^f;{\mathbb R}))\ .
\end{align}
\end{subequations}

\noindent{\bf Step 4. Energy inequality for ${w_{\epsilon}}_{tt}$.}

By using the relation
 \begin{align*}
&  \frac{1}{2}  \|{w_l}_{tt}\|^2_{L^2(\Omega;{\mathbb R}^3)} 
+{\nu} \int_0^t (\tilde a_k^r {w_l}_{tt},_r,\  \tilde a_k^s {w_l}_{tt},_s)_{L^2(\Omega_0^f;{\mathbb R}^3)} 
+ \frac{1}{2}  (c^{ijkl}   {w_l}_{t}^k,_l, {w_l}_{t}^i,_j)_{L^2(\Omega_0^s;{\mathbb R})}  \nonumber \\
 & + \epsilon \int_0^t\int_{\Omega_0^f} {q_l}^2_{tt}+\int_0^t\int_{\Omega_0^f} {q_l}_{tt}\ [ 2(\tilde a_i^j)_t 
 {w_l}_t^i,_j + (\tilde a_i^j)_{tt} 
 {w_l}^i,_j] -2 \int_0^t\int_{\Omega_0^f} (\tilde a_i^j)_t {q_l}_t{w_l}_{tt}^i,_j\nonumber \\
& - \int_0^t\int_{\Omega_0^f} (\tilde a_i^j)_{tt} {q_l}{w_l}_{tt}^i,_j + {\nu} \int_0^t ((\tilde a_k^r \tilde a_k^s)_{tt} {w_l},_r,\  {w_l}_{tt},_s)_{L^2(\Omega_0^f;{\mathbb R}^3)}\nonumber \\
& +2\ {\nu} \int_0^t ((\tilde a_k^r \tilde a_k^s)_{t} {{w_l}_t},_r,\  {w_{\epsilon}}_{tt},_s)_{L^2(\Omega_0^f;{\mathbb R}^3)}\nonumber\\
&\qquad\le C\ N(u_0,f)^2 +\int_0^t (F_{tt}, {w_l}_{tt})_{L^2 (\Omega_0^f;{\mathbb R}^3)}+\int_0^t (f_{tt},{w_l}_{tt})_{L^2 (\Omega_0^s;{\mathbb R}^3)} ,
\end{align*}
and the weak convergences (\ref{twice2}), along with the strong convergences
(\ref{twice5}) and (\ref{twice6}), we then get 
 \begin{align}
&  \frac{1}{2}  \|{w_{\epsilon}}_{tt}\|^2_{L^2(\Omega;{\mathbb R}^3)} 
+{\nu} \int_0^t (\tilde a_k^r {w_{\epsilon}}_{tt},_r,\  \tilde a_k^s {w_{\epsilon}}_{tt},_s)_{L^2(\Omega_0^f;{\mathbb R}^3)} 
+ \frac{1}{2}  (c^{ijkl}   {w_{\epsilon}}_{t}^k,_l, {w_{\epsilon}}_{t}^i,_j)_{L^2(\Omega_0^s;{\mathbb R})}  \nonumber \\
 & + \epsilon \int_0^t\int_{\Omega_0^f} {q_{\epsilon}}^2_{tt}+\int_0^t\int_{\Omega_0^f} {q_{\epsilon}}_{tt}\ [ 2(\tilde a_i^j)_t 
 {w_{\epsilon}}_t^i,_j + (\tilde a_i^j)_{tt} 
 {w_{\epsilon}}^i,_j] -2 \int_0^t\int_{\Omega_0^f} (\tilde a_i^j)_t {q_{\epsilon}}_t{w_{\epsilon}}_{tt}^i,_j\nonumber \\
& - \int_0^t\int_{\Omega_0^f} (\tilde a_i^j)_{tt} {q_{\epsilon}}{w_{\epsilon}}_{tt}^i,_j + {\nu} \int_0^t ((\tilde a_k^r \tilde a_k^s)_{tt} {w_{\epsilon}},_r,\  {w_{\epsilon}}_{tt},_s)_{L^2(\Omega_0^f;{\mathbb R}^3)}\nonumber \\
& +2\ {\nu} \int_0^t ((\tilde a_k^r \tilde a_k^s)_{t} {{w_{\epsilon}}_t},_r,\  {w_{\epsilon}}_{tt},_s)_{L^2(\Omega_0^f;{\mathbb R}^3)}\nonumber\\
&\qquad\le C\ N(u_0,f)^2 +\int_0^t (F_{tt}, {w_{\epsilon}}_{tt})_{L^2 (\Omega_0^f;{\mathbb R}^3)}+\int_0^t (f_{tt},{w_{\epsilon}}_{tt})_{L^2 (\Omega_0^s;{\mathbb R}^3)} .
\label{wepsilontt}
\end{align}

\subsection{Existence of $\tilde w$, $\tilde w_t$, $\tilde w_{tt}$, uniqueness}
 
In this section, we establish the existence of $\tilde w$, and its first and second
time derivatives by taking the limit $\epsilon\rightarrow 0$. The inequality (\ref{energywtt}) proved at the end of this section, holds for any regularized velocity
field $\tilde v=v_n$, independently of $n$,  and requires in its proof, 
strong convergence results from their penalized counterparts 
 since the regularity that we take on the data does not allow us to view
$\tilde w_{tt}$ as a weak solution of a variational problem. 

\begin{theorem}
\label{uniqueweak}
Suppose  that $u_0$ and $f$ satisfy the conditions stated in 
Theorem \ref{main}. Then, there exists a weak solution $\tilde w$
to the problem (\ref{linear}) with the mollified coefficients replacing
the actual coefficients.   Moreover, $\tilde w$ is in 
$L^2 (0,T;{\mathcal V}_{\tilde v}(\cdot))$ and is unique, and
$\tilde w_t\in {\mathcal W}([0,T])$.
\end{theorem}

\begin{proof} 
\noindent{\bf Step 1. The limit as $\epsilon \rightarrow 0$.}

Let $\epsilon=\frac{1}{m}$; we first pass to the weak limit as 
$m\rightarrow \infty$.  
The energy law (\ref{twice3}) shows that there exists a subsequence 
$\{w_{\frac{1}{m_l}}\}$ such that
\begin{equation}\label{weakconvergence}
w_{\frac{1}{m_l}} \rightharpoonup \tilde w \ \ \text{ in } \ \ {\mathcal W} ([0,T]).
\end{equation}

Moreover, since (\ref{twice3}) also shows that $\|\tilde a_i^j w_{\frac{1}{m}}^i,_j\|_{L^2(0,T;L^2(\Omega_0^f; {\mathbb R}))}\rightarrow 0$ as $m\rightarrow\infty$, we then have
$\|\tilde a_i^j \tilde w^i,_j\|_{L^2(0,T;L^2(\Omega_0^f; {\mathbb R}))}=0$, {\it i.e.}
\begin{equation}
\tilde w\in \mathcal V_{\tilde v}([0,T])\ .
\end{equation}

\noindent{\bf Step 2. The penalized time differentiated problems and 
estimates independent of $\epsilon$.}

Thanks to (\ref{twice6}) and (\ref{wepsilontt}), we have ${w_{\epsilon}}_t\in  L^2(0,T; H^1_0(\Omega; {\mathbb R}^3))\cap C^0([0,T]; L^2(\Omega; {\mathbb R}^3))$. We can thus use it as a test function in (\ref{weakwepsilont}), which gives us for a.e. $t\in (0,T)$
\begin{align*}
&  \frac{1}{2}  \|{w_{\epsilon}}_{t}\|^2_{L^2(\Omega;{\mathbb R}^3)} 
+{\nu} \int_0^t (\tilde a_k^r {w_{\epsilon}}_t,_r,\  \tilde a_k^s {w_{\epsilon}}_t,_s)_{L^2(\Omega_0^f;{\mathbb R}^3)} 
+ \frac{1}{2}  (c^{ijkl}   w_{\epsilon}^k,_l, w_{\epsilon}^i,_j)_{L^2(\Omega_0^s;{\mathbb R})}  \nonumber \\ 
& + \epsilon \int_0^t\int_{\Omega_0^f} {q_{\epsilon}}^2_t -\epsilon \int_0^t\int_{\Omega_0^f} {q_{\epsilon}}_t q_1+\int_0^t \int_{\Omega_0^f}{q_{\epsilon}}_t (\tilde a_i^j)_t 
 {w_{\epsilon}}^i,_j -\int_0^t\int_{\Omega_0^f} (\tilde a_i^j)_t q_{\epsilon}{w_{\epsilon}}_t^i,_j\\
& + {\nu} \int_0^t ((\tilde a_k^r \tilde a_k^s)_t {w_{\epsilon}},_r,\  {w_{\epsilon}}_t,_s)_{L^2(\Omega_0^f;{\mathbb R}^3)}\nonumber\\
&\qquad\qquad\le C\ \|w_1\|^2_{L^2(\Omega;{\mathbb R}^3)}+ \int_0^t ( F_t, {w_{\epsilon}}_t)_{L^2(\Omega_0^f;{\mathbb R}^3)} 
+\int_0^t (f_t,{w_{\epsilon}}_t)_{L^2 (\Omega_0^s;{\mathbb R}^3)} \ .
\end{align*}
At this stage, we remove the time derivative from the  ${q_{\epsilon}}_t$ term
in this inequality by integrating by parts:
$$\int_0^t {q_{\epsilon}}_t (\tilde a_i^j)_t {w_{\epsilon}}^i,_j=
-\int_0^t q_{\epsilon} ((\tilde a_i^j)_t {w_{\epsilon}}^i,_j)_t + (q_{\epsilon} (\tilde a_i^j)_t {w_{\epsilon}}^i,_j) (t) - q_0 (\tilde a_i^j)_t(0) u_0^i,_j$$
($q_{\epsilon}(0)=q_0$ by $\operatorname{div} u_0=0$); 
we then infer by the regularity of $\tilde a$ and of $w_{\epsilon}$ that
\begin{align*}
&  \frac{1}{2}  \|{w_{\epsilon}}_{t}\|^2_{L^2(\Omega;{\mathbb R}^3)} 
+\frac{\nu}{2} \int_0^t (\tilde a_k^r {w_{\epsilon}}_t,_r,\  \tilde a_k^s {w_{\epsilon}}_t,_s)_{L^2(\Omega_0^f;{\mathbb R}^3)} 
+ \frac{1}{2}  (c^{ijkl}   w_{\epsilon}^k,_l, w_{\epsilon}^i,_j)_{L^2(\Omega_0^s;{\mathbb R})}  \nonumber \\ 
& \le \tilde C\ \int_0^t \|q_{\epsilon}\|^2_{L^2(\Omega_0^f;{\mathbb R})}+ \delta\
\|q_{\epsilon}(t)\|^2_{L^2(\Omega_0^f;{\mathbb R})} +\tilde C\ C_{\delta}\ \|\nabla w_{\epsilon}(t)\|^2_{L^2(\Omega_0^f;{\mathbb R}^9)}\nonumber\\
&\qquad+ \tilde C\ N(u_0,f)^2 +C\int_0^t \|{w_{\epsilon}}_t\|^2_{L^2(\Omega;{\mathbb R}^3)}\ ,
\end{align*}
where $\delta>0$ is arbitrary, and $\tilde C$ denotes a generic constant depending on the smoothing parameter $n$ implicit in $\tilde a$.

Note that it is the presence of $\tilde {a}_{tt}$ which requires the use of the 
regularized coefficient matrix  $\tilde a$; this is due to the fact that 
$a_{tt}(t)$
is not in
$\L^{\infty}$ as the presence of $\nabla w_{\epsilon}$ 
and $q_{\epsilon}$ (both taken in $L^2$) would require. 
In order for us to be able to obtain consistent estimates later on which are
independent of the regularization process, we must require the pressure of the 
penalized problem to be in  $H^1$ (for a.e. $t \in (0,T)$); this requires
difference quotient methods. In order to achieve this, we first define this 
pressure function to be  in $L^2$ (a.e. $t \in (0,T)$), and then find estimates
which are  independent of the regularization of $a$.
Thus, in $(0,T)$,
\begin{align}
   &\|{w_{\epsilon}}_{t} (t)\|^2_{L^2(\Omega;{\mathbb R}^3)} +
\|{w_{\epsilon}} (t)\|^2_{H^1(\Omega_0^s;{\mathbb R}^3)}+ \int_0^t \|{w_{\epsilon}}_t\|^2_{H^1(\Omega_0^f;{\mathbb R}^3)}\nonumber\\
&\le\  \tilde C\ \int_0^t \|q_{\epsilon}\|^2_{L^2(\Omega_0^f;{\mathbb R})}+ \delta\
\|q_{\epsilon}(t)\|^2_{L^2(\Omega_0^f;{\mathbb R})} \nonumber\\
&\qquad +\tilde C\ C_{\delta}\ \|\nabla w_{\epsilon}(t)\|^2_{L^2(\Omega_0^f;{\mathbb R}^9)}+ \tilde C\ N(u_0,f)^2 +C\int_0^t \|{w_{\epsilon}}_t\|^2_{L^2(\Omega;{\mathbb R}^3)}\ .
\label{qepsilon}
\end{align}
By the Lagrange multiplier Lemma \ref{Lagrange}, we also have
\begin{align*}
\|q_{\epsilon}(t)\|^2_{L^2(\Omega_0^f;{\mathbb R})}\le\ & C\ [\ \|{w_{\epsilon}}_t (t)\|^2_{L^2(\Omega;{\mathbb R}^3)}+\|\nabla{w_{\epsilon}} (t)\|^2_{L^2(\Omega_0^f;{\mathbb R}^9)}+ \|\nabla \int_0^t w_{\epsilon}\|^2_{L^2(\Omega_0^s;{\mathbb R}^9)}\nonumber\\
&\qquad +
 N(u_0,f)^2\ ]
 \ ,
\end{align*}
which coupled with (\ref{qepsilon}) and (\ref{twice3}), gives for a choice of $\delta>0$ small enough
\begin{align*}
\|q_{\epsilon}(t)\|^2_{L^2(\Omega_0^f;{\mathbb R})}\le\  \tilde C\ [ &\int_0^t \|q_{\epsilon}\|^2_{L^2(\Omega_0^f;{\mathbb R})} + \|\nabla{w_{\epsilon}} (t)\|^2_{L^2(\Omega_0^f;{\mathbb R}^9)}+
N(u_0,f)^2 ]\\
& +C\int_0^t \|{w_{\epsilon}}_t\|^2_{L^2(\Omega;{\mathbb R}^3)}
 \ .
\end{align*}
Since $\displaystyle\int_0^t \|\nabla{w_{\epsilon}} (t)\|^2_{L^2(\Omega_0^f;{\mathbb R}^9)}\le C\ N(u_0,f)^2$, we get by
Gronwall's inequality an estimate on 
$\displaystyle \int_0^t \|q_{\epsilon}\|^2_{L^2(\Omega_0^f;{\mathbb R})} $
which in turn provides 
\begin{equation}
\|q_{\epsilon}(t)\|^2_{L^2(\Omega_0^f;{\mathbb R})}\le\  \tilde C\ [\ \|\nabla{w_{\epsilon}} (t)\|^2_{L^2(\Omega_0^f;{\mathbb R}^9)}+ N(u_0,f)^2\ ]+\tilde C\int_0^t \|{w_{\epsilon}}_t\|^2_{L^2(\Omega;{\mathbb R}^3)}
 \ .
\label{qepsilon2}
\end{equation}
Combined with (\ref{qepsilon}), still for $\delta>0$ small enough, this also gives
\begin{align*}
& \|{w_{\epsilon}}_t (t)\|^2_{L^2(\Omega;{\mathbb R}^3)}+ \|{w_{\epsilon}}(t)\|^2_{H^1(\Omega_0^s;{\mathbb R}^3)}+ \int_0^t \|{w_{\epsilon}}_t\|^2_{H^1(\Omega_0^f;{\mathbb R}^3)}\nonumber\\
&\qquad\qquad\le \  \tilde C\ N (u_0,f)^2
+\tilde C\ \|{w_{\epsilon}} (t)\|^2_{H^1(\Omega_0^f;{\mathbb R}^3)}+\tilde C\int_0^t \|{w_{\epsilon}}_t\|^2_{L^2(\Omega;{\mathbb R}^3)}
 \ .
\end{align*}
By Gronwall and (\ref{twice3}), we first deduce a bound on $\displaystyle \int_0^t \|{w_{\epsilon}}_t\|^2_{L^2(\Omega;{\mathbb R}^3)}$ which in turn provides us with 
%\begin{align*}
%& \int_0^t \|{w_{\epsilon}}_t (t)\|^2_{L^2(\Omega_0^f;{\mathbb R}^3)}\le \  \tilde C\ N (u_0,f)^2\ ,
%\end{align*}
%and thus
\begin{align}
 \|{w_{\epsilon}}_t (t)\|^2_{L^2(\Omega;{\mathbb R}^3)}+ \|{w_{\epsilon}}(t)\|^2_{H^1(\Omega_0^s;{\mathbb R}^3)}+ \int_0^t \|{w_{\epsilon}}_t\|^2_{H^1(\Omega_0^f;{\mathbb R}^3)}\le &\  \tilde C\ N (u_0,f)^2\nonumber\\
& +\tilde C\ \|{w_{\epsilon}} (t)\|^2_{H^1(\Omega_0^f;{\mathbb R}^3)} \ .
\label{wepsilont1}
\end{align}

\noindent{\bf Step 3. An estimate of ${w_{\epsilon}}_t$ on $[0,T]$ which is
independent of $\epsilon$.}

By using $\displaystyle w_{\epsilon}(t)=u_0+\int_0^t {w_{\epsilon}}_t$, we see that
\begin{align*}
 \|{w_{\epsilon}}_t (t)\|^2_{L^2(\Omega;{\mathbb R}^3)}&+ \|{w_{\epsilon}}\|^2_{H^1(\Omega_0^s;{\mathbb R}^3)}+ \int_0^t \|{w_{\epsilon}}_t\|^2_{H^1(\Omega_0^f;{\mathbb R}^3)}\nonumber\\
&\le \tilde C\ N (u_0,f)^2 +\tilde C_1\ t\ \int_0^t 
\|{{w_{\epsilon}}}_t\|^2_{H^1(\Omega_0^f;{\mathbb R}^3)}
+\tilde C\ \|u_0\|^2_{H^1(\Omega_0^f;{\mathbb R}^3)}\nonumber\\
& \le  \tilde C\ N (u_0,f)^2 +\tilde C_1\ t\ \int_0^t 
\|{w_{\epsilon}}_t\|^2_{H^1(\Omega_0^f;{\mathbb R}^3)}\ ,
\end{align*}
where we denote by $\tilde C_1$ a constant, dependent on the smoothing parameter of $\tilde a$ (but not on $\epsilon$), which will remain unchanged in the following estimates.

Now, we see that for any $0\le t\le t_1=\text{Min} (T_M,\frac{1}{2\tilde C_1})$, we have
\begin{align*}
 \|{w_{\epsilon}}_t (t)\|^2_{L^2(\Omega;{\mathbb R}^3)}&+ \|{w_{\epsilon}}\|^2_{H^1(\Omega_0^s;{\mathbb R}^3)}+ \frac{1}{2}\ \int_0^t \|{w_{\epsilon}}_t\|^2_{H^1(\Omega_0^f;{\mathbb R}^3)}\le \tilde C\ N (u_0,f)^2 \ ,
\end{align*}
which with $\displaystyle w_{\epsilon}(t_1)=u_0+\int_0^{t_1} {w_{\epsilon}}_t$ gives 
\begin{equation}
\label{wepsilonT1}
\|w_{\epsilon}(t_1)\|^2_{H^1(\Omega_0^f;{\mathbb R}^3)}\le \tilde C\ N(u_0,f)^2\ .
\end{equation} 
Next, we take $t\ge t_1$ and let  
$\displaystyle w_{\epsilon}(t)=w_{\epsilon}(t_1)+\int_{t_1}^t {w_{\epsilon}}_t$;
we have from (\ref{wepsilont1}) and (\ref{wepsilonT1}) that
\begin{align*}
 \|{w_{\epsilon}}_t (t)\|^2_{L^2(\Omega;{\mathbb R}^3)}&+ \|{w_{\epsilon}}\|^2_{H^1(\Omega_0^s;{\mathbb R}^3)}+ \int_0^t \|{w_{\epsilon}}_t\|^2_{H^1(\Omega_0^f;{\mathbb R}^3)}\nonumber\\
&\le \tilde C\ N (u_0,f)^2 +\tilde C_1\ (t-t_1)\ \int_{t_1}^t 
\|{w_{\epsilon}}_t\|^2_{H^1(\Omega_0^f;{\mathbb R}^3)}
+\tilde C\ \|w_{\epsilon}(t_1)\|^2_{H^1(\Omega_0^f;{\mathbb R}^3)}\nonumber\\
& \le  \tilde C\ N (u_0,f)^2 +\tilde C_1\ (t-t_1)\ \int_0^t 
\|{w_{\epsilon}}_t\|^2_{H^1(\Omega_0^f;{\mathbb R}^3)}\ .
\end{align*}

Now, we see that for any $t_1\le t\le 2\ t_1$, we have
\begin{align*}
 \|{w_{\epsilon}}_t (t)\|^2_{L^2(\Omega;{\mathbb R}^3)}&+ \|{w_{\epsilon}}\|^2_{H^1(\Omega_0^s;{\mathbb R}^3)}+ \frac{1}{2}\ \int_0^t \|{w_{\epsilon}}_t\|^2_{H^1(\Omega_0^f;{\mathbb R}^3)}\le \tilde C\ N (u_0,f)^2 \ ,
\end{align*}
which with $\displaystyle w_{\epsilon}(2 t_1)=u_0+\int_0^{2 t_1} {w_{\epsilon}}_t$ gives 
\begin{equation*}
\|w_{\epsilon}(2t_1)\|^2_{H^1(\Omega_0^f;{\mathbb R}^3)}\le \tilde C\ N(u_0,f)^2\ .
\end{equation*} 

We then see by an easy induction argument that for any $t\in (0,T)$,
\begin{align}
 \|{w_{\epsilon}}_t (t)\|^2_{L^2(\Omega;{\mathbb R}^3)}&+ \|{w_{\epsilon}}\|^2_{H^1(\Omega_0^s;{\mathbb R}^3)}+ \frac{1}{2}\ \int_0^t \|{w_{\epsilon}}_t\|^2_{H^1(\Omega_0^f;{\mathbb R}^3)}\le \tilde C\ N (u_0,f)^2 \ .
\label{wepsilont}
\end{align}
(As is evident in the proof, the constant $\tilde C$ grows as $T$ increases and 
thus depends on $T$.)
Thus with (\ref{qepsilon2}), for all $t\in [0,T]$,
\begin{equation}
\label{qepsilon3}
\|q_{\epsilon}(t)\|^2_{L^2(\Omega_0^f;{\mathbb R})}\le \tilde C\ N (u_0,f)^2 \ .
\end{equation}

\noindent{\bf Step 4. Weak convergence and limit problem.}

Since $w_{\epsilon}$ also satisfies (\ref{twice3}), we thus deduce that for
the choice $\epsilon=\frac{1}{m_l}$ there is a subsequence, still noted $w_{\frac{1}{m_l}}$, such that
\begin{subequations}
\label{w}
\begin{alignat}{2}
w_{\frac{1}{m_l}}&\rightharpoonup \tilde w\ \ \text{in}\ L^2(0,T; H^1(\Omega;{\mathbb R}^3))\label{w.a}
\\
{w_{\frac{1}{m_l}}}_t &\rightharpoonup  \tilde w_t\ \ \text{in}\ L^2(0,T; L^2(\Omega;{\mathbb R}^3))
\\ 
q_{\frac{1}{m_l}}&\rightharpoonup \tilde q\ \ \text{in}\ L^2(0,T; L^2(\Omega_0^f;{\mathbb R}))\ .
\end{alignat}
\end{subequations}

By the weak convergences (\ref{w}), we infer from (\ref{weakwepsilon}) that at 
the limit, for each $\phi \in L^2(0,T; H^1_0(\Omega;{\mathbb R}^3))$,	
 \begin{align}
& \int_0^T ( \tilde w_t, \phi)_{L^2(\Omega;{\mathbb R}^3)}\ dt  
+\nu \int_0^T (\tilde a_k^r  \tilde w,_r,\  \tilde a_k^s \phi,_s)_{L^2(\Omega_0^f;{\mathbb R}^3)}\ dt \nonumber\\
&+ \int_0^T (c^{ijkl}\int_0^t \tilde w^k,_l, {\phi}^i,_j)_{L^2(\Omega_0^s;{\mathbb R})}\ dt  - \int_0^T (\tilde q,\ \tilde a_k^l \phi^k,_l)_{L^2(\Omega_0^f;{\mathbb R})} dt\nonumber\\
&\qquad\qquad = \int_0^T (F, \phi)_{L^2(\Omega_0^f;{\mathbb R}^3)}+ (f, \phi)_{L^2(\Omega_0^s;{\mathbb R}^3)}  \ dt \ .
\label{weakW}
\end{align}
Now for the initial condition, we notice that 
$\tilde w\in C^0([0,T]; L^2(\Omega;{\mathbb R}^3))$. From the following identities
which  hold in $L^2(\Omega;{\mathbb R}^3)$,
$$ \tilde w(t)=w(0)+\int_0^t \tilde w_t,\ \ w_{\epsilon}(t)= u_0+\int_0^t {w_{\epsilon}}_t,$$
we deduce from the weak convergence of 
$\displaystyle\int_0^{\cdot}{w_{\epsilon}}_t$ to
$\displaystyle\int_0^{\cdot}\tilde {w}_t$ in 
$L^2(0,T; L^2(\Omega;{\mathbb R}^3))$ that $w(0)=u_0$ in 
$L^2(\Omega;{\mathbb R}^3)$. 
Combined with (\ref{weakW}), this shows that $\tilde w$ is a weak solution of 
(\ref{weakw}) associated to $\tilde v$.

Now, let us prove that the sequences in (\ref{w}) in fact converge strongly.

\noindent{\bf Step 5. Strong convergence.}

Since $\tilde w\in L^2(0,T; H^1_0(\Omega;{\mathbb R}^3))$ we can use 
$\tilde w$ as a test function in (\ref{weakW}), which provides 
%for any $t\in [0,T]$
%\begin{align}
%& \frac{1}{2} \|\tilde  w(t)\|^2_{L^2(\Omega;{\mathbb R}^3)}  
%+\nu \int_0^t (\tilde a_k^r  \tilde w,_r,\  \tilde a_k^s \tilde w,_s)_{L^2(\Omega_0^f;{\mathbb R}^3)}\ dt\nonumber\\
%& + (c^{ijkl}\int_0^{t} \tilde  w^k,_l,\int_0^{t} \tilde {w}^i,_j)_{L^2(\Omega_0^s;{\mathbb R})}\ dt  \nonumber\\
%&\qquad\qquad = \frac{1}{2} \| u_0\|^2_{L^2(\Omega;{\mathbb R}^3)}+\int_0^t (F, \tilde w)_{L^2(\Omega_0^f;{\mathbb R}^3)}+ (f, \tilde w)_{L^2(\Omega_0^s;{\mathbb R}^3)}  \ dt \ .
%\label{weakW2}
%\end{align}
an energy law that we can compare to (\ref{twice3}). By using the weak convergence
in (\ref{w}), and the fact that 
$\|\tilde a_k^l {w_{\epsilon}}^k,_l\|
_{L^2(0,T;L^2(\Omega_0^f;{\mathbb R}^3))}\rightarrow 0$ 
as $\epsilon\rightarrow 0$ from
(\ref{twice3}), we deduce from this comparison that for any $t\in [0,T]$, 
as $\epsilon\rightarrow 0$,
\begin{align*}
& \frac{1}{2} \|{w_{\epsilon}}(t)\|^2_{L^2(\Omega;{\mathbb R}^3)}  
+\nu \int_0^t (\tilde a_k^r {w_{\epsilon}},_r,\  \tilde a_k^s {w_{\epsilon}},_s)_{L^2(\Omega_0^f;{\mathbb R}^3)}\ dt \nonumber\\
&\qquad\qquad\qquad\qquad\qquad\qquad\qquad + \frac{1}{2}(c^{ijkl}\int_0^t  {w_{\epsilon}}^k,_l, \int_0^t {w_{\epsilon}}^i,_j)_{L^2(\Omega_0^s;{\mathbb R})}  \rightarrow \nonumber\\
& \frac{1}{2} \| \tilde w(t)\|^2_{L^2(\Omega;{\mathbb R}^3)}  
+\nu \int_0^t (\tilde a_k^r \tilde  w,_r,\  \tilde a_k^s \tilde w,_s)_{L^2(\Omega_0^f;{\mathbb R}^3)}\ dt\nonumber\\
&\qquad\qquad\qquad\qquad\qquad\qquad\qquad + \frac{1}{2}(c^{ijkl}\int_0^{t}  \tilde w^k,_l,\int_0^{t} \tilde {w}^i,_j)_{L^2(\Omega_0^s;{\mathbb R})}  \ ,
\end{align*}
which with (\ref{w}) precisely gives the strong convergence
\begin{subequations}
\label{wstrong}
\begin{alignat}{2}
w_{\frac{1}{m_l}}&\rightarrow \tilde w\ \ \text{in}\ L^2(0,T; H^1(\Omega_0^f;{\mathbb R}^3))\label{ws.a}
\\ 
w_{\frac{1}{m_l}}(t)&\rightarrow \tilde w(t)\ \ \text{in}\ L^2(\Omega;{\mathbb R}^3)
\ \text{for any}\ t\in[0,T]\ ,\\
\int_0^t w_{\frac{1}{m_l}}&\rightarrow  \int_0^t \tilde w\ \ \text{in}\ H^1(\Omega_0^s;{\mathbb R}^3)
\ \text{for any}\ t\in[0,T]\ .
\end{alignat}
\end{subequations}

\noindent{\bf Step 6. Uniqueness.}

Now, to prove uniqueness, let us assume that there exists another solution $w'$ to 
(\ref{weakw}), such that $w'\in L^2(0,T;{\mathcal V}_{\tilde v}(\cdot))$, 
$w'_t\in L^2(0,T; L^2(\Omega;{\mathbb R}^3))$. By denoting
$\delta w=\tilde w-w'$, we see that $\delta w\in L^2(0,T;{\mathcal V}_v)$ is a
 solution of 
\begin{align*}
 & (i)\ (\delta w_t, \phi)_{L^2(\Omega;{\mathbb R}^3)}
+{\nu} (\tilde a_k^r \delta w,_r,\ \tilde a_k^s \phi,_s)_{L^2(\Omega_0^f;{\mathbb R}^3)} 
+ (c^{ijkl}\int_0^t \delta w^k,_l, \phi^i,_j)_{L^2(\Omega_0^s;{\mathbb R})} =  0, \\&\qquad\qquad \forall \phi\in 
{\mathcal V}_v (t)\ ,\\
& (ii)\ \delta w(0)=0\ \text{in}\ \Omega\ . 
\end{align*}

Since $\delta w(t\cdot) \in L^2(0,T;{\mathcal V}_{\tilde v}(\cdot))$, we can use $\delta w$ as
a test function in (i), which gives a.e. in $(0,T)$
\begin{align*}
\frac{1}{2} \frac{d}{dt} &\left[ \|\delta w\|^2_{L^2(\Omega;{\mathbb R}^3)}+ (c^{ijkl} \int_0^t \delta w^k,_l,\ \delta w^i,_j)_{L^2(\Omega_0^s;{\mathbb R})} \right] \\
&\qquad\qquad\qquad\qquad\qquad\qquad + \nu (\tilde a_k^r \delta w,_r,\ \tilde a_k^s \delta w,_s)_{L^2(\Omega_0^f;{\mathbb R}^3)}=0 
\end{align*}
which, with the condition $\delta w(0)=0$, precisely proves that $\delta w=0$,
establishing the uniqueness of such a solution. 
\end{proof}

We will also need information on $\tilde w_{tt}$.

\begin{theorem}
\label{weaktt}
Let  $\tilde w \in L^2 (0,T;{\mathcal V}_{\tilde v}(\cdot))$  denote the unique
weak solution
of (\ref{weakW}) which is ensured to exist  by Theorem \ref{uniqueweak}.
Then  $\tilde w_t\in L^2 (0,T; H^1_0(\Omega;{\mathbb R}^3))$ and 
$\tilde w_{tt}\in\mathcal{W} ([0,T])$. Furthermore, 
$q_t\in L^2(0,T; L^2(\Omega_0^f;{\mathbb R}))$.
\end{theorem}
\begin{proof}

\noindent{\bf Step 1. Limit as $\epsilon\rightarrow 0$ in (\ref{wepsilontt}).}

 In order to get an estimate independent of $\epsilon$ from
(\ref{wepsilontt}), we integrate by parts in time to remove the
second time derivative on ${q_{\epsilon}}_{tt}$:
\begin{align}
\int_0^t {q_{\epsilon}}_{tt}\ [ 2(\tilde a_i^j)_t 
 {w_{\epsilon}}_t^i,_j + (\tilde a_i^j)_{tt} 
 {w_{\epsilon}}^i,_j] =& -\int_0^t {q_{\epsilon}}_{t}\ [ 2(\tilde a_i^j)_t 
 {w_{\epsilon}}_t^i,_j + (\tilde a_i^j)_{tt} 
 {w_{\epsilon}}^i,_j]_t \nonumber\\
& + {q_{\epsilon}}_t (t) [2 (\tilde a_i^j)_t {w_{\epsilon}}_t^i,_j + (\tilde a_i^j)_{tt} {w_{\epsilon}}^i,_j] (t)\nonumber\\
& -q_1 [2(\tilde a_i^j)_t(0) {w_1}^i,_j + (\tilde a_i^j)_{tt}(0) {u_0}^i,_j] ,
\label{qepsilontt}
\end{align} 
($\displaystyle {q_{\epsilon}}_t(0)=q_1-\frac{1}{\epsilon} [\ (\tilde a_i^j)_t(0) u_0^i,_j + \tilde a_i^j(0) w_1^i,_j\ ]=q_1$ in $\Omega_0^f$ by our compatibility conditions on the initial data), from which we then infer by the regularity of 
$\tilde a$ and of $w_{\epsilon}$ that
\begin{align*}
&  \frac{1}{2}  \|{w_{\epsilon}}_{tt}\|^2_{L^2(\Omega;{\mathbb R}^3)} 
+\nu \int_0^t (\tilde a_k^r {w_{\epsilon}}_{tt},_r,\  \tilde a_k^s {w_{\epsilon}}_{tt},_s)_{L^2(\Omega_0^f;{\mathbb R}^3)} 
+ \frac{1}{2}  (c^{ijkl}   {w_{\epsilon}}_t^k,_l, {w_{\epsilon}}_t^i,_j)_{L^2(\Omega_0^s;{\mathbb R})}  \nonumber \\ 
& \le \tilde C C_{\delta}\int_0^t \|{q_{\epsilon}}_t\|^2_{L^2(\Omega_0^f;{\mathbb R})}+ \delta\ \int_0^t \|\nabla {w_{\epsilon}}_{tt}\|^2_{L^2(\Omega_0^f;{\mathbb R}^9)}+
\delta\ \|{q_{\epsilon}}_t(t)\|^2_{L^2(\Omega_0^f;{\mathbb R})}\nonumber\\
&\qquad +\tilde C\ C_{\delta}\|\nabla {w_{\epsilon}}_t (t)\|^2_{L^2(\Omega_0^f;{\mathbb R}^9)}+\tilde C\ C_{\delta}\ \|\nabla {w_{\epsilon}} (t)\|^2_{L^2(\Omega_0^f;{\mathbb R}^9)}+ C_{\delta} \tilde C N (u_0,f)^2 \\
&\qquad\qquad + C\ \int_0^t \|{w_{\epsilon}}_{tt}\|^2_{L^2(\Omega;{\mathbb R}^3)}\ ,
\end{align*}
where $\delta>0$ is arbitrary.

Thus, for $\delta$ small enough,
\begin{align}
  \|{w_{\epsilon}}_{tt} (t)\|&^2_{L^2(\Omega;{\mathbb R}^3)} + \int_0^t \|\nabla {w_{\epsilon}}_{tt}\|^2_{L^2(\Omega_0^f;{\mathbb R}^9)} 
+  \|\nabla{w_{\epsilon}}_t (t)\|^2_{L^2(\Omega_0^s;{\mathbb R}^9)}  \nonumber \\ 
&\le\  \tilde C\ \int_0^t \|{q_{\epsilon}}_{t}\|^2_{L^2(\Omega_0^f;{\mathbb R})}+ \delta\
\|{q_{\epsilon}}_t (t)\|^2_{L^2(\Omega_0^f;{\mathbb R})} +\tilde C\ C_{\delta}\|\nabla {w_{\epsilon}}_t (t)\|^2_{L^2(\Omega_0^f;{\mathbb R}^9)}\nonumber\\
&\qquad + \tilde C\ N (u_0,f)^2 + C\ \int_0^t \|{w_{\epsilon}}_{tt}\|^2_{L^2(\Omega;{\mathbb R}^3)}\ .
\label{qepsilont}
\end{align}
By the Lagrange multiplier Lemma \ref{Lagrange}, we also have that
\begin{align*}
\|{q_{\epsilon}}_t (t)\|^2_{L^2(\Omega_0^f;{\mathbb R})}\le\ & C\ [\ \|{w_{\epsilon}}_{tt} (t)\|^2_{L^2(\Omega;{\mathbb R}^3)}+\|\nabla{w_{\epsilon}}_t (t)\|^2_{L^2(\Omega_0^f;{\mathbb R}^9)}+ \tilde C \|\nabla  w_{\epsilon}(t)\|^2_{L^2(\Omega_0^f;{\mathbb R}^9)}\nonumber\\
&\qquad + \tilde C\ \|{q_{\epsilon}} (t)\|^2_{L^2(\Omega_0^f;{\mathbb R})}+ \|\nabla  w_{\epsilon}(t)\|^2_{L^2(\Omega_0^s;{\mathbb R}^9)}+
N(u_0,f)^2
\ ]
 \ ,
\end{align*}
and thus with (\ref{wepsilont}), (\ref{qepsilon3}) and (\ref{qepsilont}) for a choice of $\delta>0$ small enough,
\begin{align*}
\|{q_{\epsilon}}_t (t)\|^2_{L^2(\Omega_0^f;{\mathbb R})}\le\ & C\ [\ \tilde C \int_0^t \|{q_{\epsilon}}_{t}\|^2_{L^2(\Omega_0^f;{\mathbb R})}+\tilde C\ \|\nabla{w_{\epsilon}}_t (t)\|^2_{L^2(\Omega_0^f;{\mathbb R}^9)}\nonumber\\
& + \|\nabla  w_{\epsilon}(t)\|^2_{L^2(\Omega_0^s;{\mathbb R}^9)} + \tilde C\ N(u_0,f)^2+ C\ \int_0^t \|{w_{\epsilon}}_{tt}\|^2_{L^2(\Omega;{\mathbb R}^3)}\ ]
 \ ,
\end{align*}
which by Gronwall's inequality and (\ref{wepsilont}), gives
\begin{align*}
\int_0^t \|{q_{\epsilon}}_t \|^2_{L^2(\Omega_0^f;{\mathbb R})}\le\ & \tilde C\ [\ N(u_0,f)^2+  \int_0^t \|{w_{\epsilon}}_{tt}\|^2_{L^2(\Omega;{\mathbb R}^3)}
\ ] \ ,
\end{align*}
and thus,
\begin{align}
\|{q_{\epsilon}}_t (t)\|^2_{L^2(\Omega_0^f;{\mathbb R})}\le\ & \tilde C\ [\ N(u_0,f)^2 + \|\nabla{w_{\epsilon}}_t (t)\|^2_{L^2(\Omega_0^f;{\mathbb R}^9)}+ \int_0^t \|{w_{\epsilon}}_{tt}\|^2_{L^2(\Omega;{\mathbb R}^3)}\ ]
 \ .
\label{qepsilont3}
\end{align}
We then infer from (\ref{qepsilont}) that
\begin{align*}
   \|{w_{\epsilon}}_{tt} (t)\|^2_{L^2(\Omega;{\mathbb R}^3)}& + \int_0^t \|\nabla {w_{\epsilon}}_{tt}\|^2_{L^2(\Omega_0^f;{\mathbb R}^9)} 
+  \|\nabla{w_{\epsilon}}_t (t)\|^2_{L^2(\Omega_0^s;{\mathbb R}^9)}  \nonumber \\ 
&\le\   \tilde C\ [\ \|\nabla {w_{\epsilon}}_{t} (t)\|^2_{L^2(\Omega_0^f;{\mathbb R}^9)}+ \ N (u_0,f)^2 + \int_0^t \|{w_{\epsilon}}_{tt}\|^2_{L^2(\Omega;{\mathbb R}^3)}\ ]\ .
\end{align*}
By the Gronwall inequality and (\ref{wepsilont}), we first get an estimate 
on $\displaystyle
 \int_0^t \|{w_{\epsilon}}_{tt} (t)\|^2_{L^2(\Omega;{\mathbb R}^3)}$ which
implies in turn that
%\begin{align*}
%  \int_0^t \|{w_{\epsilon}}_{tt} (t)\|^2_{L^2(\Omega_0^f;{\mathbb R}^3)}\le\   \tilde C\ \ N (u_0,f)^2 \ ,
%\end{align*}
%and thus
\begin{align*}
  \|{w_{\epsilon}}_{tt} (t)\|^2_{L^2(\Omega;{\mathbb R}^3)} + \int_0^t \|\nabla {w_{\epsilon}}_{tt}\|^2_{L^2(\Omega_0^f;{\mathbb R}^9)}& 
+  \|\nabla{w_{\epsilon}}_t (t)\|^2_{L^2(\Omega_0^s;{\mathbb R}^9)}  \nonumber \\ 
&\le\   \tilde C\ \ \|\nabla {w_{\epsilon}}_{t} (t)\|^2_{L^2(\Omega_0^f;{\mathbb R}^9)}+ \tilde C\ N (u_0,f)^2 \ .
\end{align*}

\noindent{\bf Step 2. $\epsilon$-independent estimate  for ${w_{\epsilon}}_{tt}$ 
on $[0,T]$.}

In the same fashion as we derived (\ref{wepsilont}) from (\ref{wepsilont1}), we
can deduce that for all $t\in [0,T]$,
\begin{align}
  \|{w_{\epsilon}}_{tt} (t)\|^2_{L^2(\Omega;{\mathbb R}^3)} + \int_0^t \|\nabla {w_{\epsilon}}_{tt}\|^2_{L^2(\Omega_0^f;{\mathbb R}^9)} 
+  \|\nabla{w_{\epsilon}}_t (t)\|^2_{L^2(\Omega_0^s;{\mathbb R}^9)} \le\  \tilde C\ N (u_0,f)^2 \ .
\label{wepsilontt3}
\end{align}
From (\ref{qepsilont3}) and (\ref{wepsilontt3}) we then infer
\begin{align}
  \|{q_{\epsilon}}_{t} (t)\|^2_{L^2(\Omega_0^f;{\mathbb R})} \le\   \tilde C\ N(u_0,f)^2\ .
\label{qepsilont2}
\end{align}

We thus deduce that for
the choice $\epsilon=\frac{1}{m_l}$ there is a subsequence, still denoted $w_{\frac{1}{m_l}}$, such that
\begin{subequations}
\label{weakcvwt}
\begin{align}
w_{\frac{1}{m_l}}&\rightharpoonup  \tilde w\ \ \text{in}\ L^2(0,T; H^1(\Omega;{\mathbb R}^3))
\\
{w_{\frac{1}{m_l}}}_t &\rightharpoonup  \tilde w_t\ \ \text{in}\ L^2(0,T; H^1(\Omega;{\mathbb R}^3))
\\ 
{w_{\frac{1}{m_l}}}_{tt} &\rightharpoonup  \tilde w_{tt}\ \ \text{in}\ L^2(0,T; L^2(\Omega;{\mathbb R}^3))\ \text{and in}\ L^2(0,T; H^1(\Omega_0^f;{\mathbb R}^3))\\
q_{\frac{1}{m_l}}&\rightharpoonup \tilde q\ \ \text{in}\ L^2(0,T; L^2(\Omega_0^f;{\mathbb R}))\\
{q_{\frac{1}{m_l}}}_t &\rightharpoonup \tilde q_t\ \ \text{in}\ L^2(0,T; L^2(\Omega_0^f;{\mathbb R}))\ .
\end{align}
\end{subequations}

\noindent{\bf Step 3. Initial condition for $\tilde w_t$.}

By the weak convergence in (\ref{weakcvwt}), we infer from (\ref{weakwepsilont}) 
that for each test function $\phi \in L^2(0,T; H^1_0(\Omega;{\mathbb R}^3))$,	
 \begin{align}
& \int_0^T ( \tilde w_{tt}, \phi)_{L^2(\Omega;{\mathbb R}^3)}\ dt  
+\nu \int_0^T ((\tilde a_k^r \tilde a_k^s  \tilde w,_r)_t,\  \phi,_s)_{L^2(\Omega_0^f;{\mathbb R}^3)}\ dt \nonumber\\
&+ \int_0^T (c^{ijkl} \tilde w^k,_l, {\phi}^i,_j)_{L^2(\Omega_0^s;{\mathbb R})}\ dt  - \int_0^T ((\tilde a_i^j \tilde q)_t,\  \phi^i,_j)_{L^2(\Omega_0^f;{\mathbb R})} dt\nonumber\\
&\qquad\qquad = \int_0^T ( F_t, \phi)_{L^2(\Omega_0^f;{\mathbb R}^3)}+ (f_t, \phi)_{L^2(\Omega_0^s;{\mathbb R}^3)} \ dt \ .
\label{weakWt}
\end{align}
Now for the initial condition, we notice that 
$\tilde w_t\in C^0([0,T]; L^2(\Omega;{\mathbb R}^3))$. From the following 
identities which hold in $L^2(\Omega;{\mathbb R}^3)$, we find that
$$ \tilde w_t(t)=\tilde w_t(0)+\int_0^t \tilde w_{tt},\ \ w_{\epsilon}(t)
= w_1+\int_0^t {w_{\epsilon}}_{tt}.$$
We deduce from the weak convergence of 
$\displaystyle\int_0^{\cdot}{w_{\epsilon}}_{tt}$ to
$\displaystyle\int_0^{\cdot}{\tilde w}_{tt}$ in $L^2(0,T; L^2(\Omega;{\mathbb R}^3))$ that $\tilde w_t(0)=w_1$ in $L^2(\Omega;{\mathbb R}^3)$.

\noindent{\bf Step 4. Strong convergence: the easy cases.}

We will also need the fact that the weak convergence in (\ref{weakcvwt}) is 
in fact strong.  Notice that since $\tilde w_t\in L^2(0,T; H^1_0(\Omega;
{\mathbb R}^3))$, we can use $\tilde w_t$ in (\ref{weakWt}) to get for any 
$t\in [0,T]$,
 \begin{align*}
&\frac{1}{2}\|\tilde w_{t}(t)\|^2_{L^2(\Omega;{\mathbb R}^3)}  
+\nu \int_0^t ((\tilde a_k^r \tilde a_k^s  \tilde w,_r)_t,\  {\tilde w_t},_s)_{L^2(\Omega_0^f;{\mathbb R}^3)} - \int_0^t ((\tilde a_i^j \tilde q)_t,\ \tilde w_t^i,_j)_{L^2(\Omega_0^f;{\mathbb R})} \nonumber\\
&+ \frac{1}{2}(c^{ijkl} \tilde w^k,_l (t), {\tilde w}^i,_j(t))_{L^2(\Omega_0^s;{\mathbb R})} = \frac{1}{2} \|w_1\|^2_{L^2(\Omega;{\mathbb R}^3)}+ \frac{1}{2}(c^{ijkl} u_0^k,_l , {u_0}^i,_j)_{L^2(\Omega_0^s;{\mathbb R})}\nonumber\\
&\qquad\qquad\qquad\qquad\qquad\qquad\qquad\qquad+ \int_0^t (F_t, \tilde w_t)_{L^2(\Omega_0^f;{\mathbb R}^3)}+ (f_t, \tilde w_t)_{L^2(\Omega_0^s;{\mathbb R}^3)}\ dt \ .
\end{align*}
Since $\tilde w(t)\in {\mathcal V}_{\tilde v}(t)$ in $[0,T]$ implies $(\tilde a_i^j)_t \tilde w^i,_j=-\tilde a_i^j \tilde w_t^i,_j$, we then deduce
 \begin{align}
&\frac{1}{2}\|\tilde w_{t}(t)\|^2_{L^2(\Omega;{\mathbb R}^3)}  
+\nu \int_0^t (\tilde a_k^r \tilde a_k^s  {\tilde w_t},_r,\  
{\tilde w_t},_s)_{L^2(\Omega_0^f;{\mathbb R}^3)} 
+ \nu \int_0^t((\tilde a_k^r \tilde a_k^s)_t  {\tilde w},_r,\  
{\tilde w_t},_s)_{L^2(\Omega_0^f;{\mathbb R}^3)} \nonumber\\
& - \int_0^t ((\tilde a_i^j)_t \tilde q,\  \tilde w_t^i,_j)_{L^2(\Omega_0^f;{\mathbb R})} +\int_0^t
((\tilde a_i^j)_t \tilde q_t,\  \tilde w^i,_j)_{L^2(\Omega_0^f;{\mathbb R})}\nonumber\\
&+ \frac{1}{2}(c^{ijkl} \tilde w^k,_l (t), {\tilde w}^i,_j(t))_{L^2(\Omega_0^s;{\mathbb R})}= \frac{1}{2} \|w_1\|^2_{L^2(\Omega;{\mathbb R}^3)}+ \frac{1}{2}(c^{ijkl} u_0^k,_l , {u_0}^i,_j)_{L^2(\Omega_0^s;{\mathbb R})}\nonumber\\
& \qquad\qquad\qquad\qquad\qquad\qquad\qquad\qquad+ \int_0^t (F_t, \tilde w_t)_{L^2(\Omega_0^f;{\mathbb R}^3)}+ (f_t, \tilde w_t)_{L^2(\Omega_0^s;{\mathbb R}^3)}\ dt  .
\label{weakWt2}
\end{align}
Similarly,
\begin{align*}
&\frac{1}{2}\|{w_{\epsilon}}_{t}(t)\|^2_{L^2(\Omega;{\mathbb R}^3)}  
+ \nu\int_0^t ((\tilde a_k^r \tilde a_k^s  {w_{\epsilon}},_r)_t,\ {w_{\epsilon}} ,_s)_{L^2(\Omega_0^f;{\mathbb R}^3)} -\int_0^t ((\tilde a_i^j q_{\epsilon})_t,\  {w_{\epsilon}}_t^i,_j)_{L^2(\Omega_0^f;{\mathbb R})} \nonumber\\
&+\frac{1}{2} (c^{ijkl} {w_{\epsilon}}^k,_l (t), {w_{\epsilon}}^i,_j(t))_{L^2(\Omega_0^s;{\mathbb R})} = \frac{1}{2} \|w_1\|^2_{L^2(\Omega;{\mathbb R}^3)}+ \frac{1}{2}(c^{ijkl} u_0^k,_l , {u_0}^i,_j)_{L^2(\Omega_0^s;{\mathbb R})}\nonumber\\
&\qquad\qquad\qquad\qquad\qquad\qquad\qquad\qquad+ \int_0^t (F_t, {w_{\epsilon}}_t)_{L^2(\Omega_0^f;{\mathbb R}^3)}+ (f_t, {w_{\epsilon}}_t)_{L^2(\Omega_0^s;{\mathbb R}^3)}\ dt \ .
\end{align*}
From the definition of $q_{\epsilon}$, $(\tilde a_i^j)_t {w_{\epsilon}}^i,_j=-\tilde a_i^j {w_{\epsilon}}_t^i,_j -\epsilon\ ({q_{\epsilon}}_t-q_1)$, and thus 
 \begin{align}
&\frac{1}{2}\|{w_{\epsilon}}_{t}(t)\|^2_{L^2(\Omega;{\mathbb R}^3)}  
+\nu \int_0^t (\tilde a_k^r \tilde a_k^s  {w_{\epsilon}}_t,_r,\  {w_{\epsilon}}_t,_s)_{L^2(\Omega_0^f;{\mathbb R}^3)}\nonumber\\
&+\nu\int_0^t ((\tilde a_k^r \tilde a_k^s)_t  {w_{\epsilon}},_r,\  {w_{\epsilon}}_t,_s)_{L^2(\Omega_0^f;{\mathbb R}^3)} - \int_0^t ((\tilde a_i^j)_t q_{\epsilon},\  {w_{\epsilon}}_t^i,_j)_{L^2(\Omega_0^f;{\mathbb R})}\nonumber\\
& +\int_0^t 
((\tilde a_i^j)_t {q_{\epsilon}}_t,\  {w_{\epsilon}}^i,_j)_{L^2(\Omega_0^f;{\mathbb R})} + \epsilon \int_0^t ({q_{\epsilon}}_t,\  {q_{\epsilon}}_t-q_1)_{L^2(\Omega_0^f;{\mathbb R})}\nonumber\\
&+ \frac{1}{2}(c^{ijkl} {w_{\epsilon}}^k,_l (t),\ {w_{\epsilon}}^i,_j(t))_{L^2(\Omega_0^s;{\mathbb R})}= \frac{1}{2} \|w_1\|^2_{L^2(\Omega;{\mathbb R}^3)}+ \frac{1}{2}(c^{ijkl} u_0^k,_l , {u_0}^i,_j)_{L^2(\Omega_0^s;{\mathbb R})}\nonumber\\
&\qquad\qquad\qquad\qquad\qquad\qquad\qquad\qquad+ \int_0^t (F_t, {w_{\epsilon}}_t)_{L^2(\Omega_0^f;{\mathbb R}^3)}+\int_0^t (f_t, {w_{\epsilon}}_t)_{L^2(\Omega_0^s;{\mathbb R}^3)} \ .
\label{weakepsilont2}
\end{align}
By integration by parts, since $q_{\epsilon} (0)=q_0$,
\begin{align}
\int_0^t ((\tilde a_i^j)_t q_{\epsilon},\  {w_{\epsilon}}_t^i,_j)_{L^2(\Omega_0^f;{\mathbb R})}\ dt=&-\int_0^t (((\tilde a_i^j)_t q_{\epsilon})_t,\  {w_{\epsilon}}^i,_j)_{L^2(\Omega_0^f;{\mathbb R})}\ dt\nonumber\\
& +((\tilde a_i^j)_t q_{\epsilon}(t),\  {w_{\epsilon}}^i,_j (t))_{L^2(\Omega_0^f;{\mathbb R})}\nonumber\\
&- ((\tilde a_i^j)_t (0) q_0,\  {u_0}^i,_j )_{L^2(\Omega_0^f;{\mathbb R})} .
\label{strongwt1}
\end{align}
Now, since $\displaystyle q_{\epsilon}(t)=q_0+\int_0^t {q_{\epsilon}}_t\ dt$ we deduce that for any $t\in [0,T]$, $\displaystyle q_{\epsilon}(t)\rightharpoonup q_0+\int_0^t \tilde q_t\ dt$ in $L^2(\Omega_0^f;{\mathbb R})$, which proves that $\displaystyle \tilde q(t)=q_0+\int_0^t \tilde q_t$, and thus that $\tilde q\in C^0 ([0,T]; L^2(\Omega_0^f;{\mathbb R}))$, with  $\tilde q(0)=q_0$. Furthermore, for any $t\in [0,T]$, 
\begin{equation}
q_{\epsilon}(t)\rightharpoonup \tilde q(t)\ \text{in}\ L^2(\Omega_0^f;{\mathbb R})\ \text{as}\ \epsilon\rightarrow 0\ .
\label{strongwt2}
\end{equation}
Similarly, since $w_{\epsilon}(0)=\tilde w(0)=u_0$ we have that for any $t\in [0,T]$, 
\begin{equation*}
w_{\epsilon}(t)\rightharpoonup \tilde w(t)\ \text{in}\ 
H^1(\Omega_0^f;{\mathbb R}^3)\ \text{as}\ \epsilon\rightarrow 0\ .
\end{equation*}
Moreover from
\begin{equation*}
\|w_{\epsilon}(t)\|^2_{H^1(\Omega_0^f;{\mathbb R}^3)}
= \|u_0\|^2_{H^1(\Omega_0^f;{\mathbb R}^3)} +2\int_0^t 
({w_{\epsilon}}_t,\ w_{\epsilon})_{H^1(\Omega_0^f;{\mathbb R}^3)}\ dt\ ,
\end{equation*}
we infer from the strong convergence in (\ref{wstrong}) and the weak 
convergence in
(\ref{weakcvwt}) that for any $t\in [0,T]$, 
$$\|w_{\epsilon}(t)\|^2_{H^1(\Omega_0^f;{\mathbb R}^3)}\rightarrow 
\|u_0\|^2_{H^1(\Omega_0^f;{\mathbb R}^3)}+2\int_0^t 
(\tilde w_t, \tilde w)_{H^1(\Omega_0^f;{\mathbb R}^3)}\ dt\ 
\text{as}\ \epsilon=\frac{1}{m_l}\rightarrow 0,$$ 
from which we obtain the strong convergence
\begin{equation}
w_{\epsilon}(t)\rightarrow \tilde w(t)\ \text{in}\ H^1(\Omega_0^f;{\mathbb R}^3)
\ \text{as}\ \epsilon=\frac{1}{m_l}\rightarrow 0\ .
\label{strongwt3}
\end{equation} 
Thus, from (\ref{strongwt1}), the strong convergence in (\ref{wstrong}) and
(\ref{strongwt3})  together with  the weak convergence  in(\ref{weakcvwt})
and (\ref{strongwt2}) shows that
\begin{equation}
\int_0^t ((\tilde a_i^j)_t q_{\epsilon},\  
{w_{\epsilon}}_t^i,_j)_{L^2(\Omega_0^f;{\mathbb R})}\ dt \rightarrow 
\int_0^t ((\tilde a_i^j)_t\, \tilde q,\  
{\tilde w}_t^i,_j)_{L^2(\Omega_0^f;{\mathbb R})}\ dt\ 
\text{as}\ \epsilon=\frac{1}{m_l}\rightarrow 0\ .
\label{strongwt4}
\end{equation}

From (\ref{strongwt4}), the weak convergence in (\ref{weakcvwt}) and the strong 
convergence in (\ref{wstrong}), we then deduce from (\ref{weakWt2}) and 
(\ref{weakepsilont2}), that as $\epsilon=\frac{1}{m_l}\rightarrow 0$, for 
any $t\in [0,T]$,
 \begin{align*}
&\frac{1}{2}\|{w_{\epsilon}}_{t}(t)\|_{L^2(\Omega;{\mathbb R}^3)}  
+\nu \int_0^t (\tilde a_k^r \tilde a_k^s  {w_{\epsilon}}_t,_r,\  {w_{\epsilon}}_t,_s)_{L^2(\Omega_0^f;{\mathbb R}^3)}\ dt\nonumber\\
&\qquad\qquad\qquad\qquad\qquad + \frac{1}{2}(c^{ijkl} {w_{\epsilon}}^k,_l (t),\ {w_{\epsilon}}^i,_j(t))_{L^2(\Omega_0^s;{\mathbb R})}\rightarrow\nonumber\\
&\frac{1}{2}\|\tilde w_{t}(t)\|_{L^2(\Omega;{\mathbb R}^3)}  
+\nu \int_0^t (\tilde a_k^r \tilde a_k^s \tilde w_t,_r,\  \tilde w_t,_s)_{L^2(\Omega_0^f;{\mathbb R}^3)}\ dt\nonumber\\
&\qquad\qquad\qquad\qquad\qquad+ \frac{1}{2}(c^{ijkl} \tilde w^k,_l (t),\ \tilde w^i,_j(t))_{L^2(\Omega_0^s;{\mathbb R})}\ ,
\end{align*}
which implies the strong convergences
\begin{subequations}
\label{ss_strong}
\begin{alignat}{2}
w_{\frac{1}{m_l'}}(t) &\rightarrow \tilde w(t)\ \ \text{in}\  H^1(\Omega;{\mathbb R}^3)\ \text{for any}\ t\in [0,T],
\\
{w_{\frac{1}{m_l'}}}_t &\rightarrow  \tilde w_t\ \ \text{in}\ L^2(0,T; H^1(\Omega_0^f;{\mathbb R}^3))
\\ 
{w_{\frac{1}{m_l'}}}_{t} (t) &\rightarrow \tilde w_{t} (t)\ \ \text{in}\ L^2(\Omega;{\mathbb R}^3)\ \text{for any}\ t\in [0,T]\ .
\label{strongcvwt}
\end{alignat}
\end{subequations}

From the strong convergence in (\ref{ss_strong}) and Lemma \ref{Lagrange}, we also 
deduce that
\begin{equation}
\label{strongcvq}
\|q_{\epsilon}-\tilde q\|_{L^2(0,T; L^2(\Omega_0^f;{\mathbb R}))}\rightarrow 0,\ \text{as}\ \epsilon\rightarrow 0\ .
\end{equation}

\noindent{\bf Step 5. Strong convergence: the more delicate case of 
$\tilde w_{tt}$.}

Our main difficulty results from the fact that we  cannot directly obtain an
energy inequality for $w_{tt}$ (from the limiting weak form of the twice
time-differentiated problem).
Rather, our starting point will be (\ref{wepsilontt}), from which we 
will get by weak lower semi-continuity the desired inequality, provided that we 
can prove that ${w_{\epsilon}}_{tt}\rightarrow w_{tt}$ 
in $L^2(0,T;L^2(\Omega_0^f;{\mathbb R}^3))$. To prove this result, let us first 
remind the reader that  
 ${w_{\epsilon}}_{tt}$ satisfies for all $\phi\in L^2(0,T;H^1_0(\Omega;{\mathbb R}^3))$,
 \begin{align}
& \int_0^T ( {w_{\epsilon}}_{ttt}, \phi)_{L^2(\Omega;{\mathbb R}^3)}\ dt  
+\nu \int_0^T ((\tilde a_k^s \tilde a_k^r   {w_{\epsilon}},_r)_{tt},\  \phi,_s)_{L^2(\Omega_0^f;{\mathbb R}^3)}\ dt \nonumber\\
&+ \int_0^T (c^{ijkl}{w_{\epsilon}}_t^k,_l, {\phi}^i,_j)_{L^2(\Omega_0^s;{\mathbb R})}\ dt  - \int_0^T ( (\tilde a_k^l q_{\epsilon})_{tt},\  \phi^k,_l)_{L^2(\Omega_0^f;{\mathbb R})} dt\nonumber\\
&\qquad\qquad = \int_0^T (F_{tt}, \phi)_{L^2(\Omega_0^f;{\mathbb R}^3)}+ (f_{tt}, \phi)_{L^2(\Omega_0^s;{\mathbb R}^3)}  \ dt \ .
\label{weakwepsilontt}
\end{align}
From the bounds associated to the weak convergence in (\ref{weakcvwt}), we then
see that
\begin{equation}
\label{twice7}
\int_0^T \|{w_{\epsilon}}_{ttt} (t)\|^2_{{\mathcal V}_{\tilde v}(t)'}\ dt\le \check C,
\end{equation}
where $\check C$ denotes a constant which depends on the data, the smoothing 
parameter implicit in $\tilde a$, but not on the penalization parameter $\epsilon$.
In the following, this letter will denote a generic constant depending on these 
variables. Let us fix
$\delta>0$ and let 
$$\Omega^f_{\delta}=\{x\in\Omega_0^f|\ \text{dist}(x,\Gamma_0)\ge \delta,\ 
\text{dist}(x,\partial\Omega)\ge \delta\}\ .$$ 
Let us then denote on $[0,T]$, 
$$\tilde\Omega_{\delta}(t)=\tilde\eta(t,\Omega^f_{\delta})\ .$$
%We remind the reader  here that $\tilde\eta(t,\Omega)$ is not necessarily included in $\Omega$, whereas $\eta(t,\Omega)=\Omega$. Thus, since $\|v_n-v\|_{V^3_f(T)}\rightarrow 0$ as $n\rightarrow\infty$, let $n_\delta$ be such that for $n\ge n_{\delta}$ and $t\in [0,T]$, $\eta_n(t,\Omega_{\delta})\subset \Omega$. From now on, $\tilde v=v_n$, with
%$n\ge n_{\delta}$. 

For each $t\in (0,T)$, we have the existence of $\delta t>0$ such that
\begin{align*}
&\forall t'\in (t-\delta t,t+\delta t),\ \tilde\Omega_{2\delta}(t)\subset \tilde\Omega_{\delta}(t') .
\end{align*}
By a simple change of variables and (\ref{twice7}), we then get
\begin{equation}
\label{twice8}
\bigl\|\det \tilde a\, {w_{\epsilon}}_{ttt}\circ\tilde\eta^{-1}
\bigr\|_{L^2(t-\delta t,t+\delta t; H^1_{0,\operatorname{div}} (\tilde\Omega_{2\delta}(t);{\mathbb R}^3)')}\le \check C\ .
\end{equation}
We set
$$u_{\epsilon}=\det \tilde a\, {w_{\epsilon}}_{tt}\circ\tilde\eta^{-1}\ .$$ 
From (\ref{twice8}) and (\ref{weakcvwt}), we then get
\begin{align*}
\|{u_{\epsilon}}_t\|_{L^2(t-\delta t,t+\delta t; H^1_{0,\operatorname{div}} (\tilde\Omega_{2\delta}(t);{\mathbb R}^3)')}&\le \check C\\ 
\|{u_{\epsilon}}\|_{L^2(t-\delta t,t+\delta t; H^1_{0} (\tilde\Omega_{2\delta}(t);{\mathbb R}^3))}&\le \check C\ .
\end{align*}
Thus, from the classical compactness results (since the domain $\tilde\Omega_{2\delta}(t)$ is fixed on $(t-\delta t,t+\delta t)$), 
$$u_{\epsilon}\rightarrow u=\det\tilde a\, {w}_{tt}\circ
\tilde\eta^{-1}\ \text{in}\ L^2(t-\delta t,t+\delta t; 
L^2 (\tilde\Omega_{2\delta}(t);{\mathbb R}^3))\ ,$$
which obviously gives (since on $(t-\delta t,t+\delta t)$, $\tilde\eta^{-1}(t',\tilde\Omega_{2\delta}(t))\subset \Omega_{\delta}^f$), 
$$ {w_{\epsilon}}_{tt}\rightarrow  {w}_{tt}\ \text{in}\ L^2(t-\delta t,t+\delta t; L^2 (
\Omega^f_{\delta};{\mathbb R}^3))\ .$$
By a finite covering argument, and the $L^{\infty}$ bound (\ref{wepsilontt3}) of
${w_{\epsilon}}_{tt}$ in $L^2(\Omega;{\mathbb R}^3)$, we then infer that
\begin{align}
\text{limsup}\ \int_0^T \|{w_{\epsilon}}_{tt}\|^2_{L^2(\Omega^f_{\delta};{\mathbb R}^3)}\ dt  = \int_0^T \| {\tilde w}_{tt}\|^2_{L^2(\Omega^f_{\delta};{\mathbb R}^3)}\ dt  
\ .
\label{qtstrong3}
\end{align} 

Successively from the Cauchy-Schwarz and Sobolev inequalities,
\begin{align}
\|{w_{\epsilon}}_{tt}\|^2_{L^2(0,T; L^2(\Omega_0^f\cap {\Omega^f_{\delta}}^c;{\mathbb R}^3))}&\le C\ |\Omega_0^f\cap{\Omega^f_{\delta}}^c|\ \|\tilde w_{tt}\|^2_{L^2(0,T; H^1(\Omega_0^f\cap {\Omega^f_{\delta}}^c;{\mathbb R}^3))} \nonumber\\
&\le\  \tilde C\ |\Omega_0^f\cap{\Omega^f_{\delta}}^c|\ N(u_0,f)^2\ .
\label{qtstrong4}
\end{align}

From (\ref{qtstrong3}) and (\ref{qtstrong4}), we then infer
\begin{align*}
\text{limsup}\ \int_0^T \|{w_{\epsilon}}_{tt}\|^2_{L^2(\Omega_0^f;{\mathbb R}^3)}\ dt  \le \int_0^T \| {\tilde w}_{tt}\|^2_{L^2(\Omega_0^f;{\mathbb R}^3)}\ dt +
\tilde C\ |\Omega_0^f\cap{\Omega^f_{\delta}}^c|\ N(u_0,f)^2\ , 
\end{align*} 
which immediately shows that
\begin{align*}
\text{limsup}\ \int_0^T \|{w_{\epsilon}}_{tt}\|^2_{L^2(\Omega_0^f;{\mathbb R}^3)}\ dt  \le \int_0^T \| {\tilde w}_{tt}\|^2_{L^2(\Omega_0^f;{\mathbb R}^3)}\ dt\ , \end{align*}
thus establishing the strong convergence as $\epsilon=\frac{1}{m_l}\rightarrow 0$,
\begin{align}
{w_{\frac{1}{m_l}}}_{tt}\rightarrow  \tilde w_{tt}\ \text{in}\ L^2(0,T; L^2(\Omega_0^f;{\mathbb R}^3))\  .
\label{qtstrong5}
\end{align} 
Next, we restrict our test function $\phi$ to be in the space
$\{ \phi \in {\mathcal V}_{\tilde v}(t) \ | \ \phi = 0 \ \text{ on } \
\overline{\Omega^s_0}\}$.  For all such test functions and
for a.e $t\in (0,T)$, 
 $\phi \in {\mathcal V}_{\tilde v}(t)$,	
 \begin{align*}
& ( \tilde w_{tt}(t)-{w_{\epsilon}}_{tt}(t), \phi)_{L^2(\Omega;{\mathbb R}^3)}  
+\nu ((\tilde a_k^r \tilde a_k^s (\tilde w-{w_{\epsilon}},_r))_t(t),\  \phi,_s)_{L^2(\Omega_0^f;{\mathbb R}^3)} \nonumber\\
& - ((\tilde a_i^j)_t (\tilde q-q_{\epsilon})(t),\  \phi^i,_j)_{L^2(\Omega_0^f;{\mathbb R})} = 0 \ ;
\end{align*}
thus, the second Lagrange multiplier Lemma \ref{Lagrangebis} ensures us, 
from the strong convergence in (\ref{qtstrong5}), (\ref{strongcvq}) and 
(\ref{strongcvwt}),that
\begin{equation}
\label{qtstrong6}
\| ({{\bar q}_{\frac{1}{m_l}}})_t- \bar q_t\|_{L^2(0,T; L^2(\Omega_0^f;{\mathbb R}))} \rightarrow 0\ ,
\end{equation}
where $\displaystyle\bar q_t
=\tilde q_t-\frac{1}{|\Omega_0^f|} 
\int_{\Omega_0^f} \tilde q_t\ \det\nabla\tilde\eta$, and a similar
definition for $\bar q_{\epsilon}$. In the following, we will denote
$\displaystyle c=\frac{1}{|\Omega_0^f|} \int_{\Omega_0^f} 
\tilde q_t\ \det\nabla\tilde\eta$ and 
$\displaystyle c_{\epsilon}=\frac{1}{|\Omega_0^f|} 
\int_{\Omega_0^f} {q_{\epsilon}}_t\ \det\nabla\tilde\eta$.

\noindent{\bf Step 6. An inequality for $\tilde w_{tt}$ with a 
constant independent of the mollification parameter.}

Now, from the weak convergence (\ref{weakcvwt}) and the compactness of the trace operator, we then infer that as $\epsilon=\frac{1}{m_l} \rightarrow 0$,
\begin{equation}
\label{qtstrong7}
\|(\tilde w-w_{\epsilon})_{tt}\|^2_{L^2(0,T;L^2(\Gamma_0;{\mathbb R}^3))}\rightarrow 0\ .
\end{equation}

 We  now note that from (\ref{wepsilontt}) and (\ref{qepsilontt}), for any $0< t< T$, and $0<\delta t<\text{Min} (t,\ T-t)$,
\begin{align}
&  \frac{1}{2}\   \int_{t-\delta t}^{t+\delta t} \|{w_{\epsilon}}_{tt}\|^2_{L^2(\Omega;{\mathbb R}^3)} 
+{\nu} \int_{t-\delta t}^{t+\delta t} \int_0^{t'} (\tilde a_k^r {w_{\epsilon}}_{tt},_r,\  \tilde a_k^s {w_{\epsilon}}_{tt},_s)_{L^2(\Omega_0^f;{\mathbb R}^3)} 
\nonumber\\
&+ \frac{1}{2}  \int_{t-\delta t}^{t+\delta t} (c^{ijkl}   {{w_{\epsilon}}_{t}}^k,_l , {{w_{\epsilon}}_{t}}^i,_j )_{L^2(\Omega_0^s;{\mathbb R})}  
 - \int_{t-\delta t}^{t+\delta t} \int_0^{t'}\int_{\Omega_0^f} {q_{\epsilon}}_{t} [ 2 (\tilde a_i^j)_t 
 {w_{\epsilon}}_t^i,_j + (\tilde a_i^j)_{tt} 
 {w_{\epsilon}}^i,_j]_t \nonumber\\
&+\int_{t-\delta t}^{t+\delta t} \int_{\Omega_0^f} {q_{\epsilon}}_{t} \ [2\ (\tilde a_i^j)_t 
 {w_{\epsilon}}_t^i,_j + (\tilde a_i^j)_{tt} 
 {w_{\epsilon}}^i,_j]-2 \int_{t-\delta t}^{t+\delta t} \int_0^{t'}\int_{\Omega_0^f} (\tilde a_i^j)_t {q_{\epsilon}}_t{w_{\epsilon}}_{tt}^i,_j\nonumber\\
& -\int_{t-\delta t}^{t+\delta t}  \int_0^{t'}\int_{\Omega_0^f} (\tilde a_i^j)_{tt} {q_{\epsilon}}{w_{\epsilon}}_{tt}^i,_j + {\nu} \int_{t-\delta t}^{t+\delta t} \int_0^{t'} ((\tilde a_k^r \tilde a_k^s)_{tt} {w_{\epsilon}},_r,\  {w_{\epsilon}}_{tt},_s)_{L^2(\Omega_0^f;{\mathbb R}^3)}\nonumber\\
& +2\ {\nu} \int_{t-\delta t}^{t+\delta t} \int_0^{t'} ((\tilde a_k^r \tilde a_k^s)_{t} {{w_{\epsilon}}_t},_r,\  {w_{\epsilon}}_{tt},_s)_{L^2(\Omega_0^f;{\mathbb R}^3)}\nonumber\\
&\le C\  \delta t\ N(u_0,f)^2+ \int_{t-\delta t}^{t+\delta t} \int_0^{t'} (F_{tt}, {w_{\epsilon}}_{tt})_{L^2 (\Omega_0^f;{\mathbb R}^3)}+  \int_{t-\delta t}^{t+\delta t} \int_0^{t'} (f_{tt},{w_{\epsilon}}_{tt})_{L^2 (\Omega_0^s;{\mathbb R}^3)} .
\label{energywtt1}
\end{align}

The first three terms of the left-hand side of this inequality will be dealt with by
weak lower semi-continuity. By the weak convergence in (\ref{weakcvwt}) and the 
strong convergence in (\ref{wstrong}), (\ref{strongcvwt}), (\ref{strongwt3}), 
and (\ref{strongcvq}), we infer that all of the remaining terms, other than
the term 
$\int_{t-\delta t}^{t+\delta t}\int_0^{t'}\int_{\Omega_0^f} (\tilde a_i^j)_t 
{q_{\epsilon}}_t{w_{\epsilon}}_{tt}^i,_j\ $,
converge as $\epsilon\rightarrow 0$
to the same expressions with the limits $\tilde w$ and $\tilde q$ replacing
$w_\epsilon$ and $q_\epsilon$. 
From the definitions of $c$ and $c_{\epsilon}$, we have that
\begin{align*}
\int_{t-\delta t}^{t+\delta t}\int_0^{t'}\int_{\Omega_0^f} (\tilde a_i^j)_t 
{q_{\epsilon}}_t{w_{\epsilon}}_{tt}^i,_j\ 
=&\int_{t-\delta t}^{t+\delta t}\int_0^{t'}
\left[ \int_{\Omega_0^f} (\tilde a_i^j)_t 
({\bar q_{\epsilon}})_t{w_{\epsilon}}_{tt}^i,_j
+c_{\epsilon}\int_{\Omega_0^f} (\tilde a_i^j)_t {w_{\epsilon}}_{tt}^i,_j
\right]\ .
\end{align*}
From the strong convergence (\ref{qtstrong6}) and the weak convergence (\ref{weakcvwt}), we then deduce that the first term of the right-hand side of this
inequality converges as $\epsilon=\frac{1}{m_l}\rightarrow 0$ to the corresponding term where $\bar{q}_t$ replaces ${q_{\epsilon}}_t$ and $\tilde w_{tt}$ 
replaces ${w_\epsilon}_{tt}$.

For the second term of this right-hand side, we notice from a spatial 
integration by parts (since $c_{\epsilon}$ depends only on the time variable) 
that
\begin{align*}
\int_{t-\delta t}^{t+\delta t}\int_0^{t'}c_{\epsilon} \int_{\Omega_0^f} 
(\tilde a_i^j)_t {w_{\epsilon}}_{tt}^i,_j
= &\ \int_{t-\delta t}^{t+\delta t}\int_0^{t'}c_{\epsilon} 
\left[ -\int_{\Omega_0^f} ((\tilde a_i^j)_t),_j {w_{\epsilon}}_{tt}^i
 +\int_{\Gamma_0} (\tilde a_i^j)_t {w_{\epsilon}}_{tt}^i N_j \right]
%\label{wtt2}
\end{align*}
and thus from the weak convergence in (\ref{weakcvwt}) and the strong 
convergence in
(\ref{qtstrong5}) and (\ref{qtstrong7}) we then get the convergence as
$\epsilon=\frac{1}{m_l}\rightarrow 0$ to the corresponding term where $c$ 
replaces 
${c_{\epsilon}}$ and $\tilde w_{tt}$ replaces 
${w_\epsilon}_{tt}$. This implies 
that as $\epsilon=\frac{1}{m_l}\rightarrow 0$,
$$
\int_{t-\delta t}^{t+\delta t}\int_0^{t'}\int_{\Omega_0^f} (\tilde a_i^j)_t {q_{\epsilon}}_t{w_{\epsilon}}_{tt}^i,_j\rightarrow 
\int_{t-\delta t}^{t+\delta t}\int_0^{t'}\int_{\Omega_0^f} (\tilde a_i^j)_t 
\tilde {q}_t{\tilde w}_{tt}^i,_j\ .$$
Consequently, all the terms, except the three first ones, appearing in the
inequality (\ref{energywtt1}) are convergent as $\epsilon=\frac{1}{m_l}
\rightarrow 0$ to the same expression, where $q_{\epsilon}$ and $w_{\epsilon}$ 
are replaced respectively by $\tilde q$ and $\tilde w$. 

By weak lower semi-continuity for the first three integrals, we then infer that as $\epsilon=\frac{1}{m_l}\rightarrow 0$ the same inequality as the previous one holds with $\tilde w$ and $\tilde q$ replacing
respectively $w_{\epsilon}$ and $q_{\epsilon}$. By dividing those integrals by
$2\delta t$ and passing to the limit as $\delta t\rightarrow 0$ (which is possible a.e. in $(0,T)$), we then get that a.e. in $(0,T)$,
\begin{align}
&  \frac{1}{2}\  \|\tilde w_{tt} (t)\|^2_{L^2(\Omega;{\mathbb R}^3)} 
+{\nu} \int_0^{t} (\tilde a_k^r {\tilde w}_{tt},_r,\  \tilde a_k^s {\tilde w}_{tt},_s)_{L^2(\Omega_0^f;{\mathbb R}^3)} 
\nonumber\\
&+ \frac{1}{2}  (c^{ijkl}   {{\tilde w}_{t}}^k,_l (t), {{\tilde w}_{t} }^i,_j (t))_{L^2(\Omega_0^s;{\mathbb R})}- \int_0^{t}\int_{\Omega_0^f} {\tilde q}_{t}\ [2 (\tilde a_i^j)_{t} 
 {\tilde w}_t^i,_j + (\tilde a_i^j)_{tt} 
 {\tilde w}^i,_j]_t \nonumber\\
&+ \int_{\Omega_0^f} {\tilde q}_{t} (t)\ [ 2(\tilde a_i^j)_t 
 {\tilde w}_t^i,_j +(\tilde a_i^j)_{tt} {\tilde w}^i,_j](t)
 - 2\int_0^t \int_{\Omega_0^f} (\tilde a_i^j)_t \tilde q_t 
 {\tilde w}_{tt}^i,_j 
- \int_0^{t}\int_{\Omega_0^f} (\tilde a_i^j)_{tt} {\tilde q}{\tilde w}_{tt}^i,_j \nonumber\\
&+ {\nu} \int_0^{t} ((\tilde a_k^r \tilde a_k^s)_{tt} {\tilde w},_r,\  {\tilde w}_{tt},_s)_{L^2(\Omega_0^f;{\mathbb R}^3)} +2\ {\nu} \int_0^{t} ((\tilde a_k^r \tilde a_k^s)_{t} {{\tilde w}_t},_r,\  {\tilde w}_{tt},_s)_{L^2(\Omega_0^f;{\mathbb R}^3)}\nonumber\\
&\qquad\qquad\qquad\le C\  N(u_0,f)^2 + \int_0^{t}( F_{tt}, {\tilde w}_{tt})_{L^2(\Omega_0^f;{\mathbb R}^3)} 
+  \int_0^{t} (f_{tt},{\tilde w}_{tt})_{L^2 (\Omega_0^s;{\mathbb R}^3)} ,
\label{energywtt}
\end{align}
where (we recall) $C$ does not depend on the smoothing parameter of $\tilde a$. 
\end{proof}

%We will use the classical results, proved in \cite{Evans1998} and
%\cite{CaCa1995}:

%\begin{lemma}\label{quotients}
%Suppose that $w \in{L^2({\mathbb R}^3_+; {\mathbb R}^3)}$.

%\noindent
%{(i)} If 
%$$\|D_h w\|_{L^2( {\mathbb R}^3_+; {\mathbb R}^3)}  \le M$$
%for a constant $M$ and for $h=t e_{\alpha}$ ($t\in \mathbb R-\{0\}$, $\alpha\in\{1,2\}$), then 
%$$\|w,_{\alpha}\|_{L^2( {\mathbb R}^3_+; {\mathbb R}^3)}  \le  M\,,$$
%and $\|D_h w\|_{L^2( {\mathbb R}^3_+; {\mathbb R}^3)}  \le C
%\|w,_{\alpha}\|_{L^2( {\mathbb R}^3_+; {\mathbb R}^6)}$ for some constant
%$C$. The reciprocal statement is also true.

%\noindent
%{(ii)} If 
%$$\|D_{-h}D_{h} w\|_{L^2( {\mathbb R}^3_+; {\mathbb R}^3)}  \le M$$
%for a constant $M$ and for $h=t e_{\alpha_i}, t\in {\mathbb R}-\{0\}, \alpha_i\in \{1,2\} $ and $h=t \frac{e_1+e_2}{\sqrt{2}}, t\in {\mathbb R}-\{0\}$, then 
%$$\| \nabla_0\nabla_0 w\|_{L^2( {\mathbb R}^3_+; {\mathbb R}^{9})}  \le M\,,$$
%and $\|D_{-h}D_{h} w\|_{L^2( {\mathbb R}^3_+; {\mathbb R}^3)}  \le
%C\ \| \nabla_0\nabla_0 w\|_{L^2( {\mathbb R}^3_+; {\mathbb R}^{9})}$. The
%reciprocal statement is also true.

% Similar statements are of course
%true for ${\mathbb R}^3_-$ as well.
%\end{lemma}

\subsection{Regularity for $\tilde w$ and its first and second time derivatives, 
dependent on the regularization parameter of $\tilde a$}

As discussed in the introduction, we shall focus on the regularity near the 
interface, which will provide us with the trace estimates on the interface. 
Elliptic regularity for the  Dirichlet problems will then yield the full regularity
result in each interior component.  In this subsection, $\tilde C$ continues to
denote a generic constant which depends on the same variables as $C$ and $C(M)$, 
and additionally on the regularization parameter. In Section \ref{9}, we obtain 
estimates independent of $n$, by interpolation mainly, which requires us to know 
{\it a priori} that
the solution is smooth (without using the estimates that we get in this subsection,
since they blow up with the regularization parameter).

Recall that we have
already shown that  $\tilde w\in L^2(0,T; {\mathcal V}_{\tilde v}( \cdot)) $, $\tilde w_t\in L^2(0,T; H^1_0(\Omega;{\mathbb R}^3))$ and $\tilde w_{tt} \in 
 {\mathcal W}([0,T])$, and that both $\tilde q$ and $\tilde q_t$ are in
$L^2(0,T; L^2(\Omega^f_0; {\mathbb R}))$. 

The missing regularity results will be recovered using difference quotients. Recall that if we consider the partition of the space ${\mathbb R}^3$ formed by the two half-spaces ${\mathbb R}^3_+ := \{(x^1,x^2,x^3)
\in{\mathbb R}^3 \ | \ x^3>0\}$ and ${\mathbb R}^3_- := \{(x^1,x^2,x^3)
\in{\mathbb R}^3 \ | \ x^3<0\}$  and the horizontal plane $\{x^3=0\}$ with the usual orthonormal 
basis $(e_1,e_2,e_3)$, then we have  

\begin{definition}\label{differences}
The first-order difference quotient of a function $u$ of size $h$ at $x$ is given by
$$
D_h u(x) = \frac{u(x+h)-u(x)}{|h|},
$$
where $h$ is any vector orthogonal to $e_3$.  The second-order difference
quotient of $u$ of size $h$ is defined as $D_{-h}D_h u(x)$, given explicitly by
$$D_{-h}D_hu(x)=\frac{u(x+h)+u(x-h)-2u(x)}{|h|^2}\,.$$
We will denote $$\nabla_0 u=(u,_1,\ u,_2).$$
\end{definition}

Letting $A^{ij}= \tilde a^i_k \tilde a^j_k$, we write the weak form as 

\begin{align*}
&(\tilde w_t, \phi) _{L^2(\Omega;{\mathbb R}^3)}  + \nu (A^{ij} \tilde w,_i, \phi,_j)
_{L^2(\Omega_0^f; {\mathbb R}^3)} 
+(c^{ijkl} \int_0^t \tilde w^k,_l , \phi^i,_j)_ {L^2(\Omega_0^s; {\mathbb R})} \\
& \qquad \qquad
- (\tilde q, \tilde a^l_k \phi^k,_l)_{L^2(\Omega_0^f; {\mathbb R})}  =
(F, \phi)_{L^2(\Omega_0^f; {\mathbb R}^3)} 
+ ( f, \phi)_ {L^2(\Omega_0^s; {\mathbb R}^3)} 
\end{align*}
for all $\phi \in H^1_0(\Omega;{\mathbb R}^3)$, for a.e. $0\le t \le T$. 

Next, assume that $\Omega=B(0,1)$, the unit ball centered at $0$, and that
$\Omega^f_0= \{ x \in B(0,1) \ | \ x^3>0\}$ and
$\Omega^s_0= \{ x \in B(0,1) \ | \ x^3<0\}$.
Select a smooth cutoff function $\zeta$ satisfying
$$
\zeta =1 \ \ \text{on} \ \ B(0,\frac{1}{2}), \ \ \zeta=0 
\ \ \text{on} \ \  {\mathbb R}^3 - B(0,1), 
 \ \ \text{and} \ \ 0 \le \zeta \le 1.
$$
Let $\phi = D_{-h} (\zeta^2 D_h \tilde w)$; then clearly $\phi \in H^1_0
(\Omega; {\mathbb R}^3)$ for a.e. $t\in [0,T]$.  We
may thus substitute  $\phi$ into the above weak form to obtain
$$
A_1+A_2+A_3-A_4 =B,
$$
where
\begin{align*}
A_1&= (D_h \tilde w_t, \zeta^2 D_h \tilde w)_{L^2(\Omega; {\mathbb R}^3)}, \\
A_2&= \nu(D_h (A^{ij} \tilde w,_i), (\zeta^2 D_h \tilde w),_j)_ {L^2(\Omega_0^f; {\mathbb R}^3)},\\
A_3&=  (D_h (c^{ijkl} \int_0^t \tilde w^k,_l), (\zeta^2 D_h \tilde w^i),_j)_
{L^2(\Omega_0^s; {\mathbb R})} , \\
A_4&=  (q, \tilde a^l_k (D_{-h} [ \zeta^2 D_h \tilde w^k]),_l)_{L^2(\Omega_0^f; {\mathbb R})},\\ 
B& =( F, D_{-h}(\zeta^2D_h\tilde w))_{L^2(\Omega_0^f; {\mathbb R}^3)}+ 
( D_h f, \zeta^2D_h\tilde w)_{L^2(\Omega_0^s; {\mathbb R}^3)}.
\end{align*}

For the first two terms, we easily find that
\begin{align*}
A_1&=\frac{1}{2} \frac{d}{dt}  \| \zeta D_h\tilde w\|^2_{L^2(\Omega; {\mathbb R}^3)},\\
A_2&\ge C\| \zeta D_h\nabla \tilde w\|^2_ {L^2(\Omega_0^f; {\mathbb R}^9)} -
\tilde C \| \nabla \tilde w \|^2_ {L^2(\Omega_0^f; {\mathbb R}^9)}. 
\end{align*}

For the remaining terms, we shall use the notation $C^h(x)$ to denote
$C(x+h)$. Whereas the coefficients of the elasticity tensor are constant, one should keep in mind that our assumption on the domain comes in fact from a change of variables which produces a non constant elasticity tensor. It is the integral
below with the $D_h C$ term which necessitates the hyperbolic scaling of our 
functional framework.

Expanding $A_3$, we have that
\begin{align*}
&A_3 = (\zeta^2 c^h: D_h \int_0^t\nabla \tilde w, D_h \nabla \tilde w)_
 {L^2(\Omega_0^s; {\mathbb R}^9)} +  (\zeta^2 D_h c : \int_0^t\nabla \tilde w,
D_h \nabla \tilde w)_ {L^2(\Omega_0^s; {\mathbb R}^9)} \\
& + (2 \zeta \nabla \zeta \otimes D_h \tilde w, c^h : D_h \int_0^t \nabla \tilde w)_
 {L^2(\Omega_0^s; {\mathbb R}^9)}
  + (2 \zeta \zeta,_j D_hc^{ijkl} \int_0^t
\tilde w^k,_l, D_h \tilde w^i)_{L^2(\Omega_0^s; {\mathbb R})}. 
\end{align*}
The second term on the right-hand-side is
\begin{align*}
& (D_{-h} (\zeta^2 D_h C : \int_0^t\nabla \tilde w), \nabla \tilde w )_
{L^2(\Omega_0^s; {\mathbb R}^9)}  \\
& \qquad = ( [\zeta^2 D_h C]^{-h} : D_{-h} \int_0^t\nabla \tilde w +
D_{-h} (\zeta^2 D_h C) : \int_0^t \nabla \tilde w, \nabla \tilde w)_
 {L^2(\Omega_0^s; {\mathbb R}^9)}. 
\end{align*}
Hence for $\theta>0$, we see that the $A_3$ term yields the following 
inequality:
\begin{align*}
& \bigr|A_3-\frac{1}{2}\frac{d}{dt} ( \zeta^2 C^h : \int_0^t D_h \nabla \tilde w,
\int_0^t D_h\nabla \tilde w )_ {L^2(\Omega_0^s; {\mathbb R}^9)}\bigl|\\
&\qquad\qquad\le 
\theta \| \zeta D_{-h} \int_0^t\nabla \tilde w\|^2_ {L^2(\Omega_0^s; {\mathbb R}^9)}+ C_\theta\| \nabla \tilde w\|^2_ {L^2(\Omega_0^s; {\mathbb R}^9)}
+ C\|\int_0^t \nabla \tilde w\|^2_ {L^2(\Omega_0^s; {\mathbb R}^9)}  .
\end{align*}
For the $A_4$ term, we have that
\begin{align*}
A_4 &= (\tilde q\ ,\ \tilde a^l_k [\zeta^2]^{-h} D_{-h} D_h \tilde w^k,_l + \tilde a^l_k (D_{-h}\zeta^2)
D_h \tilde w^k,_l + [2 \zeta \zeta,_l]^{-h} D_{-h}D_h \tilde w^k \tilde a^l_k \\
& \qquad + 2\tilde a^l_k D_{-h}(\zeta \zeta,_l) D_h \tilde w^k)_
{L^2(\Omega_0^f; {\mathbb R})}.
\end{align*}
By the divergence-free  condition, 
$\tilde w \in\mathcal V_{\tilde v}([0,T])$, we get in $\Omega_0^f$
\begin{align*}
&0= D_{-h} ( [\tilde a^l_k]^h D_h \tilde w^k,_l + D_h \tilde a^l_k \tilde w^k,_l] \\
& \qquad =
\tilde a^l_k D_{-h}D_h \tilde w^k,_l + D_{-h} [\tilde a^l_k]^hD_h \tilde w^k,_l +
[D_h \tilde a^l_k]^{-h}D_{-h}\tilde w^k,_l
+ D_{-h}D_h \tilde a^l_k \tilde w^k,_l,
\end{align*}
allowing us to eliminate the first term appearing in the expression of $A_4$, which gives for $\theta>0$,
$$
|A_4| \le \theta\| \zeta D_h \tilde w^k,_l\|^2_{L^2(\Omega_0^f; {\mathbb R}^3)} +
C_\theta \| \tilde q\|^2_ {L^2(\Omega_0^f; {\mathbb R})} + \tilde C \| \nabla \tilde w\|^2_
{L^2(\Omega_0^f; {\mathbb R}^9)}. 
$$
Finally,
\begin{align*}
|B| \le &\ \theta [\ \| \nabla \tilde w\|^2_ {L^2(\Omega_0^s; {\mathbb R}^9)}+\| \zeta D_h \nabla \tilde w\|^2_ {L^2(\Omega_0^f; {\mathbb R}^9)}\ ]\\ 
&+C \| \nabla \tilde w\|^2_ {L^2(\Omega_0^f; {\mathbb R}^9)} + C_{\theta} [\ \|f\|^2_ 
{H^1(\Omega_0^s; {\mathbb R}^3)}+\|f\|^2_ 
{L^2(\Omega_0^f; {\mathbb R}^3)}\ ].
\end{align*}

Choosing $\theta>0$ sufficiently small, we have the inequality
\begin{align*}
&\frac{d}{dt} \left(
\| \zeta D_h \tilde w \|^2_{L^2(\Omega; {\mathbb R}^3)}
+ ( \zeta^2 C^h: \int_0^t D_h \nabla \tilde w, \int_0^t D_h \nabla \tilde w)
_ {L^2(\Omega_0^s; {\mathbb R}^9)}\right) 
+ \|\zeta D_h\nabla \tilde w\|^2_{L^2(\Omega_0^f; {\mathbb R}^9)} \\
&\le \tilde C \left(
\|\zeta \int_0^t D_h\nabla \tilde w \|^2_ {L^2(\Omega_0^s; {\mathbb R}^9)}
+ \| \int_0^t \nabla \tilde w\|^2_ {L^2(\Omega_0^s; {\mathbb R}^9)} 
+ \|  \nabla \tilde w\|^2_ {L^2(\Omega_0^s; {\mathbb R}^9)}
  \right.\\
&\left.
+\|\nabla \tilde w\|^2_ {L^2(\Omega_0^f; {\mathbb R}^9)} 
+\|\tilde q\|^2_{L^2(\Omega_0^f; {\mathbb R})} + \|f\|^2_ {L^2(\Omega_0^f; {\mathbb R}^3)}+ \|f\|^2_ {H^1(\Omega_0^s; {\mathbb R}^3)} 
\right). 
\end{align*}
From Gronwall's inequality, it follows that
$\partial_ \alpha \partial_j \tilde w \in L^2(0,T; {L^2(V^f; {\mathbb R}^3)})$
where $V^f = \{ x \in B(0,\frac{1}{2})\ | \ x^3 >0 \}$, and where $\alpha=1,2$
and $j=1,2,3$.  Hence, 
$\partial_ \alpha \tilde w \in L^2(0,T; {H^1(V^f; {\mathbb R}^3)})$, so that by
the trace theorem we obtain 
$\partial_ \alpha \tilde w \in L^2(0,T; {H^{0.5}(V^f \cap \{ x^3=0\};{\mathbb R}^3)})$. Thus, \begin{equation}\label{boundaryf}
\tilde w\in { L^2(0,T; {H^{1.5}(V^f \cap \{ x^3=0\};{\mathbb R}^3)})} .
\end{equation}
(with an estimate which blows-up as the mollification parameter 
$n\rightarrow\infty$).  Similarly, 
\begin{equation}\label{boundarys}
\int_0^t \tilde w \in L^\infty(0,T; {H^{1.5}(V^s \cap \{ x^3=0\};{\mathbb R}^3)})\ ,
\end{equation}
where $V^s = \{ x \in B(0,\frac{1}{2})\ | \ x^3 <0 \}$.

We now drop the assumption that $\Omega$ is the unit ball, and once again 
assume it is an open bounded subset of ${\mathbb R}^3$ with all of the smoothness
assumption stated previously.  We choose any point $x_0 \in \Gamma_0$ and
assume that
\begin{align*}
\Omega^f_0 \cap B(x_0,r) &= \{ x \in B(x_0,r) \ | \ x^3 > \gamma(x^1,x^2) \} \\
\Omega^s_0 \cap B(x_0,r) &= \{ x \in B(x_0,r) \ | \ x^3 < \gamma(x^1,x^2) \} \\
\end{align*}
for some $r>0$ and some smooth function $\gamma:{\mathbb R}^2 \rightarrow
{\mathbb R}$.   We define the following change of variables:
\begin{align*}
y^i &= x^i =: \Phi^i(x), \ \ \ \ i=1,2 \\
y^3 &= x^3 - \gamma(x^1,x^2) = \Phi^3(x),
\end{align*}
and write
$$
y=\Phi(x).
$$
Similarly, we set
\begin{align*}
x^i &= y^i =: \Psi^i(y), \ \ \ \ i=1,2 \\
x^3 &= y^3 + \gamma(y^1,y^2) = \Psi^3(y),
\end{align*}
and write
$$
x=\Psi(y).
$$
Then $\Phi = \Psi^{-1}$, and the mapping $x \mapsto \Phi(x)=y$ {\it straightens
out} $\Gamma_0$ near $x_0$, and $\det \Phi= \det \Psi =1$. 

We can assume $0= \Phi(x_0)$.  Choose $s>0$ so small that $B(0,s) \subset
\Phi(B(x_0,r))$.  Let
\begin{align*}
w'(t,y) = \tilde w (t,\Phi(y)), \ \ 
q'(t,y) = \tilde q (t,\Phi(y)), \ \ 
f'(t,y) = f (t,\Phi(y)), 
\end{align*}
Then $w'$ and $w'_t$ are in  $L^2(0,T;H^1(\Omega';{\mathbb R}^3))$, where  $
\Omega'= B(0,s)$.  We also set
\begin{align*}
{\Omega^f_0}' = B(0,s) \cap \{y^3>0\}, \ \ 
{\Omega^s_0}' = B(0,s) \cap \{y^3<0\}, 
\end{align*}

Then, since $(\tilde w,\tilde q)$ satisfy the weak formulation, applying the change of 
variables, we see that $(w',q')$ satisfy
\begin{align}
&(w'_t, \phi') _{L^2(\Omega';{\mathbb R}^3)}  + \nu (a'^{ij} w',_i, \phi',_j)
_{L^2({\Omega_0^f}'; {\mathbb R}^3)} 
+(c'^{ijkl} \int_0^t w'^k,_l , \phi'^i,_j)_ {L^2({\Omega_0^s}'; {\mathbb R}^3)} \nonumber\\
& \qquad \qquad
- (q', [a^l_k \circ \Psi]g^r_l \phi'^k,_r)_{L^2({\Omega_0^f}'; {\mathbb R}^3)} 
=
(f', \phi')_{L^2(\Omega'; {\mathbb R}^3)} 
\label{weakchange}
\end{align}
for all $\phi' \in H^1_0(\Omega';{\mathbb R}^3)$, for a.e. $0\le t \le T$,
where
\begin{align*}
a'^{kl} = A^{ij} \circ \Psi\ g^k_i g^l_j, \ \ \ 
c'^{irks} = c^{ijkl} \ g^s_lg^r_j , \ \ \ 
g(y) = [\nabla \Psi(y)] ^{-1} .
\end{align*}
It is easy to verify that both $a'$ and $c'$ retain the uniform ellipticity
conditions of the original operators $A$ and $C$; Moreover, $w'$ satisfies the
divergence condition $a_i^j\circ\Phi w'^i,_j=0$ in $[0,T]\times \Omega'$. Thus we may apply the results
obtained above for the case that the domain is the unit ball to find that
\begin{align*}
w' &\in L^2(0,T; {H^{1.5}({V^f}' \cap \{ x^3=0\};{\mathbb R}^3)}),\\ 
\int_0^t w' &\in L^\infty(0,T; {H^{1.5}({V^s}' \cap \{ x^3=0\};{\mathbb R}^3)}),
\end{align*}
where ${V^f}' = \{ x \in B(0,\frac{s}{2})\ | \ x^3 >0 \}$ and
${V^s}' = \{ x \in B(0,\frac{s}{2})\ | \ x^3 <0 \}$.
Consequently,
\begin{equation}
\label{boundaryfs}
\tilde w\in { L^2(0,T; {H^{1.5}(\partial V^f \cap \Gamma_0;{\mathbb R}^3)})} 
, \ 
\int_0^t \tilde w\in {L^\infty(0,T; {H^{1.5}(\partial V^s \cap \Gamma_0;{\mathbb R}^3)})}, 
\end{equation}
where ${V^f} = \Psi({V^f}')$ and ${V^s} = \Psi({V^s}')$.

Since $\Gamma_0$ is compact, we can as usual cover $\Gamma_0$ with finitely
many sets of the type used above.  Summing the resulting estimates, we find
that we have for the trace on $\Gamma_0$
\begin{subequations}
\label{Boundaryfs}
\begin{gather}
\tilde w\in { L^2(0,T; {H^{1.5}(\Gamma_0;{\mathbb R}^3)})} 
, \label{Boundaryfs.a}\\
\int_0^t \tilde w\in {L^\infty(0,T; {H^{1.5}(\Gamma_0;{\mathbb R}^3)})}
. \label{Boundaryfs.b}
\end{gather}
\end{subequations}

Converting the fluid equations into Eulerian variables by composing with 
$\tilde \eta^{-1}$, we obtain a Stokes problem in the domain 
$\tilde \eta(\Omega_0^f)$:
\begin{subequations}
\label{ssStokes}
\begin{align}
-\nu  \triangle u + \nabla p &=  f-\tilde w_t\circ\tilde\eta^{-1}+\nu \tilde a_l^j,_j\circ\eta^{-1} u,_l-p\ (\tilde a_i^j),_j\circ\tilde\eta^{-1} , \\
 \operatorname{div} u &=0, 
\end{align}
\end{subequations}
with the boundary conditions that $u=0$ on $\tilde \eta(\partial\Omega)$
and that $u \in { L^2(0,T; {H^{1.5}(\tilde \eta(\Gamma_0);{\mathbb R}^3)})} $, where
$u = \tilde w \circ \tilde \eta^{-1}$ and $p = \tilde q \circ \tilde \eta^{-1}$.  Since the domain $\tilde\eta(\Omega_0^f)$ is of class $H^3$, 
by the elliptic regularity of \cite{Eben}, 
(\ref{Boundaryfs.a}) implies that
$u \in L^2(0,T; {H^2(\tilde\eta(\Omega_0^f);{\mathbb R}^3)})$ and 
$p \in L^2(0,T; {H^1(\tilde\eta(\Omega_0^f);{\mathbb R})})$.  It follows that
\begin{equation}
\label{H2_reg_f}
\tilde w\in { L^2(0,T; {H^{2}(\Omega^f_0;{\mathbb R}^3)})} 
, \ \ \ 
\tilde q\in {L^2(0,T; {H^{1}(\Omega_0^f;{\mathbb R})})}
. 
\end{equation}
Similarly, elliptic regularity of the elasticity problem shows that
\begin{equation}\label{H2_reg_s}
\int_0^t \tilde w\in L^\infty(0,T; {H^{2}(\Omega^s_0;{\mathbb R}^3)}) 
.
\end{equation}

Next, we consider the weak form for the time derivate $\tilde w_t$ for all $\phi \in H^1_0(\Omega;{\mathbb R}^3)$:
\begin{align*}
&
(\tilde w_{tt}, \phi)_{L^2(\Omega; {\mathbb R}^3)} + \nu ( [A^{rs} \tilde w,_r]_t,
\phi,_s)_{L^2(\Omega_0^f; {\mathbb R}^3)} 
+(c^{ijkl} \tilde w^k,_l \ , \ \phi^i,_j)_{L^2(\Omega_0^s; {\mathbb R})} \\
& 
- ([\tilde a^j_i \tilde q]_t, \phi^i,_j)_ {L^2(\Omega_0^f; {\mathbb R})} = 
(F_t, \phi)_{L^2(\Omega_0^f; {\mathbb R}^3)}+ (f_t, \phi)_{L^2(\Omega_0^s; {\mathbb R}^3)}
 \ \ \  a.e. \ t \in (0,T).
\end{align*}
Expanding the time derivative, we see that there are two additional terms
in the weak form given by $( A^{rs}_t \tilde w,_r , {\phi},_s)
_{L^2(\Omega_0^f; {\mathbb R}^3)}$ and 
$((\tilde a_t)^j_i q, {\phi_t}^i,_j )_ {L^2(\Omega_0^f; {\mathbb R}^3)}$.  These
additional terms are easy to handle, and by 
letting $\phi = D_{-h}( \zeta^2 D_h \tilde w_t)$, and following the identical 
procedure as above, since we also already know that $\tilde w_{tt}\in L^{\infty}(0,T;L^2(\Omega;{\mathbb R}^3))$, we find that
\begin{equation}
\label{H2_regt}
\tilde w_t\in { L^2(0,T; {H^{2}(\Omega^f_0;{\mathbb R}^3)})} 
, \  
\tilde q_t\in {L^2(0,T; {H^{1}(\Omega_0^f;{\mathbb R})})}
,  \  
\tilde w\in { L^\infty(0,T; {H^{2}(\Omega^s_0;{\mathbb R}^3)})} 
.  
\end{equation}

Because of the assumptions on the forcing and these estimates 
for $\tilde w_t$, we may improve the regularity results (\ref{H2_reg_f}) and
(\ref{H2_reg_s}).  We apply the identical procedure, but this time we use
$\phi= D_{-h}D_{h} ( \zeta^2 D_{-h}D_h \tilde w)$ as the test function.  We find
that
\begin{equation}
\label{H3_reg}
\tilde w\in { L^2(0,T; {H^{3}(\Omega^f_0;{\mathbb R}^3)})} 
, \ 
\tilde q\in {L^2(0,T; {H^{2}(\Omega_0^f;{\mathbb R})})}
, \ 
\int_0^t \tilde w\in { L^\infty(0,T; {H^{3}(\Omega^s_0;{\mathbb R}^3)})} 
.  
\end{equation}
Moreover, $\|(\tilde w, \tilde q)\|_{Z_T} \le \tilde C N(u_0,f)$, where
the constant $\tilde C \rightarrow \infty$ as the mollification parameter
$n \rightarrow \infty$.

In the following section, we will use a different form of (\ref{weakchange}). If we denote by $\zeta$ a smooth cut-off function, equal to $1$ in a neighborhood of $0$ contained in $\Omega'$ and $0$ outside $\Omega'$, and denote $W=\zeta^2 w'$, $Q=\zeta^2 q'$, $\tilde b_l^j=\tilde a_l^k\circ\Psi\ g^j_k$, $C^{ijkl}=c'^{ijkl}$, we then obtain for any $\varphi\in H^1({\mathbb R}^3;{\mathbb R}^3)$, 
\begin{align}
&(W_t,\ \varphi) _{L^2({\mathbb R}^3;{\mathbb R}^3)}  + \nu\ (\tilde b_l^j \tilde b_l^k W,_k,\ \varphi,_j)
_{L^2({\mathbb R}^3_+; {\mathbb R}^3)} 
+(C^{ijkl} \int_0^t W^k,_l ,\ \varphi^i,_j)_ {L^2({\mathbb R}^3_-; {\mathbb R}^3)} \nonumber\\
& - (Q,\tilde b_k^r \varphi^k,_r)_{L^2({\mathbb R}^3_+; {\mathbb R})}
\nonumber\\ 
&=
(F_1,\ \varphi)_{L^2({\mathbb R}^3_+; {\mathbb R}^3)}+ 
(H_j,\ \varphi,_j)_{L^2({\mathbb R}^3_+; {\mathbb R}^3)}
+ (F_2,\ \varphi)_{L^2({\mathbb R}^3_-; {\mathbb R}^3)} +
(K_j,\ \varphi,_j)_{L^2({\mathbb R}^3_+; {\mathbb R}^3)}\ ,
\label{variationalr3} 
\end{align}
where
\begin{subequations}
\begin{align}
F_1^i&=\zeta^2 F'^i-\nu\ \zeta^2,_j \tilde b_l^j \tilde b_l^k w'^i,_k  + q' \tilde b_i^r \zeta^2,_r,\\
H_j&=\nu\ \tilde b_l^j \tilde b_l^k \zeta^2,_k w',\\
F_2^i&= \zeta^2 f'^i -C^{ijkl} \zeta^2,_j \int_0^t w'^k,_l,\\
K^i_j&=C^{ijkl} \zeta^2,_l \int_0^t w'^k.
\end{align}  
\end{subequations}

Moreover, $W$ satisfies the divergence condition 
\begin{equation}
\label{divchange}
\tilde b_i^j W^i,_j=\mathfrak{a}=\tilde b_i^j \zeta^2,_j W^i\ \text{in}\ [0,T]\times{\mathbb R}^3\ .
\end{equation}

Note that we consider the above inner-products over all of $\mathbb R^3$ 
since $W$ and 
its derivatives are compactly supported in $\Omega'$; the contribution outside
of $\Omega'$ is zero regardless of the way in which we extend $\tilde b$ and 
$g$ to $[\Omega']^c$.
This same remark also applies to (\ref{divchange}).

\section{Estimate for (\ref{linear}): the case of the actual coefficients}
\label{9}

\subsection{Energy estimate for $\tilde w_{tt}$ independent of the 
regularization parameter for $\tilde a$} 

We are now going to use the regularity results (\ref{H2_regt}) and 
(\ref{H3_reg}) in the energy inequality (\ref{energywtt}) (which was 
bounded by a constant that does not dependent
on the mollification parameter). 
Our approach now will be to use interpolation inequalities to obtain
an estimate which is independent of the regularization parameter. 

This section will be divided into eight steps, each of which is devoted to 
the estimation of the various integral terms in (\ref{energywtt}). 

In what follows, $\delta>0$ is a given positive number;  the choice of 
$\delta$ will be made precise later, as it will have to be chosen sufficiently
small.

\noindent {\bf Step 1.} 
Let $\displaystyle I_1=\int_0^{t}\int_{\Omega_0^f} {\tilde q}_{t} 
(\tilde a_i^j)_{t} {\tilde w}_{tt}^i,_j$. 
Then, by the Cauchy-Schwarz inequality and by interpolation,
\begin{align*}
I_1&\le \delta \int_0^t \|\nabla \tilde w_{tt}\|^2_{L^2 (\Omega_0^f;{\mathbb R}^9)} +C_{\delta}\int_0^t \|\tilde a_{t}\|^2_{L^4(\Omega_0^f;{\mathbb R}^9)}\|\tilde q_t\|^2_{L^4 (\Omega_0^f;{\mathbb R})}\\
&\le \delta \int_0^t \|\nabla \tilde w_{tt}\|^2_{L^2 (\Omega_0^f;{\mathbb R}^9)} +C_{\delta} C(M)\int_0^t \|\tilde q_t\|^{0.5}_{L^2 (\Omega_0^f;{\mathbb R})} \|\tilde q_t\|^{1.5}_{H^1 (\Omega_0^f;{\mathbb R})}\ ,
\end{align*}
where we have used (\ref{c0}) for the $L^{\infty}$ control of $\tilde a_t$ in $H^1$. 
Thus, 
\begin{align*}
I_1\le \delta \int_0^t \|\nabla \tilde w_{tt}\|^2_{L^2 (\Omega_0^f;{\mathbb R}^9)} +C_{\delta} C(M) \sup_{(0,t)}\|\tilde q_t\|^{0.5}_{L^2 (\Omega_0^f;{\mathbb R})}\ [\int_0^t \|\tilde q_t\|^{2}_{H^1 (\Omega_0^f;{\mathbb R})}]^{\frac{3}{4}}\ T^{\frac{1}{4}}\ . 
\end{align*}
By Lemma \ref{Lagrange} applied to (\ref{weakWt}), and (\ref{c0}),
\begin{align}
\|\tilde q_t\|^{2}_{L^2 (\Omega_0^f;{\mathbb R})}\le C&\ [\ \|\tilde w_{tt}\|^{2}_{L^2 (\Omega;{\mathbb R}^3)}+ \|\tilde q\ \tilde a_t\|^{2}_{L^2 (\Omega_0^f;{\mathbb R}^9)} +  \|F_t\|^{2}_{L^2 (\Omega_0^f;{\mathbb R}^3)}\nonumber\\
& +\|f_t\|^{2}_{L^2 (\Omega_0^s;{\mathbb R}^3)}+ \|\tilde w_t\|^{2}_{H^1 (\Omega_0^f;{\mathbb R}^3)}+  \|\tilde w\|^{2}_{H^2 (\Omega_0^f;{\mathbb R}^3)}+\|\tilde w\|^{2}_{H^1 (\Omega_0^s;{\mathbb R}^3)}\ ]\ .
\label{eliminateqt}
\end{align}
Thus, with $\displaystyle\tilde w_t(t)=w_1+\int_0^t \tilde w_{tt}$, 
$\displaystyle\tilde w(t)=u_0 +\int_0^t \tilde w_{t}$ and $\displaystyle\tilde q=q_0+\int_0^t \tilde q_t$ respectively in
$H^1 (\Omega_0^f;{\mathbb R}^3)$, $H^2 (\Omega_0^f;{\mathbb R}^3)$,
 and $H^1(\Omega_0^f;\mathbb R)$,
\begin{align}
I_1\le &\ \delta \int_0^t \|\nabla \tilde w_{tt}\|^2_{L^2 (\Omega_0^f;{\mathbb R}^9)}\nonumber\\
& +C_{\delta} C(M)\ T^{\frac{1}{4}}\ [\ N(u_0,f)^2+ T \int_0^t \|\tilde q_{t}\|^{2}_{H^1 (\Omega_0^f;{\mathbb R})} + \int_0^t \|\tilde q_t\|^{2}_{H^1 (\Omega_0^f;{\mathbb R})}\nonumber\\
&\qquad\qquad\qquad\qquad + T\ [\int_0^t \|\tilde w_{tt}\|^{2}_{H^1 (\Omega_0^f;{\mathbb R}^3)}+ \int_0^t \|\tilde w_{t}\|^{2}_{H^2 (\Omega_0^f;{\mathbb R}^3)}]\nonumber\\ 
&\qquad\qquad\qquad\qquad + \sup_{[0,T]}\|\tilde w\|^{2}_{H^1 (\Omega_0^s;{\mathbb R}^3)}+ \sup_{[0,T]}\|\tilde w_{tt}\|^{2}_{L^2 (\Omega;{\mathbb R}^3)}\ ] .
\label{I1} 
\end{align}

\noindent {\bf Step 2.} Let $\displaystyle I_2=\int_0^{t}\int_{\Omega_0^f} {\tilde q}_{t} (\tilde a_i^j)_{tt} {\tilde w}_{t}^i,_j$. Then,
\begin{align*}
I_2&\le \int_0^t \|(\tilde a_i^j)_{tt} {\tilde w}_{t}^i,_j\|_{L^{\frac{6}{5}}(\Omega_0^f;{\mathbb R})} \|\tilde q_t\|_{L^6 (\Omega_0^f;{\mathbb R})}\\
&\le \delta \int_0^t \|\tilde q_{t}\|^2_{H^1 (\Omega_0^f;{\mathbb R})} +C_{\delta} \int_0^t \|\nabla \tilde w_t\|^2_{L^3 (\Omega_0^f;{\mathbb R}^9)} \|\tilde a_{tt}\|^{2}_{L^2 (\Omega_0^f;{\mathbb R}^9)}\ ,
\end{align*}
Thus, using (\ref{c0}) for the $L^{\infty}$ control of $\tilde a_{tt}$ in $L^2$,
\begin{align*}
I_2\le \delta \int_0^t \|\tilde q_{t}\|^2_{H^1 (\Omega_0^f;{\mathbb R})} +C_{\delta} \  C(M)\ \int_0^t \|\nabla \tilde w_t\|^{0.5}_{L^2 (\Omega_0^f;{\mathbb R}^9)} \|\nabla \tilde w_t\|^{1.5}_{H^1 (\Omega_0^f;{\mathbb R}^9)}\ ,
\end{align*}
which with $\displaystyle\tilde w_t(t)=\tilde w_1+\int_0^t \tilde w_{tt}$ gives
\begin{align}
I_2\le  &\ \delta \int_0^t \|\tilde q_{t}\|^2_{H^1 (\Omega_0^f;{\mathbb R})}\nonumber\\
& +C_{\delta} \ T^{\frac{1}{4}}\ [\ N(u_0,f)^2+ \int_0^t \|\nabla \tilde w_{t}\|^{2}_{H^1 (\Omega_0^f;{\mathbb R}^9)}+ {T}\ \int_0^t \|\nabla \tilde w_{tt}\|^2_{L^2 (\Omega_0^f;{\mathbb R}^9)}\ ]\ .
\label{I2}
\end{align}

\noindent {\bf Step 3.} Let $\displaystyle I_3=\int_0^{t}\int_{\Omega_0^f} {\tilde q}_{t} (\tilde a_i^j)_{ttt} {\tilde w}^i,_j$. Then,
\begin{align*}
I_3&\le \delta \int_0^t \| \tilde a_{ttt}\|^2_{L^2 (\Omega_0^f;{\mathbb R}^9)}\|\nabla \tilde w\|^2_{L^4 (\Omega_0^f;{\mathbb R}^9)} +C_{\delta}\int_0^t \| \tilde q_{t}\|^2_{L^4(\Omega_0^f;{\mathbb R})}\\
&\le \delta \int_0^t \|\tilde a_{ttt}\|^2_{L^2 (\Omega_0^f;{\mathbb R}^9)} [N(u_0,f)^2+ T\ \int_0^t \|\tilde w_t\|^2_{H^2 (\Omega_0^f;{\mathbb R}^3)}]\\
&\qquad +C_{\delta} \int_0^t \|\tilde q_t\|^{0.5}_{L^2 (\Omega_0^f;{\mathbb R})} \|\tilde q_t\|^{1.5}_{H^1 (\Omega_0^f;{\mathbb R})}\ ,
\end{align*}
where we have used $\displaystyle \tilde w=u_0+\int_0^t \tilde w_t$. Thus, by (\ref{c0}), and (\ref{eliminateqt}),
\begin{align}
I_3\le &\ \delta C(M)\ [\ N(u_0,f)^2+ T\int_0^t \|\tilde w_{t}\|^2_{H^2 (\Omega_0^f;{\mathbb R}^3)}\ ] +C_{\delta} \ T^{\frac{1}{4}}\ \int_0^t \|\tilde q_t\|^2_{H^1(\Omega_0^f;{\mathbb R}^3)}\nonumber\\
& +C_{\delta} C(M)\ T^{\frac{1}{4}}\ [\ N(u_0,f)^2+ T \int_0^t \|\tilde q_{t}\|^{2}_{H^1 (\Omega_0^f;{\mathbb R})} + T\int_0^t \|\tilde w_{tt}\|^{2}_{H^1 (\Omega_0^f;{\mathbb R}^3)} \nonumber\\
&\qquad\qquad\qquad  + T\int_0^t\|\tilde w_t\|^{2}_{H^2 (\Omega_0^f;{\mathbb R}^3)}+ \sup_{[0,T]}\|\tilde w\|^{2}_{H^1 (\Omega_0^s;{\mathbb R}^3)}
+ \sup_{[0,T]}\|\tilde w_{tt}\|^{2}_{L^2 (\Omega;{\mathbb R}^3)}\ ] .
\label{I3} 
\end{align}

\noindent {\bf Step 4.} Let $\displaystyle I_4=\int_0^{t}\int_{\Omega_0^f} {\tilde q} (\tilde a_i^j)_{tt} {\tilde w}_{tt}^i,_j$. 
Then,
\begin{align}
I_4&\le \delta \int_0^t \|\nabla \tilde w_{tt}\|^2_{L^2 (\Omega_0^f;{\mathbb R}^9)} +C_{\delta}\int_0^t \|\tilde a_{tt}\|^2_{L^2(\Omega_0^f;{\mathbb R}^9)}\|\tilde q\|^2_{W^{1,4} (\Omega_0^f;{\mathbb R})}\nonumber\\
&\le \delta \int_0^t \|\nabla \tilde w_{tt}\|^2_{L^2 (\Omega_0^f;{\mathbb R}^9)} +C_{\delta} C(M)\int_0^t \|\tilde q\|^{0.5}_{H^1 (\Omega_0^f;{\mathbb R})} \|\tilde q\|^{1.5}_{H^1 (\Omega_0^f;{\mathbb R})}\nonumber\\
&\le \delta \int_0^t \|\nabla \tilde w_{tt}\|^2_{L^2 (\Omega_0^f;{\mathbb R}^9)}\nonumber\\
&\qquad  +C_{\delta} C(M)\ T^{\frac{1}{4}}\ [\ N(u_0,f)^2+T \int_0^t\|\tilde q_t\|^{2}_{H^1 (\Omega_0^f;{\mathbb R})} + \int_0^t \|\tilde q\|^{2}_{H^2 (\Omega_0^f;{\mathbb R})}\ ]\ .
\label{I4} 
\end{align}

The next two steps will require the introduction of $\delta_1>0$ which is
different than $\delta$ and will also be made precise later.

\noindent {\bf Step 5.} Let $\displaystyle I_5=-\int_{\Omega_0^f} {\tilde q}_{t}(t) (\tilde a_i^j)_{tt} {\tilde w}^i,_j(t)$. We first notice that
$$I_5=-\int_{\Omega_0^f} {\tilde q}_{t}(t) (\tilde a_i^j)_{tt} ({\tilde w}^i,_j(t)-u_0^i,_j)-
      \int_{\Omega_0^f} {\tilde q}_{t}(t) (\tilde a_i^j)_{tt} u_0^i,_j\ .$$
For the second term of the right-hand side of this equality, 
\begin{align*}
 -\int_{\Omega_0^f} {\tilde q}_{t}(t) (\tilde a_i^j)_{tt} u_0^i,_j\le 
\delta_1 \|\tilde q_{t}(t)\|^2_{L^2 (\Omega_0^f;{\mathbb R})} +C_{\delta_1} \|\tilde a_{tt}(t)\|^2_{L^2(\Omega_0^f;{\mathbb R}^9)} \|u_0\|^2_{H^3 (\Omega_0^f;{\mathbb R}^3)}\ ,
\end{align*}
and thus by (\ref{etah3}), since $T\le T_M$,
\begin{equation}
\label{I51}
-\int_{\Omega_0^f} {\tilde q}_{t}(t) (\tilde a_i^j)_{tt} u_0^i,_j\le 
\delta_1 \|\tilde q_{t}(t)\|^2_{L^2 (\Omega_0^f;{\mathbb R})} +C_{\delta_1}C(M)
\ N(u_0,f)^2.
\end{equation}
For the other term,  
\begin{align}
 -\int_{\Omega_0^f} {\tilde q}_{t}(t) (\tilde a_i^j)_{tt} (\tilde w^i,_j (t)-u_0^i,_j)\le &\  
\delta\  \|\tilde q_{t}(t)\|^2_{L^2 (\Omega_0^f;{\mathbb R})}\\
& +C_{\delta} \|\tilde a_{tt}(t)\|^2_{L^3(\Omega_0^f;{\mathbb R}^9)}\|\nabla \tilde w-\nabla u_0\|^2_{L^6 (\Omega_0^f;{\mathbb R}^9)}\ ,
\label{exponent}
\end{align}
and thus, by the $L^{\infty}$ control in $L^3$ provided by (\ref{c4}),
\begin{align}
 \int_{\Omega_0^f} {\tilde q}_{t}(t) (\tilde a_i^j)_{tt} (\tilde w^i,_j (t)-u_0^i,_j)\le &\  
\delta\  \|\tilde q_{t}(t)\|^2_{L^2 (\Omega_0^f;{\mathbb R})}\nonumber\\
& +C_{\delta}\ C(M)\ T\ \int_0^t \|\nabla \tilde w_t\|^2_{H^1 (\Omega_0^f;{\mathbb R}^9)}\ .
\label{I52}
\end{align}
By (\ref{eliminateqt}), (\ref{I51}) and (\ref{I52}), we finally have
\begin{align}
I_5\le &\ (\delta+\delta_1)\  [\ N(u_0,f)^2+ C(M) T\ [\int_0^t \|\tilde q_{t}\|^{2}_{H^1 (\Omega_0^f;{\mathbb R})} + \int_0^t \|\tilde w_{tt}\|^{2}_{H^1 (\Omega_0^f;{\mathbb R}^3)}]\nonumber\\
&\qquad\qquad + C(M)T\int_0^T \|\tilde w_t\|^{2}_{H^2 (\Omega_0^f;{\mathbb R}^3)} + \sup_{[0,T]}\|\tilde w\|^{2}_{H^1 (\Omega_0^s;{\mathbb R}^3)}+ \sup_{[0,T]}\|\tilde w_{tt}\|^{2}_{L^2 (\Omega;{\mathbb R}^3)} ]\nonumber\\
&+C_{\delta}\ C(M)\ T\ \int_0^t \|\nabla \tilde w_t\|^2_{H^1 (\Omega_0^f;{\mathbb R}^9)}+ C_{\delta_1} C\ N(u_0,f)^2\ .
\label{I5}
\end{align}

\begin{remark} Note that $L^3$ and $L^6$ in (\ref{exponent}) are limit cases for both (\ref{c4}) and the Sobolev embeddings in dimension three. In dimension $\ge 4$, this
would no longer be possible and we would be required to introduce a smoother
functional framework.
\end{remark}

\noindent {\bf Step 6.} 
Let $\displaystyle I_6=-\int_{\Omega_0^f} {\tilde q}_{t}(t) (\tilde a_i^j)_{t} 
{w}^i_t,_j(t)$. Similarly to our previous step, we first notice that
$$I_6=-\int_{\Omega_0^f} {\tilde q}_{t}(t) ((\tilde a_i^j)_{t}(t)
-(\tilde a_i^j)_{t}(0)) {\tilde w}_t^i,_j(t)-
\int_{\Omega_0^f} {\tilde q}_{t}(t) (\tilde a_i^j)_{t}(0) 
{\tilde w}_t^i,_j(t)\ .$$
For the second term of the right-hand side of this equality, 
\begin{align*}
- \int_{\Omega_0^f} {\tilde q}_{t}(t) (\tilde a_i^j)_{t}(0) \tilde w_t^i,_j (t)\le 
\delta_1\ \|\tilde a_{t}(0)\|^2_{H^2(\Omega_0^f;{\mathbb R}^9)} \|\tilde q_{t}(t)\|^2_{L^2 (\Omega_0^f;{\mathbb R})} +C_{\delta_1} \|\nabla \tilde w_t (t)\|^2_{L^2 (\Omega_0^f;{\mathbb R}^9)}\ ,
\end{align*}
and thus,
\begin{align}
\label{I61}
-\int_{\Omega_0^f} {\tilde q}_{t}(t) (\tilde a_i^j)_{t}(0) \tilde w_t^i,_j (t)\le & \  
C \delta_1\  N(u_0,f)^2\ \|\tilde q_{t}(t)\|^2_{L^2 (\Omega_0^f;{\mathbb R})}\nonumber\\
& +C_{\delta_1} [\ \|\nabla \tilde w_1\|^2_{L^2 (\Omega_0^f;{\mathbb R}^9)} + T\ \int_0^t\|\nabla \tilde w_{tt}\|^2_{L^2 (\Omega_0^f;{\mathbb R}^9)}\ ]\ .
\end{align}
For the other term,  
\begin{align*}
- \int_{\Omega_0^f} {\tilde q}_{t}(t) ((\tilde a_i^j)_{t}(t)-(\tilde a_i^j)_{t}(0)) \tilde w_t^i,_j (t)\le &\  
\delta\  \|\tilde q_{t}(t)\|^2_{L^2 (\Omega_0^f;{\mathbb R})}\\
& +C_{\delta} \|\tilde a_{t}(t)-\tilde a_t (0)\|^2_{L^6(\Omega_0^f;{\mathbb R}^9)}\|\nabla \tilde w_t\|^2_{L^3 (\Omega_0^f;{\mathbb R}^9)} ,
\end{align*}
and by (\ref{c6}),
\begin{align*}
- \int_{\Omega_0^f} {\tilde q}_{t}(t) ((\tilde a_i^j)_{t}(t)-(\tilde a_i^j)_{t}(0)) \tilde w_t^i,_j (t)\le &\  
\delta\  \|\tilde q_{t}(t)\|^2_{L^2 (\Omega_0^f;{\mathbb R})}\nonumber\\
& +C_{\delta}\ C(M)\ T\  \|\nabla \tilde w_t (t)\|^2_{L^3 (\Omega_0^f;{\mathbb R}^9)}\ .
%\label{I62}
\end{align*}
In the same fashion as we proved (\ref{c4}), we use the $L^{\infty}$ control 
in $L^3$:
$$\|\nabla \tilde w_t (t)\|^2_{L^3 (\Omega_0^f;{\mathbb R}^9)}\le \|\nabla w_1\|^2_{L^3 (\Omega_0^f;{\mathbb R}^9)} + C\ [\int_0^t \|\nabla \tilde w_t \|^2_{H^1 (\Omega_0^f;{\mathbb R}^9)}+\int_0^t \|\nabla \tilde w_{tt} \|^2_{L^2 (\Omega_0^f;{\mathbb R}^9)}],$$
which combined with the previous inequality provides us with
\begin{align}
- \int_{\Omega_0^f} {\tilde q}_{t}(t) ((\tilde a_i^j)_{t}(t)&-(\tilde a_i^j)_{t}(0)) \tilde w_t^i,_j (t)\nonumber\\
&\le \  
\delta\  \|\tilde q_{t}(t)\|^2_{L^2 (\Omega_0^f;{\mathbb R})}+C_{\delta}\ C(M)\ T\  \|\nabla \tilde w_1\|^2_{H^1 (\Omega_0^f;{\mathbb R}^9)}\nonumber\\
&\qquad + C_{\delta} C(M) T\ [\int_0^t \|\nabla \tilde w_t \|^2_{H^1 (\Omega_0^f;{\mathbb R}^9)}+\int_0^t \|\nabla \tilde w_{tt} \|^2_{L^2 (\Omega_0^f;{\mathbb R}^9)}]\ .
\label{I62}
\end{align}

By (\ref{I61}) and (\ref{I62}), we finally have that
\begin{align}
I_6\le &\ (\delta+\delta_1)\ [\ C N(u_0,f)^2+  C(M) T\ \int_0^t \|\tilde w_{tt}\|^{2}_{H^1 (\Omega_0^f;{\mathbb R}^3)} \nonumber\\
&\qquad\qquad+ C(M)T [\int_0^t\|\tilde q_t\|^{2}_{H^1 (\Omega_0^f;{\mathbb R})}+
\int_0^t\|\tilde w_t\|^{2}_{H^2 (\Omega_0^f;{\mathbb R}^3)}] 
+ \sup_{[0,t]}\|\tilde w_{tt}\|^{2}_{L^2 (\Omega;{\mathbb R}^3)}\nonumber\\
&\qquad\qquad + \sup_{[0,t]}\|\tilde w\|^{2}_{H^1 (\Omega_0^s;{\mathbb R}^3)}\ ]+C_{\delta_1}\ [\ N(u_0,f)^2 + T\ \int_0^t\|\nabla \tilde w_{tt}\|^2_{L^2 (\Omega_0^f;{\mathbb R}^9)}\ ]\nonumber\\
&  +C_{\delta} C(M)\ T\ [\ N(u_0,f)^2 + \int_0^t\|\nabla \tilde w_{tt}\|^2_{L^2 (\Omega_0^f;{\mathbb R}^9)}+
\int_0^t\|\nabla \tilde w_{t}\|^2_{H^1 (\Omega_0^f;{\mathbb R}^9)} \ ]\ .
\label{I6}
\end{align}

\noindent {\bf Step 7.} 
Let $\displaystyle I_7=-\int_0^{t} ((\tilde a_k^r \tilde a_k^s)_{tt} {\tilde w},_r,\  {\tilde w}_{tt},_s)_{L^2(\Omega_0^f;{\mathbb R}^3)}$. Then,
\begin{align*}
I_7&\le \delta \int_0^t \|\nabla \tilde w_{tt}\|^2_{L^2 (\Omega_0^f;{\mathbb R}^9)} +C_{\delta}\int_0^t \|\tilde (\tilde a_k^r \tilde a_k^s)_{tt}\|^2_{L^2(\Omega_0^f;{\mathbb R})}\|\nabla \tilde w\|^2_{W^{1,4} (\Omega_0^f;{\mathbb R}^9)}\\
&\le \delta \int_0^t \|\nabla \tilde w_{tt}\|^2_{L^2 (\Omega_0^f;{\mathbb R}^9)} +C_{\delta} C(M)\int_0^t \|\nabla \tilde w\|^{0.5}_{H^1 (\Omega_0^f;{\mathbb R}^9)} \|\nabla \tilde w\|^{1.5}_{H^2 (\Omega_0^f;{\mathbb R}^9)}\ ,
\end{align*}
where we have used (\ref{c0}) for the $L^{\infty}$ control of $\tilde a_{tt}$, $\tilde a_t$ and $\tilde a$ respectively in $L^2$, $H^1$ and $H^2$. 
Thus,
\begin{align}
I_7\le &\ \delta \int_0^t \|\nabla \tilde w_{tt}\|^2_{L^2 (\Omega_0^f;{\mathbb R}^9)}\nonumber\\
& +C_{\delta} C(M)\ T^{\frac{1}{4}}\ [\ N(u_0,f)^2+ T\ \int_0^t \|\nabla \tilde w_t\|^{2}_{H^1 (\Omega_0^f;{\mathbb R}^9)} + \int_0^t \|\nabla \tilde w\|^{2}_{H^2 (\Omega_0^f;{\mathbb R}^9)}\ ]\ .
\label{I7} 
\end{align}

\noindent {\bf Step 8.} 
Let $\displaystyle I_8=-\int_0^{t} ((\tilde a_k^r \tilde a_k^s)_{t} {\tilde w_t},_r,\  {\tilde w}_{tt},_s)_{L^2(\Omega_0^f;{\mathbb R}^3)}$. Then,
\begin{align*}
I_8&\le \delta \int_0^t \|\nabla \tilde w_{tt}\|^2_{L^2 (\Omega_0^f;{\mathbb R}^9)} +C_{\delta}\int_0^t \|\tilde (\tilde a_k^r \tilde a_k^s)_{t}\|^2_{L^{4}(\Omega_0^f;{\mathbb R})}\|\nabla \tilde w_t\|^2_{L^4(\Omega_0^f;{\mathbb R}^9)}\\
&\le \delta \int_0^t \|\nabla \tilde w_{tt}\|^2_{L^2 (\Omega_0^f;{\mathbb R}^9)} +C_{\delta} C(M)\int_0^t \|\nabla \tilde w_t\|^{0.5}_{L^2 (\Omega_0^f;{\mathbb R}^9)} \|\nabla \tilde w_t\|^{1.5}_{H^1 (\Omega_0^f;{\mathbb R}^9)}\ .
\end{align*}
 Consequently, 
%\begin{align*}
%I_8\le &\ \delta \int_0^t \|\nabla \tilde w_{tt}\|^2_{L^2 (\Omega_0^f;{\mathbb R}^9)}+C_{\delta} C(M)  [\ \|\nabla \tilde w_1\|^2_{H^1 (\Omega_0^f;{\mathbb R}^9)} + {T} \int_0^t \|\nabla \tilde w_{tt}\|^{2}_{L^2 (\Omega_0^f;{\mathbb R}^9)} ]^{\frac{1}{4}}\times\nonumber\\
%&\qquad\qquad\qquad\qquad\qquad\qquad\qquad\qquad\qquad\ T^{\frac{1}{4}}\ [\ \int_0^t \|\nabla \tilde w_t\|^{2}_{H^1 (\Omega_0^f;{\mathbb R}^9)}\ ]^{\frac{3}{4}}\  ,
%\end{align*}
%and finally,
\begin{align}
I_8\le &\ \delta \int_0^t \|\nabla \tilde w_{tt}\|^2_{L^2 (\Omega_0^f;{\mathbb R}^9)}\nonumber\\
& +C_{\delta} C(M)\ T^{\frac{1}{4}}\  [\ N(u_0,f) + {T} \int_0^t \|\nabla \tilde w_{tt}\|^{2}_{L^2 (\Omega_0^f;{\mathbb R}^9)}+ \int_0^t \|\nabla \tilde w_t\|^{2}_{H^1 (\Omega_0^f;{\mathbb R}^9)}\ ]\ .
\label{I8} 
\end{align}

\noindent{\bf Step 9.} 
Thus, from (\ref{energywtt}), and estimates (\ref{I1})-(\ref{I8}), we then 
obtain the inequality
\begin{align}
 &\ \sup_{[0,T]} \|\tilde w_{tt} (t)\|^2_{L^2(\Omega;{\mathbb R}^3)} 
+\int_0^{T} \|{\tilde w}_{tt}\|^2_{H^1(\Omega_0^f;{\mathbb R}^3)} 
+ \sup_{[0,T]} \|{{\tilde w}_{t}}\|^2_{H^1(\Omega_0^s;{\mathbb R}^3)}\nonumber\\
& \le [\ C \delta (1+C(M))N(u_0,f)^2+
C \delta_1 (1+T C(M)+N(u_0,f)^2)\ ]\ \|(\tilde w,\tilde q)\|^2_{Z_T}\nonumber\\
&\qquad + C_{\delta_1} N(u_0,f)^2 +C_{\delta_1} T^{\frac{1}{4}}\ \|(\tilde w,\tilde q)\|^2_{Z_T} + C_{\delta}C(M) T^{\frac{1}{4}}\ (\|(\tilde w,\tilde q)\|^2_{Z_T}+N(u_0,f)^2) .
\label{energywttlimit}
\end{align}

\subsection{Estimate of $\tilde w_t$ independent of the regularization of $\tilde a$}\hfill\break

In this subsection, let $\Psi_i$ be one of the $H^3$ charts defining a neighborhood of $\Omega_0^s$ and let $W_i=\zeta_i^2\ \tilde w\circ \Psi_i$. Since the estimates that follow do not depend on the choice of $\Psi_i$, we will 
denote $W_i$ simply by $W$.

Recall that for all $\phi \in L^2(0,T; H^1({\mathbb R}^3;{\mathbb R}^3))$,	
 \begin{align*}
& \int_0^T ( W_{tt}, \phi)_{L^2({\mathbb R}^3;{\mathbb R}^3)}\ dt  
+\nu \int_0^T ((\tilde b_k^r \tilde b_k^s  W,_r)_t,\  \phi,_s)_{L^2({\mathbb R}^3_+;{\mathbb R}^3)}\ dt \nonumber\\
&+ \int_0^T (C^{irks} W^k,_r, {\phi}^i,_s)_{L^2({\mathbb R}^3_-;{\mathbb R})}\ dt  - \int_0^T ((\tilde b_i^j Q)_t,\  \phi^i,_j)_{L^2({\mathbb R}^3_+;{\mathbb R})} dt\nonumber\\
&\qquad\qquad = \int_0^T ({F_1}_t, \phi)_{L^2({\mathbb R}^3_+;{\mathbb R}^3)}+ 
({H_i}_t, \phi,_i)_{L^2({\mathbb R}^3_+;{\mathbb R}^3)}\ dt \nonumber\\
&\qquad\qquad\qquad +\int_0^T ({F_2}_t, \phi)_{L^2({\mathbb R}^3_-;{\mathbb R}^3)}+ 
({K_i}_t, \phi,_i)_{L^2({\mathbb R}^3_-;{\mathbb R}^3)}\ dt\ .
\end{align*}
With the choice of $\phi=D_{-h}D_h W_t$ in this variational formulation, which is possible since $\tilde w_t\in L^2(0,T;H^1_0(\Omega;{\mathbb R}^3))$, we then get
 \begin{align}
&  \frac{1}{2} \|D_h W_{t}(T)\|^2_{L^2({\mathbb R}^3;{\mathbb R}^3)}  
+\nu \int_0^T (\tilde b_k^r \tilde b_k^s D_h {W_t},_r,\ D_h {W_t},_s)_{L^2({\mathbb R}^3_+;{\mathbb R}^3)} \nonumber\\
&+ \frac{1}{2} (C^{irks}D_h W^k,_r (T),\ D_h W^i,_s (T))_{L^2({\mathbb R}^3_-;{\mathbb R})}  - \int_0^T (D_h(\tilde b_i^j Q)_t,\  D_h W_t^i,_j)_{L^2({\mathbb R}^3_+;{\mathbb R})}\nonumber\\
&+\nu \int_0^T (D_h(\tilde b_k^r \tilde b_k^s) {W_t},_r ^h,\ D_h {W_t},_s)_{L^2({\mathbb R}^3_+;{\mathbb R}^3)}+\nu \int_0^T ((\tilde b_k^r \tilde b_k^s)_t D_h {W},_r,\ D_h {W_t},_s)_{L^2({\mathbb R}^3_+;{\mathbb R}^3)}\nonumber\\
&+\nu \int_0^T (D_h (\tilde b_k^r \tilde b_k^s)_t {W},_r^h,\ D_h {W_t},_s)_{L^2({\mathbb R}^3_+;{\mathbb R}^3)}+ \int_0^T (D_h C^{irks} W^k,_r,  D_h W_t^i,_s)_{L^2({\mathbb R}^3_-;{\mathbb R})}\nonumber\\ 
&\qquad \le C\ N(u_0,f)^2+ \int_0^T ({F_1}_t, D_{-h}D_h W_t)_{L^2({\mathbb R}^3_+;{\mathbb R}^3)}+ 
(D_h {H_i}_t,\ D_h W_t,_i)_{L^2({\mathbb R}^3_+;{\mathbb R}^3)} \nonumber\\
&\qquad\qquad\qquad +\int_0^T ({F_2}_t,D_{-h}D_h W_t)_{L^2({\mathbb R}^3_-;{\mathbb R}^3)}+ 
(D_h {H_i}_t, D_h W_t,_i)_{L^2({\mathbb R}^3_-;{\mathbb R}^3)}\ .
\label{energywtH2}
\end{align}
Since the estimation of the integrals with the indefinite sign in this 
inequality does not create any new difficulty with respect to the estimates 
that we have obtained  in the previous
subsection (they are even easier since the more difficult integrals $I_5$ and 
$I_6$ do not have an analogue here), we provide the details in the appendix.
With $\delta>0$ to be fixed later, this leads us to
 \begin{align*}
& \int_0^T \|D_h \nabla {W_t}\|^2_{L^2({\mathbb R}^3_+;{\mathbb R}^9)}+ \sup_{[0,T]}\|D_h \nabla W \|^2_{L^2({\mathbb R}^3_+;{\mathbb R}^9)}\nonumber\\ 
& \le [\ C \delta (1+C(M))N(u_0,f)^2+C \delta_1 (1+T C(M)+N(u_0,f)^2)\ ]\ \|(\tilde w,\tilde q)\|^2_{Z_T}\nonumber\\
&\qquad + C_{\delta_1} N(u_0,f)^2 +C_{\delta_1} T^{\frac{1}{4}}\ \|(\tilde w,\tilde q)\|^2_{Z_T} + C_{\delta}C(M) T^{\frac{1}{4}}\ (\|(\tilde w,\tilde q)\|^2_{Z_T}+N(u_0,f)^2) .
\end{align*}
We remark here that the estimates obtained in this section could have been 
performed with $t\in (0,T)$ generically replacing $T$; this explains the 
presence of 
$$\sup_{[0,T]}\|D_h \nabla W \|_{L^2({\mathbb R}^3_+;{\mathbb R}^9)}$$ 
on the left-hand side of this inequality. As this inequality is independent of
$h$, we then deduce that
 \begin{align*}
& \int_0^T \|\nabla_0\nabla {W_t}\|^2_{L^2({\mathbb R}^3_+;{\mathbb R}^{18})}+ \sup_{[0,T]}\|\nabla_0 \nabla W \|^2_{L^2({\mathbb R}^3_+;{\mathbb R}^{18})}\nonumber\\
& \le [\ C \delta (1+C(M))N(u_0,f)^2+C \delta_1 (1+T C(M)+N(u_0,f)^2)\ ]\ \|(\tilde w,\tilde q)\|^2_{Z_T}\nonumber\\
&\qquad + C_{\delta_1} N(u_0,f)^2 +C_{\delta_1} T^{\frac{1}{4}}\ \|(\tilde w,\tilde q)\|^2_{Z_T} + C_{\delta}C(M) T^{\frac{1}{4}}\ (\|(\tilde w,\tilde q)\|^2_{Z_T}+N(u_0,f)^2) ,
\end{align*}
and thus for the trace, where we will denote for notational convenience ${\mathbb R}^2=\{x_3=0\}$, 
 \begin{align*}
& \int_0^T \|\nabla_0 {W_t}\|^2_{H^{0.5}({\mathbb R}^2;{\mathbb R}^{6})}+ \sup_{[0,T]}\|\nabla_0  W \|^2_{H^{0.5}({\mathbb R}^2;{\mathbb R}^{6})}\nonumber\\ 
& \le [\ C \delta (1+C(M))N(u_0,f)^2+C \delta_1 (1+T C(M)+N(u_0,f)^2)\ ]\ \|(\tilde w,\tilde q)\|^2_{Z_T}\nonumber\\
&\qquad + C_{\delta_1} N(u_0,f)^2 +C_{\delta_1} T^{\frac{1}{4}}\ \|(\tilde w,\tilde q)\|^2_{Z_T} + C_{\delta}C(M) T^{\frac{1}{4}}\ (\|(\tilde w,\tilde q)\|^2_{Z_T}+N(u_0,f)^2) ,
\end{align*}
which implies that
 \begin{align*}
 &\int_0^T \|{W_t}\|^2_{H^{1.5}({\mathbb R}^2;{\mathbb R}^{3})}+ \sup_{[0,T]}\|  W \|^2_{H^{1.5}({\mathbb R}^2;{\mathbb R}^{3})}\nonumber\\ 
& \le [\ C \delta (1+C(M))N(u_0,f)^2+C \delta_1 (1+T C(M)+N(u_0,f)^2)\ ]\ \|(\tilde w,\tilde q)\|^2_{Z_T}\nonumber\\
&\qquad + C_{\delta_1} N(u_0,f)^2 +C_{\delta_1} T^{\frac{1}{4}}\ \|(\tilde w,\tilde q)\|^2_{Z_T} + C_{\delta}C(M) T^{\frac{1}{4}}\ (\|(\tilde w,\tilde q)\|^2_{Z_T}+N(u_0,f)^2) .
\end{align*}
Since this has been done for any $W=\zeta_i^2\ \tilde w\circ \Psi_i$, we then deduce by the finite covering argument and the fact that each $\Psi_i$ is of class $H^3$ that
 \begin{align}
& \int_0^T \|{\tilde w_t}\|^2_{H^{1.5}(\Gamma_0;{\mathbb R}^{3})}+ \sup_{[0,T]}\|  \tilde w\|^2_{H^{1.5}(\Gamma_0;{\mathbb R}^{3})}\nonumber\\ 
& \le [\ C \delta (1+C(M))N(u_0,f)^2+C \delta_1 (1+T C(M)+N(u_0,f)^2)\ ]\ \|(\tilde w,\tilde q)\|^2_{Z_T}\nonumber\\
&\qquad + C_{\delta_1} N(u_0,f)^2 +C_{\delta_1} T^{\frac{1}{4}}\ \|(\tilde w,\tilde q)\|^2_{Z_T} + C_{\delta}C(M) T^{\frac{1}{4}}\ (\|(\tilde w,\tilde q)\|^2_{Z_T}+N(u_0,f)^2) .
\label{tracewt}
\end{align}
Elliptic regularity for the Stokes problem (see \cite{Eben}) (for $t\in [0,T]$ considered as fixed)
\begin{align*}
 -\nu\triangle [\tilde w^i_t\circ\tilde\eta^{-1}] + (\tilde q_t\circ\tilde\eta^{-1}),_i&= -\tilde w^i_{tt}\circ\tilde\eta^{-1} +F^i_t\circ{\tilde\eta}^{-1}+\nu (\tilde a_l^j,_j \tilde a_l^k \tilde w^i,_k)_t\circ\tilde\eta^{-1}\\
&\qquad-[\tilde a_i^j,_j \tilde q]_t\circ\tilde\eta^{-1}
 -[(\tilde a_i^k)_t \tilde q,_k]\circ{\tilde\eta}^{-1}\\
&\qquad+
\nu [(\tilde a_l^j \tilde a_l^k)_t \tilde w^i,_k],_j\circ{\tilde\eta}^{-1}
\ \text{in}\ \tilde\eta(t,\Omega_0^f)\\
\operatorname{div}(\tilde w_t\circ\tilde\eta^{-1}) (t,\cdot)&=-[(\tilde a_i^j)_t w^i,_j]\circ {\tilde\eta}^{-1} \ \text{in}\ \tilde\eta(t,\Omega_0^f)\\
\tilde w_t\circ\tilde\eta^{-1} (t,\cdot)&=0\ \text{on}\ \tilde\eta(t,\partial\Omega)\\
\tilde w_t\circ\tilde\eta^{-1} (t,\cdot)&=\tilde w_t\circ\tilde\eta^{-1} (t,\cdot)\ \text{on}\ \tilde\eta(t,\Gamma_0)\ ,
\end{align*}
then implies with the $L^{\infty}$ in time estimate (\ref{etah3}) of $\eta$ 
(and thus of $\tilde\eta$) into $H^3$
\begin{align*}
&\|\tilde w_t\circ\tilde\eta^{-1} (t)\|_{H^2(\tilde\eta(t,\Omega_0^f);{\mathbb R}^3)}+ \|\tilde q_t\circ\tilde\eta^{-1} (t)\|_{H^1(\tilde\eta(t,\Omega_0^f);{\mathbb R})}\\
&\le C  
[ \|-\tilde w^i_{tt}\circ\tilde\eta^{-1} +F^i_t\circ\tilde\eta^{-1}-[(\tilde a_i^k)_t \tilde q,_k]\circ{\tilde\eta}^{-1}+
\nu [(\tilde a_l^j \tilde a_l^k)_t \tilde w^i,_k],_j\circ{\tilde\eta}^{-1}\|_{L^2(\tilde\eta(t,\Omega_0^f);{\mathbb R}^3)}\nonumber\\
& \qquad\qquad +\|(\tilde a_l^j,_j \tilde a_l^k \tilde w,_k)_t\|^2_{L^2(\tilde\eta(t,\Omega_0^f);{\mathbb R}^3)} +
\|\tilde w_t\circ\tilde\eta^{-1}(t) \|_{H^{1.5}(\tilde\eta(t,\Gamma_0);{\mathbb R}^3)}\nonumber\\
&\qquad\qquad\qquad+\|[\tilde a_i^j,_j \tilde q]_t\circ\tilde\eta^{-1}\|^2_{L^2(\tilde\eta(t,\Omega_0^f);{\mathbb R})} ]\ .
\end{align*}
Thus, still with (\ref{etah3}) and (\ref{c6}),
\begin{align*}
&\|\tilde w_t(t)\|_{H^2(\Omega_0^f;{\mathbb R}^3)}- C\ \| \nabla\tilde w_t(t)\|_{L^4(\Omega_0^f;{\mathbb R}^9)} + \|\tilde q_t(t)\|_{H^1(\Omega_0^f;{\mathbb R})}\\
&\le\ C\   
[\ \|\tilde w_{tt}\|_{L^2(\Omega_0^f;{\mathbb R}^3)}+\|\nabla \tilde w\|_{W^{1,4}(\Omega_0^f;{\mathbb R}^9)}+\|\nabla \tilde w\|_{L^{4}(\Omega_0^f;{\mathbb R}^9)}+\|\nabla \tilde w_t\|_{L^{4}(\Omega_0^f;{\mathbb R}^9)}\\
&\qquad\qquad+
\|\nabla \tilde q\|_{L^4(\Omega_0^f;{\mathbb R}^3)} +\|f_t\|_{L^2(\Omega;{\mathbb R}^3)}+
\|\tilde w_t\|_{H^{1.5}(\Gamma_0;{\mathbb R}^3)}+\sqrt{T} \|\tilde q_t\|_{L^2(\Omega_0^f;{\mathbb R})} 
 ]\ ,
\end{align*}
from which we immediately infer that
\begin{align*}
&\int_0^T \|\tilde w_t(t)\|^2_{H^2(\Omega_0^f;{\mathbb R}^3)} + \int_0^T \|\tilde q_t(t)\|^2_{H^1(\Omega_0^f;{\mathbb R})}\\
& \le\ C\   
[\ \int_0^T\|\tilde w_{tt}\|^2_{L^2(\Omega_0^f;{\mathbb R}^3)}
+ \int_0^T\| f_{t}\|^2_{L^2(\Omega;{\mathbb R}^3)}
+ \int_0^T\|\tilde w_{t}\|^2_{H^{1.5}(\Gamma_0;{\mathbb R}^3)}
\\
&\qquad\qquad+ T^{\frac{1}{4}} [\ \|\nabla \tilde w_1\|^2_{H^1(\Omega_0^f;{\mathbb R}^9)}+ T\int_0^T \|\nabla \tilde w_{tt}\|^2_{H^1(\Omega_0^f;{\mathbb R}^9)}+\int_0^T \|\nabla \tilde w_{t}\|^2_{H^1(\Omega_0^f;{\mathbb R}^9)}\ ]\\ 
&\qquad\qquad+ T^{\frac{1}{4}} [\ \|\nabla u_0\|^2_{H^1(\Omega_0^f;{\mathbb R}^9)}+ T \int_0^T \|\nabla \tilde w_{t}\|^2_{H^1(\Omega_0^f;{\mathbb R}^9)}+\int_0^T \|\nabla \tilde w\|^2_{H^2(\Omega_0^f;{\mathbb R}^9)}\ ]\\
&\qquad\qquad+ T^{\frac{1}{4}} [\ N(u_0,f)^2+ T \int_0^T \|\tilde q_{t}\|^2_{L^2(\Omega_0^f;{\mathbb R})}+\int_0^T \|\tilde q\|^2_{H^1(\Omega_0^f;{\mathbb R})}\ ]\ ]
\ .
\end{align*}
Thus, with the trace estimate (\ref{tracewt}) and (\ref{energywttlimit}),
\begin{align}
& \int_0^T \|{\tilde w_t}\|^2_{H^{2}(\Omega_0^f;{\mathbb R}^{3})}+ \int_0^T \|{\tilde q_t}\|^2_{H^{1}(\Omega_0^f;{\mathbb R})}\nonumber\\ 
&\le [\ C \delta (1+C(M))N(u_0,f)^2+C {\delta_1} (1+T C(M)+N(u_0,f)^2)\ ]\ \|(\tilde w,\tilde q)\|^2_{Z_T} \nonumber\\
&\qquad+ C_{\delta_1} N(u_0,f)^2 +C_{\delta}C(M)T^{\frac{1}{4}}\ (\|(\tilde w,\tilde q)\|^2_{Z_T}+N(u_0,f)^2) +C_{\delta_1} T^{\frac{1}{4}} \|(\tilde w,\tilde q)\|^2_{Z_T}\ .
\label{tracewt1}
\end{align}

Similarly, the classical elliptic regularity theory for the elasticity problem 
(for $t\in [0,T]$)
\begin{align*}
 -[c^{ijkl} \tilde w^k,_l],_j &= -{\tilde w}_{tt} +f_t\ \text{in}\ \Omega_0^s\\
\tilde w (t,\cdot)&=\tilde w (t,\cdot)\ \text{on}\ \Gamma_0=\partial\Omega_0^s\ ,
\end{align*}
implies that
$\|\tilde w\|_{H^{2}(\Omega_0^s;{\mathbb R}^3)}\le C\ [\ \|-\tilde w_{tt}+f_t\|_{L^{2}(\Omega_0^s;{\mathbb R}^3)}+\|\tilde w\|_{H^{1.5}(\Gamma_0;{\mathbb R}^3)}\ ]\ ,$
which with (\ref{tracewt}) and (\ref{energywttlimit}) provides us with the 
estimate
\begin{align}
&\sup_{[0,T]}\|{\tilde w}\|^2_{H^{2}(\Omega_0^s;{\mathbb R}^{3})}\nonumber\\
&\le [\ C \delta (1+C(M))N(u_0,f)^2+C {\delta_1} (1+T C(M)+N(u_0,f)^2)\ ]\ \|(\tilde w,\tilde q)\|^2_{Z_T} \nonumber\\
&\qquad+ C_{\delta_1} N(u_0,f)^2 +C_{\delta}C(M)T^{\frac{1}{4}}\ (\|(\tilde w,\tilde q)\|^2_{Z_T}+N(u_0,f)^2) +C_{\delta_1} T^{\frac{1}{4}} \|(\tilde w,\tilde q)\|^2_{Z_T}\ .
\label{tracewt2}
\end{align} 

\subsection{Estimate of $\tilde w$ independent of the regularization 
of $\tilde a$}\hfill\break

Just as in the previous subsection, $W$ again denotes $W_i=\zeta^2\ 
\tilde w\circ \Psi_i$, where recall that $\Psi_i$ denotes the $i$th chart.

Choosing $\phi=D_{-h}D_h D_{-h}D_h W$ in the variational formulation 
(\ref{variationalr3}), we then find that
 \begin{align}
&  \frac{1}{2} \|D_{-h}D_h W(T)\|^2_{L^2({\mathbb R}^3;{\mathbb R}^3)}  
+\nu \int_0^T (\tilde b_k^r \tilde b_k^s D_{-h} D_h {W},_r,\ D_{-h} D_h {W},_s)_{L^2({\mathbb R}^3_+;{\mathbb R}^3)} \nonumber\\
&+ \frac{1}{2} (C^{irks}D_{-h}D_h \int_0^{T} W^k,_r ,\ D_{-h}D_h \int_0^{T} W^i,_s )_{L^2({\mathbb R}^3_-;{\mathbb R})}\nonumber\\
&  - \int_0^T (D_{-h}D_h(\tilde b_i^j Q),\ D_{-h} D_h W^i,_j)_{L^2({\mathbb R}^3_+;{\mathbb R})}\nonumber\\
&+\nu \int_0^T (D_{-h}D_h(\tilde b_k^r \tilde b_k^s) {W},_r,\ D_{-h}D_h {W},_s)_{L^2({\mathbb R}^3_+;{\mathbb R}^3)}\nonumber\\
&-\sum_{p=0}^1 \nu \int_0^T (D_{(-1)^p h} (\tilde b_k^r \tilde b_k^s) D_{(-1)^{p} h}{W},_r,\ D_{-h} D_h {W},_s)_{L^2({\mathbb R}^3_+;{\mathbb R}^3)}\nonumber\\
&+ \int_0^T (D_{-h} D_h[C^{irks}]  \int_0^{T}W^k,_r, D_{-h} D_h W^i,_s)_{L^2({\mathbb R}^3_-;{\mathbb R})}\nonumber\\ 
&- \sum_{p=0}^1 \int_0^T ( D_{(-1)^p h}[C^{irks}]  \int_0^{T} D_{(-1)^{p}h} W^k,_r, D_{-h} D_h W^i,_s)_{L^2({\mathbb R}^3_-;{\mathbb R})}\nonumber\\ 
\le &\  C\ N(u_0,f)^2+ \int_0^T (D_{-h}F_1, D_{h}D_{-h}D_h W)_{L^2({\mathbb R}^3_+;{\mathbb R}^3)}\nonumber\\
& \qquad\qquad\qquad + 
\int_0^T (D_{-h}D_h H_i,\ D_{-h}D_h W,_i)_{L^2({\mathbb R}^3_+;{\mathbb R}^3)} \nonumber\\
&+\int_0^T (D_{-h}D_{h}F_2,\ D_{-h}D_h W)_{L^2({\mathbb R}^3_-;{\mathbb R}^3)}+ 
(D_{-h}D_h K_i,\  D_h W,_i)_{L^2({\mathbb R}^3_-;{\mathbb R}^3)}\ .
\label{energywlimit}
\end{align}
Similarly as in the previous
subsection, the estimates provided in the appendix yield (with $\delta>0$ to 
be fixed later)
 \begin{align*}
 \int_0^T \|D_{-h}D_h \nabla {W}\|^2&_{L^2({\mathbb R}^3_+;{\mathbb R}^9)}+ \sup_{[0,T]}\|D_{-h}D_h \nabla \int_0^{\cdot} W \|^2_{L^2({\mathbb R}^3_+;{\mathbb R}^9)}\nonumber\\ &\le \ C\delta\ (1+C(M))N(u_0,f)^2\ \|(\tilde w,\tilde q)\|^2_{Z_T} + C N(u_0,f)^2\\
& \qquad +C_{\delta}C(M)T^{\frac{1}{4}}\ (\|(\tilde w,\tilde q)\|^2_{Z_T}+N(u_0,f)^2)\  .
\end{align*}
As this inequality is independent of $h$, we deduce just as in the previous 
section that
\begin{align}
 \int_0^T \|{\tilde w}\|^2_{H^{2.5}(\Gamma_0;{\mathbb R}^{3})}&+ \sup_{[0,T]}\|  \int_0^{\cdot} \tilde w\|^2_{H^{2.5}(\Gamma_0;{\mathbb R}^{3})}\nonumber\\ &\le \ C\delta\ (1+C(M))N(u_0,f)^2\ \|(\tilde w,\tilde q)\|^2_{Z_T} + C N(u_0,f)^2\nonumber\\
&\qquad+C_{\delta}C(M)T^{\frac{1}{4}}\ (\|(\tilde w,\tilde q)\|^2_{Z_T}+N(u_0,f)^2)\  .
\label{tracew}
\end{align}
Elliptic regularity for the Stokes problem (see \cite{Eben}) for $t\in [0,T]$
\begin{align*}
 -\nu \triangle [\tilde w\circ\tilde\eta^{-1}](t,\cdot) +\nabla (\tilde q\circ\tilde\eta^{-1})(t,\cdot)&= -\tilde w_t\circ\tilde\eta^{-1}
 +f +\nu \tilde a_l^j,_j\circ\eta^{-1} (\tilde w\circ\eta^{-1}),_l\\
&\qquad -(\tilde a_i^j,_j \tilde q)\circ\tilde\eta^{-1}\text{in}\ \tilde\eta(t,\Omega_0^f)\\
\operatorname{div}(\tilde w\circ\tilde\eta^{-1}) (t,\cdot)&=0 \ \text{in}\ \tilde\eta(t,\Omega_0^f)\\
(\tilde w\circ\tilde\eta^{-1}) (t,\cdot)&=0\ \text{on}\ \tilde\eta(t,\partial\Omega)\\
(\tilde w\circ\tilde\eta^{-1}) (t,\cdot)&=(\tilde w\circ\tilde\eta^{-1}) (t,\cdot)\ \text{on}\ \tilde\eta(t,\Gamma_0)\ ,
\end{align*}
 then implies with (\ref{etah3}) 
\begin{align*}
&\|\tilde w\circ\tilde\eta^{-1} (t)\|_{H^3(\tilde\eta(t,\Omega_0^f);{\mathbb R}^3)}+ \|\tilde q\circ\tilde\eta^{-1}(t)\|_{H^2(\tilde\eta(t,\Omega_0^f);{\mathbb R})}\\
&\qquad\qquad \le\ C\ 
[\ \|-\tilde w_t\circ\tilde\eta^{-1} +f+\nu \tilde a_l^j,_j\circ\eta^{-1} (\tilde w\circ\eta^{-1}),_l\|_{H^{1}(\tilde\eta(t,\Omega_0^f);{\mathbb R}^3)}\\
&\qquad\qquad\qquad\qquad+\|(\tilde a_i^j,_j\tilde q)\circ\tilde\eta^{-1}\|_{H^{1}(\tilde\eta(t,\Omega_0^f);{\mathbb R})}+
\|\tilde w\circ\tilde\eta^{-1}(t) \|_{H^{2.5}(\tilde\eta(t,\Gamma_0);{\mathbb R}^3)}]\ .
\end{align*}
Thus, with (\ref{etah3})
\begin{align*}
&\|\tilde w(t)\|_{H^3(\Omega_0^f;{\mathbb R}^3)}- C\ \|  \tilde w(t)\|_{W^{2,4}(\Omega_0^f;{\mathbb R}^3)} + \|\tilde q(t)\|_{H^2(\Omega_0^f;{\mathbb R})}- C\ \| \tilde q(t)\|_{W^{1,4}(\Omega_0^f;{\mathbb R})}\\
&\qquad \le\ C\   
[\  \|\tilde w_t\|_{H^1(\Omega_0^f;{\mathbb R}^3)}
+\|\tilde w\|_{W^{2,4}(\Omega_0^f;{\mathbb R}^3)}+\|f\|_{H^1(\Omega;{\mathbb R}^3)}+\|\tilde w\|_{H^{2.5}(\Gamma_0;{\mathbb R}^3)}]\ ,
\end{align*}
from which we immediately infer,
\begin{align*}
&\int_0^T\|\tilde w(t)\|^2_{H^3(\Omega_0^f;{\mathbb R}^3)} + \int_0^T\|\tilde q(t)\|^2_{H^2(\Omega_0^f;{\mathbb R})}\\
&\le\ C\   
[\ T N(u_0,f)^2+T^2 \int_0^T\|\tilde w_{t}\|^2_{H^1(\Omega_0^f;{\mathbb R}^3)}+\int_0^T\|f\|^2_{H^1(\Omega;{\mathbb R}^3)}+ \int_0^T\|\tilde w\|^2_{H^{2.5}(\Gamma_0;{\mathbb R}^3)}\\
&\qquad\qquad\qquad+C\ T^{\frac{1}{4}}\ [\ N(u_0,f)^2+ T \int_0^T\|\tilde q_t\|^2_{H^1(\Omega_0^f;{\mathbb R}^3)}+ \int_0^T \|\tilde q\|^2_{H^2(\Omega_0^f;{\mathbb R}^3)}
]\\
&\qquad\qquad\qquad+C\ T^{\frac{1}{4}}\ [\|u_0\|^2_{H^2(\Omega_0^f;{\mathbb R}^3)}+ T \int_0^T\|\tilde w_t\|^2_{H^2(\Omega_0^f;{\mathbb R}^3)}+ \int_0^T \|\tilde w\|^2_{H^3(\Omega_0^f;{\mathbb R}^3)}
]\ ]\ .
\end{align*}
Thus, with the trace estimate (\ref{tracew}),
\begin{align}
&  \int_0^T \|{\tilde w}\|^2_{H^{3}(\Omega_0^f;{\mathbb R}^{3})}+ \int_0^T \|{\tilde q}\|^2_{H^{2}(\Omega_0^f;{\mathbb R})}\ \nonumber\\
&\le [\ C \delta (1+C(M))N(u_0,f)^2+C {\delta_1} (1+T C(M)+N(u_0,f)^2)\ ]\ \|(\tilde w,\tilde q)\|^2_{Z_T} \nonumber\\
&\qquad+ C_{\delta_1} N(u_0,f)^2 +C_{\delta}C(M)T^{\frac{1}{4}}\ (\|(\tilde w,\tilde q)\|^2_{Z_T}+N(u_0,f)^2) +C_{\delta_1} T^{\frac{1}{4}} \|(\tilde w,\tilde q)\|^2_{Z_T}\ . 
\label{tracew1}
\end{align}

Similarly, elliptic regularity for the elasticity problem (for $t\in [0,T]$)
\begin{align*}
 -[c^{ijkl} \int_0^t\tilde w^k,_l],_j &= -{\tilde w}_{t} +f\ \text{in}\ \Omega_0^s\\
\int_0^t \tilde w (t,\cdot)&=\int_0^t \tilde w (t,\cdot)\ \text{on}\ \Gamma_0=\partial\Omega_0^s\ ,
\end{align*}
implies that
$$\|\int_0^t \tilde w\|_{H^{3}(\Omega_0^s;{\mathbb R}^3)}\le C\ [\ \|-\tilde w_{t}+f\|_{H^1(\Omega_0^s;{\mathbb R})}+\|\int_0^t \tilde w\|_{H^{2.5}(\Gamma_0;{\mathbb R}^3)}\ ]\ ,$$
which with (\ref{tracew}) and (\ref{energywttlimit}) provides the inequality
\begin{align}
&\sup_{[0,T]}\|\int_0^{\cdot} {\tilde w}\|^2_{H^{3}(\Omega_0^s;{\mathbb R}^{3})}\nonumber\\
&\le [\ C \delta (1+C(M))N(u_0,f)^2+C {\delta_1} (1+T C(M)+N(u_0,f)^2)\ ]\ \|(\tilde w,\tilde q)\|^2_{Z_T} \nonumber\\
&\qquad+ C_{\delta_1} N(u_0,f)^2 +C_{\delta}C(M)T^{\frac{1}{4}}\ (\|(\tilde w,\tilde q)\|^2_{Z_T}+N(u_0,f)^2) +C_{\delta_1} T^{\frac{1}{4}} \|(\tilde w,\tilde q)\|^2_{Z_T}\ .
\label{tracew2}
\end{align} 

\subsection{Existence and uniqueness of a smooth solution for the non-regularized system (\ref{linear}).} 

We now infer from (\ref{energywttlimit}), (\ref{tracewt1}), (\ref{tracewt2}),
(\ref{tracew1}) and (\ref{tracew2}) that
\begin{align*}
&   \|(\tilde w,\tilde q)\|^2_{Z_T}\\ 
&\qquad\le [\ C \delta (1+C(M)+N(u_0,f)^2)+C {\delta_1} (1+T C(M)+N(u_0,f)^2)\ ]\ \|(\tilde w,\tilde q)\|^2_{Z_T} \nonumber\\
&\qquad\ + C_{\delta_1} N(u_0,f)^2 +C_{\delta}C(M)T^{\frac{1}{4}}\ (\|(\tilde w,\tilde q)\|^2_{Z_T}+N(u_0,f)^2) +C_{\delta_1} T^{\frac{1}{4}} \|(\tilde w,\tilde q)\|^2_{Z_T},
\end{align*}
this inequality being independent of the smoothing parameter of $\tilde a$.

We will call the constant $C$ in this inequality $C_1$ to indicate that at 
this stage it is a constant given by our successive estimates which, for the
sake of conciseness, we have yet to make explicit.

First, we  fix $\delta_1$ so that
$$C_1\delta_1\le\frac{1}{8}\ 
\text{and}\ C_1\delta_1 N(u_0,f)^2\le\frac{1}{8}.$$ 
The constant $C_{\delta_1}$ becomes thus determined, and we have that
\begin{align*}
\|(\tilde w,\tilde q)\|^2_{Z_T} \le&\ [C_1\delta\ (1+C(M))N(u_0,f)^2+ C_1 \delta_1 T C(M)] \ \|(\tilde w,\tilde q)\|^2_{Z_T}\nonumber\\
& + [C_{\delta_1}+C_{\delta}C(M)]\ T^{\frac{1}{4}}\|(\tilde w,\tilde q)\|^2_{Z_T} + C_{\delta_1}  N(u_0,f)^2 \\
& + C_{\delta} C(M) N(u_0,f)^2 T^{\frac{1}{4}} +\frac{1}{4}\|(\tilde w,\tilde q)\|^2_{Z_T} .
\end{align*}

Now let 
\begin{equation}
\label{M}
M=\sup(M_0,\ 2\ [C_1+C_{\delta_1}]\ N(u_0,f)^2)\ .
\end{equation}
Consequently, $C(M)$ becomes a fixed constant. Now, let us fix $\delta>0$ small
enough so that  
\begin{align*}
  \|(\tilde w,\tilde q)\|^2_{Z_T} 
\le &\ [\frac{1}{8}+ C_1\delta_1 T C(M)+[C_{\delta_1}+C_{\delta} C(M)] T^{\frac{1}{4}}]\ \|(\tilde w,\tilde q)\|^2_{Z_T}\\
& + C_{\delta}C(M)\ T^{\frac{1}{4}}\ N(u_0,f)^2  + C_{\delta_1}\ N(u_0,f)^2 +\frac{1}{4}\|(\tilde w,\tilde q)\|^2_{Z_T}\  .
\end{align*}
Now, let $T\in (0,T_M)$ be small enough so that
\begin{align*}
  \frac{3}{4} \|(\tilde w,\tilde q)\|^2_{Z_T} \le& \ \frac{1}{4} \ \|(\tilde w,\tilde q)\|^2_{Z_T} +  C_1\ N(u_0,f)^2 + C_{\delta_1}\ N(u_0,f)^2,
\end{align*}
which implies 
\begin{equation}
\label{stabilityconvex}
 \|(\tilde w,\tilde q)\|^2_{Z_T}\le M\ .
\end{equation}

Henceforth, we revert to our original notation, denoting
$\tilde w$ and  $\tilde a$ by the sequential notation $w_n$ and $a_n$, 
respectively. 
The uniform bound (\ref{stabilityconvex}) ensures the existence of a 
weakly convergent subsequence $(w_{\sigma(n)},q_{\sigma(n)})$ in the 
reflexive Hilbert space $Y_T$ such that 
$$(w_{\sigma(n)},q_{\sigma(n)})\rightharpoonup (w,q)\ \text{in}\ Y_T\ .$$
The usual compactness arguments then provide the strong convergence 
\begin{equation}
\label{strongnonregular}
 (w_{\sigma(n)},q_{\sigma(n)})\rightarrow (w,q)\ \text{in}\ L^2(0,T;H^2(\Omega_0^f;{\mathbb R}^3))\times L^2(0,T;H^1(\Omega_0^f;{\mathbb R}))\ .
\end{equation}
 Combined with the strong convergence
$$ a_n\rightarrow a\ \text{in}\ L^2(0,T;H^2(\Omega_0^f;{\mathbb R}^9))\ $$
(which follows from the mollification process),
the Sobolev embeddings provide the strong convergence
\begin{align*}
a_n a_n^T \nabla w_n &\rightarrow a a^T \nabla w\ \text{in}\ L^2(0,T;L^2(\Omega_0^f;{\mathbb R}^9))\ ,\\
a_n^T q_n&\rightarrow a^T q\ \text{in}\ L^2(0,T;L^2(\Omega_0^f;{\mathbb R}^9))\ .
\end{align*}

We then deduce from (\ref{weakW}) that for each $\phi \in L^2(0,T; H^1_0(\Omega_0^f;{\mathbb R}^3)) $,	
 \begin{align}
& \int_0^T ( w_t, \phi)_{L^2(\Omega;{\mathbb R}^3)}\ dt  
+\nu \int_0^T (a_k^r   w,_r,\  a_k^s \phi,_s)_{L^2(\Omega_0^f;{\mathbb R}^3)}\ dt \nonumber\\
&+ \int_0^T (c^{ijkl}\int_0^t w^k,_l, {\phi}^i,_j)_{L^2(\Omega_0^s;{\mathbb R})}\ dt  + \int_0^T (q,\ a_k^l \phi^k,_l)_{L^2(\Omega_0^f;{\mathbb R})} dt\nonumber\\
&\qquad\qquad = \int_0^T (f\circ\eta, \phi)_{L^2(\Omega_0^f;{\mathbb R}^3)}+ (f, \phi)_{L^2(\Omega_0^s;{\mathbb R}^3)}  \ dt \ ,
\label{weakwlimit}
\end{align}
which combined with the fact that, from (\ref{strongnonregular}), for 
all $t\in [0,T]$, $w(t)\in {\mathcal V}_v(t)$, proves that $w$ is a weak solution 
of (\ref{linear}).

Since $w\in L^2(0,T;H^1_0(\Omega;{\mathbb R}^3))$ we infer the uniqueness of a solution in $Y_T$ to this system in the same classical fashion as for the solution $\tilde w$ of the regularized problem.

Moreover, it is also immediate that we have from (\ref{stabilityconvex}) the estimate
\begin{equation}
\label{stabilityconvexbis}
 \|w\|^2_{W_T}\le M\ .
\end{equation}
This concludes the proof of Theorem \ref{thm1}.

Henceforth,  $M$ is given by (\ref{M}) and $T$ is chosen such that
(\ref{stabilityconvex}) holds.

\section{The fixed-point scheme for the nonlinear problem}
\label{10}

We will make use of the Tychonoff Fixed-Point Theorem  in our fixed point procedure (see, for example, \cite{Deimling}).  Recall that
this states that
for a reflexive separable Banach space $X$, and $C\subset X$ a closed,
convex, bounded subset,  if $F:C \rightarrow C$ is weakly
sequentially continuous into $X$, then $F$ has at least one fixed-point.

With the  quantities $M$ and $T$ being defined as in the previous section, we 
make the following
\begin{definition}
We define a mapping $\Theta_T$ from $C_T (M)$ into itself (from estimate (\ref{stabilityconvexbis})), which to a given
element $v \in  W_T $ associates $w\in  W_T$, the unique solution in $Y_T$ of
(\ref{linear}).
\end{definition}

We next have the following weak sequential continuity result.

\begin{lemma}
\label{weakcontinuity}
The mapping $\Theta_T$ associating $w$ to $v\in C_T (M)$ is weakly sequentially continuous from
$C_T (M)$ into $C_T (M)$ (endowed with the norm of $X_T$).
\end{lemma}
\begin{proof}
Let $(v_p)_{p\in\mathbb N}$ be a given sequence of elements of $C_T (M)$
weakly convergent (in $X_T$) toward a given element $v\in C_T (M)$ (
$C_T (M)$ is sequentially weakly closed as a closed convex set) and let  $(v_{\sigma(p)})_{p\in\mathbb N}$ be
any subsequence of this sequence.

Since $V_f^3(T)$ is compactly embedded into
$L^2(0,T;H^2(\Omega_0^f;{\mathbb R}^3))$, 
we deduce the following strong convergence results in 
$L^2(0,T;L^{2}(\Omega_0^f;{\mathbb R}))$  as $p$ goes to $\infty$:
\begin{subequations}
\label{scv}
\begin{align}
(a^j_l)_p (a^k_l)_p &\rightarrow a^j_l a^k_l\,, 
\label{scv.a}\\
[(a^j_l)_p (a^k_l)_p]_{,j} &\rightarrow (a^j_l a^k_l)_{,j}\,,
\label{scv.b}\\
(a^k_i)_p &\rightarrow a^k_i \,, \\
f^i\circ\eta_p &\rightarrow f^i\circ\eta\ .
\end{align}
\end{subequations}

Now, let $w_p=\Theta_T (v_p)$ and let $q_p$ be the associated pressure, so that $(q_p)_{p\in\mathbb N}$ is in a bounded set of $V^2_f (T)$.
Since $Y_T=X_T\times V^2_f (T)$ is a reflexive Hilbert space, let 
$(w_{\sigma'(p)},q_{\sigma'(p)})_{p\in\mathbb N}$ be a subsequence 
weakly converging in $Y_T$ toward an element $(w, q) \in Y_T$. 
Since $C_T(M)$ is weakly closed in $X_T$, we also have $w\in C_T(M)$.

We can then infer in a similar fashion as for the proof of (\ref{weakwlimit}) in the previous section that for each $\phi \in L^2(0,T; H^1_0(\Omega_0^f;{\mathbb R}^3)) )$,	
 \begin{align*}
& \int_0^T ( w_t, \phi)_{L^2(\Omega;{\mathbb R}^3)}\ dt  
+\nu \int_0^T (a_k^r   w,_r,\  a_k^s \phi,_s)_{L^2(\Omega_0^f;{\mathbb R}^3)}\ dt \nonumber\\
&+ \int_0^T (c^{ijkl}\int_0^t w^k,_l, {\phi}^i,_j)_{L^2(\Omega_0^s;{\mathbb R})}\ dt  + \int_0^T (q,\ a_k^l \phi^k,_l)_{L^2(\Omega_0^f;{\mathbb R})} dt\nonumber\\
&\qquad\qquad = \int_0^T (f\circ\eta, \phi)_{L^2(\Omega_0^f;{\mathbb R}^3)}+ (f, \phi)_{L^2(\Omega_0^s;{\mathbb R}^3)}  \ dt \ ,
\end{align*}
which combined with the fact that, from (\ref{scv}), for all $t\in [0,T]$, $w(t)
\in {\mathcal V}_v(t)$, shows that $w$ is a weak solution of (\ref{linear}) 
in $C_T(M)$, {\it i.e.} $w=\Theta_T (v)$.

Therefore, we deduce that the whole sequence 
${(\Theta_T (v_n))}_{n\in\mathbb N}$ weakly converges
in $C_T (M)$ toward $ \Theta_T (v) $, which concludes the proof.
\end{proof} 

\section{Proof of theorem \ref{main}}
\label{11}

The mapping $\Theta$ being weakly continuous from the closed bounded convex
set $C_T(M)$ into itself from Lemma \ref{weakcontinuity}, we
infer from the Tychonoff fixed point theorem ({\it see} for instance 
\cite{Deimling}) that it admits (at least) one fixed point 
$v=\Theta(v)$ in $C_T$. Moreover, since $T\le T_M$, Lemma \ref{collision} ensures us that there is no collision between solids or between a solid and $\partial\Omega$. Thus, 
$(v,q)$ is a solution of (\ref{nsl}).
Note that the continuity of the Lagrangian velocities $v^f=v^s$ at the 
interface $\Gamma_0$ is ensured by our functional framework, since
$(v,q)\in X_T$ implies $v\in L^2 (0,T; H^1_0(\Omega;{\mathbb R}^3))$, 
which provides the equality 
$v^f=v^s$ in $H^{\frac{1}{2}}(\Gamma_0; {\mathbb R}^3)$. 

\section{Uniqueness}
\label{12}

Uniqueness will be obtained under stronger assumptions than the ones used to establish existence, for reasons that will be explained hereafter.

If $(\tilde v,\tilde q)\in Y_T$ is another solution of (\ref{nsl}), then, 
\begin{subequations}
\label{uniquesystem}
\begin{align}
({{v-\tilde{v}}})^i_t - \nu (a^j_l a^k_l ({{v}}^i,_k-{\tilde{v}}^i,_k )),_j
 + a^k_i ({q},_k-\tilde{q},_k ) 
&= \delta f^i\ \ \text{in} \ \ (0,T)\times \Omega_0^f \,, \\
 a^k_i {(v-\tilde{v})}^i,_k &= \delta a  \ \ \text{in}
 \ \ (0,T)\times \Omega_0^f \,, \\
 (v-\tilde v)^i_t - [c^{ijkl}\int_0^t (v-\tilde v)^k,_l],_j &= 0\ \ \text{in} \ \ (0,T)\times \Omega_0^s \,, \\
\nu (v^f-\tilde v^f)^i,_k a^k_l a_l^j N_j-(q-\tilde q) a_i^j N_j &=
c^{ijkl} \int_0^t ({{v^s-\tilde{v}^s}})^k,_l\ N_j \nonumber\\
&\qquad +\delta g^i\ \text{on} \ \ (0,T)\times \Gamma_0 \,, \\
  {v}-\tilde{v} &= 0    
 \ \ \text{on} \ \ (0,T)\times \partial\Omega \,,\\
  {v}-\tilde{v} &= 0    
 \ \ \text{on} \ \ \{0\}\times \Omega \,, 
\end{align}
\end{subequations}
with
\begin{subequations}
\label{movingforces}
\begin{align}
\delta f^i&=- \nu ((a^j_l a^k_l -\tilde{a}^j_l \tilde{a}^k_l)\ {\tilde{v}}^i,_k )),_j + f\circ\eta-f\circ\tilde{\eta} +(-a_i^k+\tilde a_i^k)\tilde q,_k
\ \text{in} \ \ (0,T)\times \Omega_0^f \,, \label{mf.a} \\
\delta a&= (\tilde{a}^k_i-a^k_i) {\tilde{v}}^i,_k 
\ \text{in} \ \ (0,T)\times \Omega_0^f \,, \label{mf.b}\\
\delta g^i&= \nu (\tilde v^f,_k^i \tilde a^k_l \tilde a_l^j N_j- \tilde v^f,_k^i a^k_l a_l^j N_j) - \tilde q (\tilde a_i^j-a^j_i) N_j\ \text{on} \ \ (0,T)\times \Gamma_0 . 
\label{u6}
\end{align}
\end{subequations}

If we view this problem with $v-\tilde v$ as the unknown velocity and $q-\tilde q$ as the associated pressure in the fluid, this problem looks similar to
the linear problem (\ref{linear}); it is tempting to conclude that similar
estimates as in the study of the regularity of (\ref{linear})
would yield a differential inequality that would provide uniqueness.
It appears, however, that such a procedure fails because  due to the Dirichlet
boundary condition on $\partial \Omega$ for the velocity, we are not
able to get the necessary information on 
the second derivative of the pressure function.
Such information is crucial since $\delta f_{tt}$ contains $\tilde q_{tt}$, 
which makes the second time differentiated problem impossible to estimate.

For this reason, we will need to impose more regularity on the data and 
forcing functions, so that we, in turn, have enough information on 
$\tilde q_{tt}$, which will then be viewed as a coefficient in the study of 
the regularity of (\ref{uniquesystem}) in $Y_T$.

We first update the functional framework. Let us denote 
\begin{align*}
V^4_f (T)=\{u\in L^2(0,T;H^4(\Omega_0^f;{\mathbb R}^3))|\ u_t\in V^3_f (T)\},\\
V^4_s (T)=\{u\in L^2(0,T;H^4(\Omega_0^s;{\mathbb R}^3))|\ u_t\in V^3_s (T)\}\ .
\end{align*}
Let us the denote the reflexive separable Hilbert space 
$$X_T^u=\{u\in L^2(0,T;H^1_0(\Omega;{\mathbb R}^3))|\ (u^f,\int_0^{\cdot} u^s)\in V^4_f(T)\times V^4_s(T)\}\ ,$$
endowed with its natural norm $$\|u\|^2_{X_T^u}=\|u_t\|^2_{X_T}+\|u\|^2_{L^2(0,T;H^4(\Omega_0^f;{\mathbb R}^3))}+ \|\int_0^{\cdot} u\|^2_{L^2(0,T;H^4(\Omega_0^s;{\mathbb R}^3))}\ .$$ 
In a similar fashion, we introduce
\begin{align*}
Y_T^u=\{(u,p)|\ &u\in X_T^u,\ p\in L^2(0,T;H^3(\Omega_0^f;{\mathbb R})),\
p_t\in L^2(0,T;H^2(\Omega_0^f;{\mathbb R})),\\
& p_{tt}\in L^2(0,T;H^1(\Omega_0^f;{\mathbb R}))
\}\ ,
\end{align*}
endowed with its natural norm 
\begin{align*}
\|(u,p)\|^2_{Y_T^u}=&\|u\|^2_{X_T^u}+\|p\|^2_{L^2(0,T;H^3(\Omega_0^f;{\mathbb R}))} + \|p_{t}\|^2_{L^2(0,T;H^2(\Omega_0^f;{\mathbb R}))}\\
&+ \|p_{tt}\|^2_{L^2(0,T;H^1(\Omega_0^f;{\mathbb R}))} .
\end{align*} 

We will also need 
\begin{align*}
W_T^u=\{u\in X_T^u|\ &u_{ttt}\in L^{\infty}(0,T;L^2(\Omega;{\mathbb R}^3)),\ \int_0^{\cdot} u\in L^{\infty}(0,T;H^4(\Omega_0^s;{\mathbb R}^3)\\
&  u\in L^{\infty}(0,T;H^3(\Omega_0^s;{\mathbb R}^3),\ u_t \in L^{\infty}(0,T;H^2(\Omega_0^s;{\mathbb R}^3)\\
& u_{tt}\in L^{\infty}(0,T;H^1(\Omega_0^s;{\mathbb R}^3)\}\ ,
\end{align*}
endowed with its natural norm 
\begin{align*}
\|u\|^2_{W_T^u}=\ &\|u\|^2_{X_T^u}+\|u_{ttt}\|^2_{L^{\infty}(0,T;L^2(\Omega;{\mathbb R}^3))}+\|\int_0^{\cdot} u\|^2_{L^{\infty}(0,T;H^4(\Omega_0^s;{\mathbb R}^3))}\\
&+\|u\|^2_{L^{\infty}(0,T;H^3(\Omega_0^s;{\mathbb R}^3))}+\| u_t\|^2_{L^{\infty}(0,T;H^2(\Omega_0^s;{\mathbb R}^3))}+ \|u_{tt}\|^2_{L^{\infty}(0,T;H^1(\Omega_0^s;{\mathbb R}^3))}\ ,
\end{align*}
as well as $Z_T^u (T)=\{(u,p)\in Y_T^u|\ u\in W_T^u\} $ endowed with its natural norm, 
\begin{align*}
\|(u,p)\|^2_{Z_T^u}=&\|u\|^2_{W_T^u}+\|p\|^2_{L^2(0,T;H^3(\Omega_0^f;{\mathbb R}))} + \|p_{t}\|^2_{L^2(0,T;H^2(\Omega_0^f;{\mathbb R}))}\\
&+ \|p_{tt}\|^2_{L^2(0,T;H^1(\Omega_0^f;{\mathbb R}))} .
\end{align*}

We can then define the convex set $C_M^u(T)$ in the same fashion as $C_M(T)$, with $W_T^U$ replacing $W_T$ and with the additional condition $w_{tt}(0)=w_2$ where $w_2$ has been defined in (\ref{w2def}).

We can then prove, in a way similar to the proof of Theorem \ref{thm1} (with the introduction of the penalized problems, time differentiated three times now) that the following holds:

\begin{theorem}\label{existenceupdated}
With the same assumptions as in Theorem \ref{main} and under the supplementary 
conditions $\Omega$ of class $H^4$, $\Omega_0^s$ of class $H^5$, the initial 
data $u_0 \in H^7(\Omega_0^f;{\mathbb R}^3)\cap 
H^4(\Omega_0^s;{\mathbb R}^3)\cap H^1_0(\Omega;{\mathbb R}^3)\cap 
L^2_{{div},f}$, $f(0)\in H^5(\Omega;{\mathbb R}^3)$,  
satisfying the supplementary compatibility conditions (recall the assumption 
of Section \ref{1bis}) 
\begin{align*}
 [ (\nu [\nabla w_2^f\ N]^i + 2\nu w_1^f,_k^i ({a^k_l}a_l^j)_t (0) N_j +& \nu u_0^f,_k^i ({a^k_l} a_l^j)_{tt} (0))_{i=1}^3 ]_{\operatorname{tan}}\\
&= [\ (c^{ijkl} w_1^s,_l^k N_j)_{i=1}^3\ ]_{\operatorname{tan}}\ \text{on}\ \Gamma_0,\\
% \ \operatorname{div} w_2+2{a_i^j}_t (0) {w_1^i},_j+ a_i^j u_0^i,_j&=0\ \text{in}\ \Omega_0^f,\\
 w_2&\in H^1_0(\Omega;{\mathbb R}^3) ,
\end{align*}
and the supplementary assumption on the forcing function that 
\begin{align*}
f \in L^2(0,\bar T; H^3(\Omega;{\mathbb R}^3)),&\ f_{t} \in L^2(0,\bar T; H^2(\Omega;{\mathbb R}^3)),\ f_{tt} \in L^2(0,\bar T; H^1(\Omega;{\mathbb R}^3))\\
&f_{ttt} \in L^2(0,\bar T; L^2(\Omega;{\mathbb R}^3)),
\end{align*}
 we have the existence of $T>0$ such that there exists a solution $(v,q)\in Y_T^u$ of (\ref{nsl}). Furthermore, $v\in C_M^u(T)$ for $M$ appropriate.
\end{theorem}

We can now get estimates for (\ref{uniquesystem}) which will give an appropriate differential inequality, in the space $Z_T$ used to prove Theorem \ref{thm1}.
We notice that this problem is similar to (\ref{linear}) except for the 
divergence-type condition which is not set to zero, and the boundary forcing 
on the interface. 

The Neumann forcing does not give any specific difficulty, and can be handled 
without modification of our previous estimates. 

The divergence-type condition does not bring any difficulty either because we 
do not need to establish the existence of a solution to (\ref{uniquesystem}), 
since it comes {\it de facto} from the definition of $v$ and $\tilde v$; we
can directly use this condition in the steps where we obtained
$\epsilon$-independent estimates for 
$w_{tt}$, $w_t$ and $w$. We also do not have to regularize the coefficients, 
since the 
regularity of $w$ is a given. Those three steps would provide us in the same 
fashion as for the proof of Theorem \ref{linear} with the appropriate estimates 
to be made precise later. Note that this process works because the right-hand 
side of the divergence condition for $w$ in (\ref{uniquesystem}) has (roughly 
speaking) the term $\nabla\eta-\nabla\tilde\eta$, which has one time 
derivative less than the right-hand side
$\nabla v-\nabla\tilde v$ (the term $\nabla v$ on the right-hand side being 
viewed as a coefficient whose regularity is given by Theorem 
\ref{existenceupdated}).

We are now in  a position to state the uniqueness result.
\begin{theorem}\label{unique}
With the same assumptions as in Theorem \ref{existenceupdated}, and with
the additional assumption that there exists $K>0$ such that 
\begin{equation}\label{Lip}
\begin{array}{l}
\forall t\le \bar T,\  \forall (x,y)\in \Omega\times \Omega, \\
|f(t,x)-f(t,y)|+ |\nabla f(t,x)-\nabla f(t,y)|+|f_t(t,x)-f_t(t,y)|\\
+ |\nabla f_t(t,x)-\nabla f_t(t,y)|+|f_{tt}(t,x)-f_{tt}(t,y)|+
|\nabla^2 f(t,x)-\nabla^2 f(t,y)|\le K\ |x-y|\,,
\end{array}
\end{equation}
i.e, 
$f$, $\nabla f$, $\nabla^2 f$, $f_t$, $\nabla f_t$ and $f_{tt}$ are uniformly Lipschitz continuous in the spatial
variable, then the solution is unique.
\end{theorem}

\begin{proof}
With those assumptions on $f$, we have for the forcing $f\circ\eta-f\circ\tilde\eta$ appearing in (\ref{uniquesystem}) an estimate 
\begin{align*}
\|f\circ\eta-f\circ\tilde\eta\|&_{L^2(0,T;H^2(\Omega_0^f;{\mathbb R}^3))}
+ \|(f\circ\eta-f\circ\tilde\eta)_{t}\|_{L^2(0,T;H^1(\Omega_0^f;{\mathbb R}^3))}\\
&+ \|(f\circ\eta-f\circ\tilde\eta)_{tt}\|_{L^2(0,T;L^2(\Omega_0^f;{\mathbb R}^3))}\le C\ \|\eta-\tilde\eta\|_{V^3_f(T)}\ .
\end{align*}
The other terms associated to $\delta f$, $\delta g$, $\delta a$ have the same
effect in the integral estimates for $w_{tt}$, $w_t$ and $w$. This leads us to
\begin{align}
\forall t\in (0,T),\ \|(v-\tilde v,q-\tilde q)\|_{Z_t}\le C_1\ \|\eta-\tilde\eta\|_{V^3_f(t)}\ .
\end{align}
where the constant $C_1$ depends here on the same variables as the generic constant $C$ as well as on the initial data.
This thus implies $$\forall t\in (0,T),\ \|v-\tilde v\|_{V^3_f (t)}\le C_1\ \|\eta-\tilde\eta\|_{V^3_f(t)},$$
from which we infer 
$$\forall t\in (0,T),\ \|v-\tilde v\|_{V^3_f (t)}\le C_1\ \int_0^t\|v-\tilde v\|_{V^3_f(t)}\le C_1 t\ \|v-\tilde v\|_{V^3_f (t)},$$
 which shows that for $T_1=\frac{1}{2 {C_1}}$, we have $v=\tilde v$ on $[0,T_1]$. We can then iterate this, starting from the initial time set at $T_1$, which gives in a similar fashion, since $v(T_1)=\tilde v(T_1)$, that $v=\tilde v$ on $[T_1,2 T_1]$. By induction, we get 
$v=\tilde v$ on $[0,T]$. 
\end{proof}

\section{concluding remarks}
Whereas the fluid-solid interaction is indeed a moving interface problem, it 
appears that the methods for its analysis differ drastically from the 
classical methods developed for the  Navier-Stokes fluid interfaces 
independently by Solonnikov (see \cite{Sol1992} and references therein) and 
Beale \cite{Beale1983}.

First, our functional framework scales in a {\it hyperbolic} fashion in both 
the parabolic (fluid) and hyperbolic (solid) phases.

Second, whereas the fixed-point problem (\ref{linear}) is inspired by the 
classical fixed-point problem used in parabolic-type interface problems, 
the Fourier transform 
technique used to get regularity in parabolic theories requires the 
introduction of the problem with constant coefficients (for which one does get 
explicit solutions), with the forcing functions containing the difference 
(small in a neighborhood of a point on the interface) between the actual 
coefficient and this constant coefficient. 
Whereas this procedure is contractive for parabolic problems, the hyperbolic 
part is problematic in the sense that the difference between the actual and 
the constant hyperbolic viscosity is not regular enough to get these 
contractive estimates 
(those coefficients are not constant after the truncation and change of 
variables to the full space problem). 

Third, whereas energy methods without the use of Fourier techniques are indeed 
known for incompressible fluid interfaces (see for instance \cite{CoSh2003}), 
the highest-order time derivative of the pressure is known in 
$L^2(0,T;L^2(\Omega;{\mathbb R}^3))$ in that case, which allows the use of an 
iterative method from the constant-coefficient problem in energy spaces similar
to the ones described in \cite{Sol1992}. In the fluid-solid problem, the 
knowledge of the highest order time derivative of the pressure is not known, 
which prevents such an iterative procedure from the constant-coefficient 
problem to get regularity. Instead, we are forced to work directly with the
Lagrangian formulation (\ref{linear}), which requires the introduction of
the penalized problems for reasons explained previously about the pressure. 
In turn, working with the Lagrangian formulation (\ref{linear}) requires us
to first smooth the coefficients, and then to obtain estimates 
independent of the smoothing  parameter by using interpolation inequalities. 

Fourth, we clearly identify in our method the central and sufficient role of 
the {\it trace} of the velocity on the material  interface $\Gamma_0$, 
whereas classical
regularity results in interface problems involve the study of the regularity 
in the interior.

Fifth, once again regarding  the pressure estimate, obtaining  a 
contractive fixed-point scheme does not seem possible for the 
hyperbolic-parabolic 
problem (even with data arbitrarily smooth), whereas it is indeed the most 
well-known method for the parabolic interface case. Note, however, that this 
later point is associated to the {\it incompressibility} of the fluid and 
does not seem to appear without this constraint. 

This last remark is not without consequences for the numerical analysis of the
problem, which we shall  develop in future work. As for the question of the 
convergence of solutions of certain regularized models considered by other 
authors, it seems that the evidently natural approach of taking an 
elasticity law with a finite number $N$ of modes introduced in \cite{DeEsGrLe}
and letting $N \rightarrow \infty$ leads to some difficulties, as there is no 
elliptic
operator for the discrete elasticity problem for which one may use $H^3$ 
regularity independently of the number of modes. On the other hand, it can be 
shown that the addition of a hyperviscosity to the solid problem 
(similar in spirit to the hyperviscous plate problem introduced in 
\cite{ChDeEsGr}) would indeed converge to the solution of the actual  problem 
as the hyperviscosity parameter tends to zero, since we can apply 
the methods constructed here to this family of problems and obtain estimates 
that are independent of the hyperviscous parameter in the correct norms.  

\appendix
\section{Some additional estimates}
\label{13}
\subsection{Estimates for (\ref{energywtH2})}
In this section, $\delta>0$ is assumed given and we now proceed to the estimate
of the terms of (\ref{energywtH2}) whose sign is indefinite. Recall that from
(\ref{c0}), 
$\tilde a$, $\tilde a_t$, and thus $\tilde b$, $\tilde b_t$ are controlled respectively in $L^{\infty}(H^2)$ and $L^{\infty}(H^1)$, 
independently of the regularizing parameter $n$ associated to $\tilde a$. 

\noindent {\bf Step 1.} 
Let $\displaystyle J_1= \int_0^T (D_h \tilde b_i^j\ Q_t,\  D_h W_t^i,_j)_{L^2({\mathbb R}^3_+;{\mathbb R})}$. Then,
\begin{align*}
|J_1|&\le \delta \int_0^T \|D_h \nabla W_{t}\|^2_{L^2 ({\mathbb R}^3_+;{\mathbb R}^9)} +C_{\delta}\int_0^T \|D_h\tilde b\|^2_{L^4({\mathbb R}^3_+;{\mathbb R}^9)}\|Q_t\|^2_{L^{4} ({\mathbb R}^3_+;{\mathbb R})}\nonumber\\
&\le \delta \int_0^T \|D_h\nabla \tilde W_{t}\|^2_{L^2 ({\mathbb R}^3_+;{\mathbb R}^9)} +C_{\delta} C(M)\int_0^T \|Q_t\|^{0.5}_{L^2 ({\mathbb R}^3_+;{\mathbb R})} \|Q_t\|^{1.5}_{H^1 ({\mathbb R}^3_+;{\mathbb R})}\nonumber\\
&\le  \delta \int_0^T \|D_h \nabla W_{t}\|^2_{L^2 ({\mathbb R}^3_+;{\mathbb R}^9)}\nonumber\\
&\qquad +C_{\delta} C(M)\ T^{\frac{1}{4}}\ [\sup_{(0,T)}\|Q_t\|^{2}_{L^2 ({\mathbb R}^3_+;{\mathbb R})} + \int_0^T \|Q_t\|^{2}_{H^1 ({\mathbb R}^3_+;{\mathbb R})} ]\ .
\end{align*}
From (\ref{eliminateqt}) and the definitions of $W$ and $Q$, we then infer
\begin{align}
|J_1|\le &\ C\delta \int_0^T \|\tilde w_t\|^2_{H^2 (\Omega_0^f;{\mathbb R}^3)}\nonumber\\
& +C_{\delta} C(M)\ T^{\frac{1}{4}}\ [\ N(u_0,f)^2+ T \int_0^T \|\tilde q_t\|^{2}_{H^1 (\Omega_0^f;{\mathbb R})} +\sup_{[0,T]} \|\tilde w_{tt}\|^{2}_{L^2 (\Omega;{\mathbb R}^3)}\nonumber\\
&\qquad\qquad\qquad\qquad + T\ [\int_0^T\|\tilde w_t\ \|^{2}_{H^2 (\Omega_0^f;{\mathbb R}^3)} + \int_0^T\|\tilde w_{tt}\|^{2}_{H^1 (\Omega_0^f;{\mathbb R}^3)}]\nonumber\\
&\qquad\qquad\qquad\qquad+ \sup_{[0,T]}\|\tilde w\|^{2}_{H^1 (\Omega_0^s;{\mathbb R}^3)} + \int_0^T \|\tilde q_t\|^{2}_{H^1 (\Omega_0^f;{\mathbb R})}\ ]\ .
\label{J1} 
\end{align}

\noindent {\bf Step 2.} Let $\displaystyle J_2= \int_0^T (D_h \tilde {b_t}_i^j\ Q,\  D_h W_t^i,_j)_{L^2({\mathbb R}^3_+;{\mathbb R})}$. Similarly as for $J_1$,
\begin{align}
|J_2|&\le \delta \int_0^T \|D_h \nabla W_{t}\|^2_{L^2 ({\mathbb R}^3_+;{\mathbb R}^9)} +C_{\delta}\int_0^T \|D_h\tilde b_t\|^2_{L^2({\mathbb R}^3_+;{\mathbb R}^9)}\|Q\|^2_{W^{1,4} ({\mathbb R}^3_+;{\mathbb R})}\nonumber\\
&\le \delta \int_0^T \|D_h\nabla  W_{t}\|^2_{L^2 ({\mathbb R}^3_+;{\mathbb R}^9)} +C_{\delta} C(M)\int_0^T \|Q\|^{0.5}_{H^1 ({\mathbb R}^3_+;{\mathbb R})} \|Q\|^{1.5}_{H^2 ({\mathbb R}^3_+;{\mathbb R})}\nonumber\\ 
%&\le \delta \int_0^T \|D_h \nabla W_{t}\|^2_{L^2 ({\mathbb R}^3_+;{\mathbb R}^9)} +C_{\delta} C(M)\ T^{\frac{1}{4}}\ [\sup_{(0,t)}\|Q\|^{2}_{H^1 ({\mathbb R}^3_+;{\mathbb R})} + \int_0^T \|Q\|^{2}_{H^2 ({\mathbb R}^3_+;{\mathbb R})} ]\nonumber\\
&\le  C\delta \int_0^T \|\tilde w_{t}\|^2_{H^2 (\Omega_0^f;{\mathbb R}^3)} \nonumber\\
& +C_{\delta} C(M)\ T^{\frac{1}{4}}\ [\ N(u_0,f)^2+ T\int_0^T \|\tilde q_t\|^{2}_{H^1 (\Omega_0^f;{\mathbb R})} + \int_0^T \|\tilde q\|^{2}_{H^2 (\Omega_0^f;{\mathbb R})} ]\ .
\label{J2} 
\end{align}

\noindent {\bf Step 3.} Let $\displaystyle J_3= \int_0^T (\tilde {b_t}_i^j\ D_h Q,\  D_h W_t^i,_j)_{L^2({\mathbb R}^3_+;{\mathbb R})}$. Then,
\begin{align*}
|J_3|&\le \delta \int_0^T \|D_h \nabla W_{t}\|^2_{L^2 ({\mathbb R}^3_+;{\mathbb R}^9)} +C_{\delta}\int_0^T \|\tilde b_t\|^2_{L^4({\mathbb R}^3_+;{\mathbb R}^9)}\|D_h Q\|^2_{L^{4} ({\mathbb R}^3_+;{\mathbb R})}\ .
\end{align*}
 Thus, similarly as for (\ref{J2}),
\begin{align}
|J_3|\le &\ C\delta \int_0^T \|\tilde w_{t}\|^2_{H^2 (\Omega_0^f;{\mathbb R}^3)} \nonumber\\
& +C_{\delta} C(M)\ T^{\frac{1}{4}}\ [\ N(u_0,f)^2+ T\int_0^T \|\tilde q_t\|^{2}_{H^1 (\Omega_0^f;{\mathbb R})} + \int_0^T \|\tilde q\|^{2}_{H^2 (\Omega_0^f;{\mathbb R})} ]\ .
\label{J3} 
\end{align}

\noindent {\bf Step 4.} Let $\displaystyle J_4= \int_0^T ( \tilde {b}_i^j\ D_h Q_t,\  D_h W_t^i,_j)_{L^2({\mathbb R}^3_+;{\mathbb R})}$. This term will require more care.
We first notice that
$$J_4=\int_0^T (  D_h Q_t,\  D_h [\tilde {b}_i^j W_t^i,_j])_{L^2({\mathbb R}^3_+;{\mathbb R})}- \int_0^T (  D_h Q_t,\  D_h \tilde {b}_i^j\  W_t^i,_j (\cdot+h))_{L^2({\mathbb R}^3_+;{\mathbb R})}\ ,$$
which with the divergence relation (\ref{divchange}) leads us to
\begin{align}
J_4=&\int_0^T (  D_h Q_t,\  D_h \mathfrak a_t)_{L^2({\mathbb R}^3_+;{\mathbb R})}-\int_0^T (  D_h Q_t,\  D_h [\tilde {b_t}_i^j W^i,_j])_{L^2({\mathbb R}^3_+;{\mathbb R})}\nonumber\\
& - \int_0^T (  D_h Q_t,\  D_h \tilde {b}_i^j\  W_t^i,_j (\cdot+h))_{L^2({\mathbb R}^3_+;{\mathbb R})}\ .
\label{J41}
\end{align}

For the first integral of
this identity, $\displaystyle J_4^1=\int_0^T (  D_h Q_t,\  D_h \mathfrak a_t)_{L^2({\mathbb R}^3_+;{\mathbb R})}$, we have
\begin{align}
|J_4^1|& \le  \delta \int_0^T \|D_h Q_{t}\|^2_{L^2 ({\mathbb R}^3_+;{\mathbb R})} +C_{\delta} \int_0^T \|D_h\tilde b_t\|^2_{L^2({\mathbb R}^3_+;{\mathbb R}^9)}\|\tilde w\|^2_{W^{1,4} (\Omega_0^f;{\mathbb R}^3)}\nonumber \\
&\qquad+C_{\delta}\int_0^T \|D_h\tilde b\|^2_{L^4({\mathbb R}^3_+;{\mathbb R}^9)}\|\tilde w_t\|^2_{L^{4} (\Omega_0^f;{\mathbb R}^3)}+C_{\delta}\int_0^T \|\tilde b\|^2_{W^{1,4}({\mathbb R}^3_+;{\mathbb R}^9)}\|\nabla \tilde w_t\|^2_{L^2(\Omega_0^f;{\mathbb R}^9)}\nonumber \\
&\qquad+C_{\delta}\int_0^T \|\tilde b_t\|^2_{L^4({\mathbb R}^3_+;{\mathbb R}^9)}\|\nabla \tilde w\|^2_{L^4 (\Omega_0^f;{\mathbb R}^9)}\ ]\nonumber\\
&\le \delta \int_0^T \|D_h Q_{t}\|^2_{L^2 ({\mathbb R}^3_+;{\mathbb R})} +C_{\delta} C(M)\ T\ [ \sup_{[0,T]} \|\tilde w\|^{2}_{H^2 (\Omega_0^f;{\mathbb R}^3)} + \sup_{[0,T]} \|\tilde w_t\|^{2}_{H^1 (\Omega_0^f;{\mathbb R}^3)} ]\ ,\nonumber\\
&\le C\delta \int_0^T \|\tilde q_t\|^2_{H^1 (\Omega_0^f;{\mathbb R})} \nonumber\\
&\qquad +C_{\delta} C(M)\ T\ [N(u_0,f)^2+ T\int_0^T \|\tilde w_t\|^{2}_{H^2 (\Omega_0^f;{\mathbb R}^3)} + T\int_0^T \|\tilde w_{tt}\|^{2}_{H^1 (\Omega_0^f;{\mathbb R}^3)} ]\ .
\label{J42} 
\end{align}
Next, for $\displaystyle J_4^2=\int_0^T (  D_h Q_t,\  D_h [\tilde {b_t}_i^j W^i,_j])_{L^2({\mathbb R}^3_+;{\mathbb R})}$,
\begin{align*}
|J_4^2|\le &\ \delta \int_0^T \|D_h Q_{t}\|^2_{L^2 ({\mathbb R}^3_+;{\mathbb R})} +C_{\delta}\int_0^T \|D_h\tilde b_t\|^2_{L^2({\mathbb R}^3_+;{\mathbb R}^9)}\|\nabla W\|^2_{W^{1,4} ({\mathbb R}^3_+;{\mathbb R}^9)}\\
&+C_{\delta}\int_0^T \|\tilde b_t\|^2_{L^4({\mathbb R}^3_+;{\mathbb R}^9)}\|D_h \nabla W\|^2_{L^{4} ({\mathbb R}^3_+;{\mathbb R}^9)} \ .
\end{align*}
Therefore, 
\begin{align*}
|J_4^2|\le\ & \delta \int_0^t \|D_h Q_{t}\|^2_{L^2 ({\mathbb R}^3_+;{\mathbb R})}\\
& +C_{\delta} C(M)\ T^{\frac{1}{4}}\ [ \sup_{[0,T]} \|\nabla W\|^{2}_{H^1 ({\mathbb R}^3_+;{\mathbb R}^9)} + \int_0^T \|\nabla W\|^{2}_{H^2 ({\mathbb R}^3_+;{\mathbb R}^9)} ]\ ,
\end{align*}
which with the definition of $W$ and $Q$ provides
\begin{align}
|J_4^2|\le\ & C\delta \int_0^t \|\tilde q_{t}\|^2_{H^1 (\Omega_0^f;{\mathbb R})}\nonumber\\
& +C_{\delta} C(M)\ T^{\frac{1}{4}}\ [N(u_0,f)^2+ T\int_0^T \|\nabla \tilde w_t\|^{2}_{H^1 (\Omega_0^f;{\mathbb R}^9)} + \int_0^T \|\nabla\tilde w\|^{2}_{H^2 (\Omega_0^f;{\mathbb R}^9)} ]\ .
\label{J43} 
\end{align}

Similarly, for $\displaystyle J_4^3=\int_0^T (  D_h Q_t,\  D_h \tilde {b}_i^j\ {W_t}^i,_j (\cdot+h))_{L^2({\mathbb R}^3_+;{\mathbb R})}$,
\begin{align}
|J_4^3|&\le \ \delta \int_0^T \|D_h Q_{t}\|^2_{L^2 ({\mathbb R}^3_+;{\mathbb R})} +C_{\delta}\int_0^T \|D_h\tilde b\|^2_{L^4({\mathbb R}^3_+;{\mathbb R}^9)}\|\nabla W_t\|^2_{L^4 ({\mathbb R}^3_+;{\mathbb R}^9)}\nonumber\\
%& \le\  \delta \int_0^T \|D_h Q_{t}\|^2_{L^2 ({\mathbb R}^3_+;{\mathbb R})} +C_{\delta} C(M)\ \int_0^T \|\nabla W_t\|^{0.5}_{L^2 ({\mathbb R}^3_+;{\mathbb R}^9)} \|\nabla W_t\|^{1.5}_{H^1 ({\mathbb R}^3_+;{\mathbb R}^9)}\nonumber\\ .
%\end{align*}
%where we have used (\ref{c0}) for the $L^{\infty}$ control of $\tilde a$ in $H^2$. Thus,
%\begin{align}
%|J_4^3|
&\le \ C\delta \int_0^t \|\nabla \tilde q_{t}\|^2_{L^2 (\Omega_0^f;{\mathbb R}^3)} \nonumber\\
&\qquad +C_{\delta} C(M)\ T^{\frac{1}{4}}\ [N(u_0,f)^2+ T\int_0^T \|\nabla \tilde w_{tt}\|^{2}_{L^2 (\Omega_0^f;{\mathbb R}^9)} + \int_0^T \|\nabla \tilde w_t\|^{2}_{H^1 (\Omega_0^f;{\mathbb R}^9)} ]\ .
\label{J44} 
\end{align}

\noindent{\bf Step 5.} 
Let $\displaystyle J_5=\int_0^T (D_h(\tilde b_k^r \tilde b_k^s) {W_t},_r (\cdot+h),\ D_h {W_t},_s)_{L^2({\mathbb R}^3_+;{\mathbb R}^3)}$. Then,
\begin{align}
|J_5|&\le \ \delta \int_0^T \|D_h \nabla W_{t}\|^2_{L^2 ({\mathbb R}^3_+;{\mathbb R}^9)} +C_{\delta}\int_0^T \|D_h(\tilde b\tilde b^T)\|^2_{L^4({\mathbb R}^3_+;{\mathbb R}^9)}\|\nabla W_t\|^2_{L^4 ({\mathbb R}^3_+;{\mathbb R}^9)}\nonumber\\
%& \le\  \delta \int_0^T \|D_h \nabla W_{t}\|^2_{L^2 ({\mathbb R}^3_+;{\mathbb R}^9)} +C_{\delta} C(M)\ \int_0^T \|\nabla W_t\|^{0.5}_{L^2 ({\mathbb R}^3_+;{\mathbb R}^9)} \|\nabla W_t\|^{1.5}_{H^1 ({\mathbb R}^3_+;{\mathbb R}^9)}\\\ .
%\end{align*}
%where we have used (\ref{c0}) for the $L^{\infty}$ control of $\tilde a$ in $H^2$. Thus,
%\begin{align}
%|J_5|
&\le  \  C\delta \int_0^t \|\tilde w_t\|^2_{ H^2 (\Omega_0^f;{\mathbb R}^3)} \nonumber\\
&\qquad +C_{\delta} C(M)\ T^{\frac{1}{4}}\ [N(u_0,f)^2+ T\int_0^T \|\nabla \tilde w_{tt}\|^{2}_{L^2 (\Omega_0^f;{\mathbb R}^9)} + \int_0^T \|\nabla \tilde w_t\|^{2}_{H^1 (\Omega_0^f;{\mathbb R}^9)} ]\ .
\label{J5} 
\end{align}

\noindent{\bf Step 6.} Let $\displaystyle J_6=\int_0^T (D_h(\tilde b_k^r \tilde b_k^s)_t {W},_r (\cdot+h),\ D_h {W_t},_s)_{L^2({\mathbb R}^3_+;{\mathbb R}^3)}$. Similarly,
\begin{align}
|J_6|&\le \ \delta \int_0^T \|D_h \nabla W_{t}\|^2_{L^2 ({\mathbb R}^3_+;{\mathbb R}^9)} +C_{\delta}\int_0^T \|D_h(\tilde b\tilde b^T)_t\|^2_{L^2({\mathbb R}^3_+;{\mathbb R}^9)}\|\nabla W\|^2_{W^{1,4} ({\mathbb R}^3_+;{\mathbb R}^9)}\nonumber\\
%& \le\  \delta \int_0^T \|D_h \nabla W_{t}\|^2_{L^2 ({\mathbb R}^3_+;{\mathbb R}^9)} +C_{\delta} C(M)\ \int_0^T \|\nabla W\|^{0.5}_{H^1 ({\mathbb R}^3_+;{\mathbb R}^9)} \|\nabla W\|^{1.5}_{H^2 ({\mathbb R}^3_+;{\mathbb R}^9)}\\\ .
%\end{align*}
%where we have used (\ref{c0}) for the $L^{\infty}$ control of $\tilde a$ and $\tilde a_t$ respectively in $H^2$ and $H^1$. Thus,
%\begin{align}
%|J_6|
&\le \ C\delta \int_0^t \|\tilde w_{t}\|^2_{H^2 (\Omega_0^f;{\mathbb R}^3)} \nonumber\\
&\qquad +C_{\delta} C(M)\ T^{\frac{1}{4}}\ [N(u_0,f)^2+ T\int_0^T \|\nabla \tilde w_{t}\|^{2}_{H^1 (\Omega_0^f;{\mathbb R}^9)} + \int_0^T \|\nabla \tilde w\|^{2}_{H^2 (\Omega_0^f;{\mathbb R}^9)} ]\ .
\label{J6} 
\end{align}

\noindent{\bf Step 7.} For $\displaystyle J_7=\int_0^T ((\tilde b_k^r \tilde b_k^s)_t\ D_h{W},_r,\ D_h {W_t},_s)_{L^2({\mathbb R}^3_+;{\mathbb R}^3)}$, we have
\begin{align}
|J_7|&\le \ \delta \int_0^T \|D_h \nabla W_{t}\|^2_{L^2 ({\mathbb R}^3_+;{\mathbb R}^9)} +C_{\delta}\int_0^T \|(\tilde b\tilde b^T)_t\|^2_{L^4({\mathbb R}^3_+;{\mathbb R}^9)}\|D_h \nabla W\|^2_{L^4 ({\mathbb R}^3_+;{\mathbb R}^9)}\nonumber\\
%& \le\  \delta \int_0^T \|D_h \nabla W_{t}\|^2_{L^2 ({\mathbb R}^3_+;{\mathbb R}^9)} +C_{\delta} C(M)\ \int_0^T \|\nabla W\|^{0.5}_{H^1 ({\mathbb R}^3_+;{\mathbb R}^9)} \|\nabla W\|^{1.5}_{H^2 ({\mathbb R}^3_+;{\mathbb R}^9)}\\\ .
%\end{align*}
%where we have used (\ref{c0}) for the $L^{\infty}$ control of $\tilde a$ and $\tilde a_t$ respectively in $H^2$ and $H^1$. Thus,
%\begin{align}
%|J_7|
&\le \ C\delta \int_0^T \|w_{t}\|^2_{H^2 (\Omega_0^f;{\mathbb R}^3)} \nonumber\\
& +C_{\delta} C(M)\ T^{\frac{1}{4}}\ [N(u_0,f)^2+ T\int_0^T \|\nabla \tilde w_t\|^{2}_{H^1 (\Omega_0^f;{\mathbb R}^9)} + \int_0^T \|\nabla \tilde w\|^{2}_{H^2 (\Omega_0^f;{\mathbb R}^9)} ]\ .
\label{J7} 
\end{align}

For the next step, we introduce  $\delta_1>0$, which is different from
$\delta$.

\noindent{\bf Step 8.} Let $\displaystyle J_8=\int_0^T (D_h C^{irks}{W^k},_r,\ D_h {W^i_t},_s)_{L^2({\mathbb R}^3_-;{\mathbb R})}$. 

An integration by parts in time gives
\begin{align*}
J_8=\ - \int_0^T (D_h C^{irks}&\ {W^k_t},_r,\ D_h {W^i},_s)_{L^2({\mathbb R}^3_-;{\mathbb R})} \\
& + [\ (\ D_h C^{irks}\ {W^k},_r (\cdot),\ D_h {W^i}(\cdot),_s)_{L^2({\mathbb R}^3_-;{\mathbb R})}\ ]_0^T\ . 
\end{align*}

Since $\Omega_0^s$ is of class $H^4$,
\begin{align*}
|J_8|\le &\ C T\ [\ \sup_{[0,T]} \|\nabla {W_t}\|^2_{L^2({\mathbb R}^3_-;{\mathbb R}^9)} + \sup_{[0,T]} \|{W}\|^2_{H^2({\mathbb R}^3_-;{\mathbb R}^3)}\ ] + C N(u_0,f)^2\\
& + C_{\delta_1}\ [\ \|\nabla {W} (T)-\nabla W (0)\|^2_{L^2({\mathbb R}^3_-;{\mathbb R}^9)} + \|\nabla W (0)\|^2_{L^2({\mathbb R}^3_-;{\mathbb R}^9)}\ ]\\
& + \delta_1 \sup_{[0,T]} \| D_h \nabla W\|^2_{L^2({\mathbb R}^3_-;{\mathbb R}^9)}, 
\end{align*}
and thus,
\begin{align}
|J_8|\le &\ C T\ [\ \sup_{[0,T]} \|\nabla {\tilde w_t}\|^2_{L^2(\Omega_0^s;{\mathbb R}^9)} + \sup_{[0,T]} \|{\tilde w}\|^2_{H^2(\Omega_0^s;{\mathbb R}^3)}\ ] + C_{\delta_1}\ N(u_0,f)^2\nonumber\\
& + C_{\delta_1}\  T^2\ \sup_{[0,T]} \|\nabla {\tilde w_t}\|^2_{L^2(\Omega_0^s;{\mathbb R}^9)} + C\delta_1 \sup_{[0,T]} \| \tilde w\|^2_{H^2(\Omega_0^s;{\mathbb R}^3)}\ . 
\label{J8}
\end{align}

\noindent{\bf Step 9.} Let $\displaystyle J_9=\int_0^T ( {F_1}_t,\ D_{-h} D_h {W_t})_{L^2({\mathbb R}^3_+;{\mathbb R}^3)}\ .$ Then
\begin{align*}
|J_9|& \le  \ \delta \int_0^T \| W_{t}\|^2_{H^2 ({\mathbb R}^3_+;{\mathbb R}^3)} +C_{\delta}\int_0^T \|(\tilde b\tilde b^T)_t\|^2_{L^4({\mathbb R}^3_+;{\mathbb R}^9)}\|\nabla W\|^2_{L^4 ({\mathbb R}^3_+;{\mathbb R}^9)}\nonumber\\
 &\qquad + C_{\delta}\int_0^T \|(\tilde b\tilde b^T)\|^2_{W^{1,4}({\mathbb R}^3_+;{\mathbb R}^9)}\|\nabla W_t\|^2_{L^2 ({\mathbb R}^3_+;{\mathbb R}^9)}+
C_\delta\ N(u_0,f)^2\nonumber\\
&\qquad + C_{\delta}\int_0^T \|\tilde b_t\|^2_{L^4({\mathbb R}^3_+;{\mathbb R}^9)}\|q\|^2_{L^4 (\Omega_0^f;{\mathbb R})} +
C_{\delta}\int_0^T \|\tilde b\|^2_{L^{\infty}({\mathbb R}^3_+;{\mathbb R}^9)}\|q_t\|^2_{L^2 (\Omega_0^f;{\mathbb R})}\nonumber\\ 
&\le \ C \delta \int_0^T \| \tilde w_{t}\|^2_{H^2 (\Omega_0^f;{\mathbb R}^3)} +
C_\delta\ N(u_0,f)^2\nonumber\\
&\qquad +C_{\delta} C(M)\ T\ [N(u_0,f)^2+ T\int_0^T \|\tilde w_t\|^{2}_{H^2 (\Omega_0^f;{\mathbb R}^3)}+T\ \int_0^T \|\tilde w_{tt}\|^{2}_{H^1 (\Omega_0^f;{\mathbb R}^3)}\ ]\nonumber\\
&\qquad +C_{\delta} C(M)\ T\ [N(u_0,f)^2+ T\ \int_0^T \|\tilde q_t\|^{2}_{H^1 (\Omega_0^f;{\mathbb R})}+   \sup_{[0,T]} \|\tilde q_t\|^{2}_{L^2 (\Omega_0^f;{\mathbb R})}\ ]\ . 
\end{align*}
Thus, with (\ref{eliminateqt}),
\begin{align}
|J_9|\le &\ C \delta 
\int_0^T \| \tilde w_{t}\|^2_{H^2 (\Omega_0^f;{\mathbb R}^3)} 
+ C_\delta\ N(u_0,f)^2\nonumber\\
& +C_{\delta} C(M)\ T\ [N(u_0,f)^2+ T\int_0^T \|\tilde w_t\|^{2}_{H^2 (\Omega_0^f;{\mathbb R}^3)}+T\ \int_0^T \|\tilde w_{tt}\|^{2}_{H^1 (\Omega_0^f;{\mathbb R}^3)}\ ]\nonumber\\
& +C_{\delta} C(M)\ T\ [\ T\ \int_0^T \|\tilde q_t\|^{2}_{H^1 (\Omega_0^f;{\mathbb R})}+   \sup_{[0,T]} \|\tilde w_{tt}\|^{2}_{L^2 (\Omega;{\mathbb R}^3)} +   \sup_{[0,T]} \|\tilde w\|^{2}_{H^1 (\Omega_0^s;{\mathbb R}^3)}\ ] .
\label{J9} 
\end{align}
\noindent{\bf Step 10.} For $\displaystyle J_{10}=\int_0^T ( D_h {H_i}_t,\  D_h {W_t},_i)_{L^2({\mathbb R}^3_+;{\mathbb R}^3)}$, we have
\begin{align}
|J_{10}|& \le  \delta \int_0^T \| W_{t}\|^2_{H^2 ({\mathbb R}^3_+;{\mathbb R}^3)} +C_{\delta}\int_0^T \|(\tilde b\tilde b^T)_t\|^2_{L^4({\mathbb R}^3_+;{\mathbb R}^9)}\|\nabla W\|^2_{L^4 ({\mathbb R}^3_+;{\mathbb R}^9)}\nonumber\\
&\qquad +C_{\delta}\int_0^T \|D_h(\tilde b\tilde b^T)_t\|^2_{L^2({\mathbb R}^3_+;{\mathbb R}^9)}\|W\|^2_{W^{1,4} ({\mathbb R}^3_+;{\mathbb R}^3)}\nonumber\\
&\qquad+ C_{\delta}\int_0^T \|\tilde b\tilde b^T\|^2_{W^{1,4}({\mathbb R}^3_+;{\mathbb R}^9)}\|\nabla W_t\|^2_{L^2 ({\mathbb R}^3_+;{\mathbb R}^9)}\nonumber\\
&\qquad + C_{\delta}\int_0^T \|D_h(\tilde b\tilde b^T)\|^2_{L^4({\mathbb R}^3_+;{\mathbb R}^9)}\|W_t\|^2_{L^4 ({\mathbb R}^3_+;{\mathbb R}^3)}\nonumber\\ 
&\le C\delta \int_0^T \| \tilde w_{t}\|^2_{H^2 (\Omega_0^f;{\mathbb R}^3)} \nonumber\\
&\qquad +C_{\delta} C(M)\ T\ [N(u_0,f)^2+ T\int_0^T \|\tilde w_t\|^{2}_{H^2 (\Omega_0^f;{\mathbb R}^3)}+T\int_0^T \|\tilde w_{tt}\|^{2}_{H^1 (\Omega_0^f;{\mathbb R}^3)}\ ]\ .
\label{J10} 
\end{align}

\noindent{\bf Step 11.} Let $\displaystyle J_{11}=\int_0^T ( D_h {F_2}_t,\ D_h {W_t})_{L^2({\mathbb R}^3_-;{\mathbb R}^3)} + \int_0^T ( D_{-h} D_h {K_i}_t,\ {W_t},_i)_{L^2({\mathbb R}^3_-;{\mathbb R}^3)}\ .$ Then,
\begin{align}
|J_{11}|\le &\ C\ {T}\ [\sup_{[0,T]}\|\tilde w_t\|^{2}_{H^1 (\Omega_0^s;{\mathbb R}^3)}+ \sup_{[0,T]} \|\tilde w\|^{2}_{H^2 (\Omega_0^s;{\mathbb R}^3)}]+C N(u_0,f)^2\ .
\label{J11} 
\end{align}

\subsection{Estimates for (\ref{energywlimit})}\hfill\break
As in the previous section, recall  that from (\ref{c0}), $\tilde a$, and thus 
$\tilde b$, is bounded in $L^{\infty}(H^2)$ independently of the parameter 
$n$ associated to $\tilde a$.

\noindent {\bf Step 1.} For 
$\displaystyle K_1= \int_0^T (D_{-h}D_h (\tilde b_i^j)\ Q,\  D_{-h}D_h W^i,_j)_{L^2({\mathbb R}^3_+;{\mathbb R})}$, we have
 \begin{align}
|K_1|&\le \delta \int_0^T \|D_{-h}D_h \nabla W\|^2_{L^2 ({\mathbb R}^3_+;{\mathbb R}^9)} +C_{\delta}\int_0^T \|D_{-h}D_h\tilde b\|^2_{L^2({\mathbb R}^3_+;{\mathbb R}^9)}\|Q\|^2_{W^{1,4} ({\mathbb R}^3_+;{\mathbb R})}\nonumber\\
&\le \delta \int_0^T \|D_{-h}D_h\nabla W\|^2_{L^2 ({\mathbb R}^3_+;{\mathbb R}^9)} +C_{\delta} C(M)\int_0^T \|Q\|^{0.5}_{H^1 ({\mathbb R}^3_+;{\mathbb R})} \|Q\|^{1.5}_{H^2 ({\mathbb R}^3_+;{\mathbb R})}\nonumber\\ 
&\le \ C \delta \int_0^T \|\tilde w\|^2_{H^3 (\Omega_0^f;{\mathbb R}^3)}\nonumber\\
& +C_{\delta} C(M)\ T^{\frac{1}{4}}\ [\ N(u_0,f)^2+ T\int_0^T\|q_t\|^{2}_{H^1 (\Omega_0^f;{\mathbb R})} + \int_0^T \|q\|^{2}_{H^2 (\Omega_0^f;{\mathbb R})} ]\ .
\label{K1} 
\end{align}

\noindent {\bf Step 2.} Let $\displaystyle K_2= \sum_{p=0}^1\int_0^T (D_{(-1)^p h} \tilde {b}_i^j\ D_{(-1)^{p} h} Q,\ D_{-h} D_h W^i,_j)_{L^2({\mathbb R}^3_+;{\mathbb R})}$. Then,
\begin{align}
|K_2|&\le \delta \int_0^T \|D_{-h}D_h \nabla W\|^2_{L^2 ({\mathbb R}^3_+;{\mathbb R}^9)} +C_{\delta}\int_0^T \|D_h\tilde b\|^2_{L^4({\mathbb R}^3_+;{\mathbb R}^9)}\|D_h Q\|^2_{L^4 ({\mathbb R}^3_+;{\mathbb R})}\nonumber\\
%&\le \delta \int_0^T \|D_{-h} D_h\nabla \tilde W\|^2_{L^2 ({\mathbb R}^3_+;{\mathbb R}^9)} +C_{\delta} C(M)\int_0^T \|Q\|^{0.5}_{H^1 ({\mathbb R}^3_+;{\mathbb R})} \|Q\|^{1.5}_{H^2 ({\mathbb R}^3_+;{\mathbb R})}\ ,
&\le C\delta \int_0^T \|\tilde w\|^2_{H^3 (\Omega_0^f;{\mathbb R}^3)} \nonumber\\
& +C_{\delta} C(M)\ T^{\frac{1}{4}}\ [\ N(u_0,f)^2+ T\int_0^T \|q_t\|^{2}_{H^1 (\Omega_0^f;{\mathbb R})} + \int_0^T \|q\|^{2}_{H^2 (\Omega_0^f;{\mathbb R})} ]\ .
\label{K2} 
\end{align}

\noindent {\bf Step 3.} Let $\displaystyle K_3= \int_0^T ( \tilde {b}_i^j\ D_{-h} D_h Q,\  D_{-h} D_h W^i,_j)_{L^2({\mathbb R}^3_+;{\mathbb R})}$. 

We first notice that
\begin{align*}
K_3=&\ \int_0^T (  D_{-h} D_h Q,\  D_{-h} D_h [\tilde {b}_i^j W^i,_j])_{L^2({\mathbb R}^3_+;{\mathbb R})}\nonumber\\
&+\sum_{p=0}^1  \int_0^T ( D_{-h} D_h Q,\  D_{(-1)^p h} \tilde {b}_i^j \  D_{(-1)^p h} W^i,_j)_{L^2({\mathbb R}^3_+;{\mathbb R})}\\
&-  \int_0^T ( D_{-h} D_h Q,\  D_h D_{-h}[\tilde {b}_i^j] \  W^i,_j)_{L^2({\mathbb R}^3_+;{\mathbb R})}\ ,
\end{align*}
which with the divergence relation (\ref{divchange}) leads us to
\begin{align*}
K_3=&\ \int_0^T (  D_{-h} D_h Q,\  D_{-h} D_h \mathfrak a)_{L^2({\mathbb R}^3_+;{\mathbb R})}\nonumber\\
&+\sum_{p=0}^1  \int_0^T ( D_{-h} D_h Q,\  D_{(-1)^p h} \tilde {b}_i^j \  D_{(-1)^p h} W^i,_j)_{L^2({\mathbb R}^3_+;{\mathbb R})}\\
&-  \int_0^T ( D_{-h} D_h Q,\  D_h D_{-h} [\tilde {b}_i^j] \  W^i,_j)_{L^2({\mathbb R}^3_+;{\mathbb R})}\ .
\end{align*}

We then have that
\begin{align}
|K_3|&\le \ \delta \int_0^T \|D_{-h} D_h Q\|^2_{L^2 ({\mathbb R}^3_+;{\mathbb R})}  +C_{\delta}\int_0^T \|\tilde b\|^2_{H^2({\mathbb R}^3_+;{\mathbb R}^9)}\|\tilde w\|^2_{H^2(\Omega_0^f;{\mathbb R}^3)}\nonumber\\
&\qquad+C_{\delta}\int_0^T \|D_h\tilde b\|^2_{L^4({\mathbb R}^3_+;{\mathbb R}^9)}\|D_h\nabla W\|^2_{L^{4} ({\mathbb R}^3_+;{\mathbb R}^9)}\nonumber\\
&\qquad +C_{\delta}\int_0^T \|\tilde b\|^2_{H^2({\mathbb R}^3_+;{\mathbb R}^9)}\|\nabla W \|^2_{W^{1,4} ({\mathbb R}^3_+;{\mathbb R}^9)}\nonumber\\
&\le\ C \delta \int_0^T \|\tilde q\|^2_{H^2 (\Omega_0^f;{\mathbb R})} +C_{\delta} C(M)\ T\  [\ N(u_0,f)^2+T \int_0^T\|\tilde w_t\|^{2}_{H^2 (\Omega_0^f;{\mathbb R}^3)}\ ]\nonumber\\
&\qquad +C_{\delta} C(M)\ T^{\frac{1}{4}}\  [\ N(u_0,f)^2+T\int_0^T\|\tilde w_t\|^{2}_{H^2 (\Omega_0^f;{\mathbb R}^3)} +\int_0^T \|\tilde w\|^{2}_{H^3 (\Omega_0^f;{\mathbb R}^3)}\ ]
\ .
\label{K3}
\end{align}

\noindent{\bf Step 4.} Let $\displaystyle K_4=\int_0^T (D_{-h}D_h(\tilde b_k^r \tilde b_k^s) {W},_r,\ D_{-h}D_h {W},_s)_{L^2({\mathbb R}^3_+;{\mathbb R}^3)}$. Then,
\begin{align}
|K_4|&\le \ \delta \int_0^T \|D_{-h} D_h \nabla W\|^2_{L^2 ({\mathbb R}^3_+;{\mathbb R}^9)} +C_{\delta}\int_0^T \|\tilde b\tilde b^T\|^2_{H^2({\mathbb R}^3_+;{\mathbb R}^9)}\|\nabla W\|^2_{W^{1,4} ({\mathbb R}^3_+;{\mathbb R}^9)}\nonumber\\
%& \le\  \delta \int_0^T \|D_{-h} D_h \nabla W\|^2_{L^2 ({\mathbb R}^3_+;{\mathbb R}^9)} +C_{\delta} C(M)\ \int_0^T \|\nabla W\|^{0.5}_{L^2 ({\mathbb R}^3_+;{\mathbb R}^9)} \|\nabla W\|^{1.5}_{H^1 ({\mathbb R}^3_+;{\mathbb R}^9)}\nonumber\\ 
&\le  C\delta \int_0^T \|\tilde w\|^2_{H^3 (\Omega_0^f;{\mathbb R}^3)} \nonumber\\
&\qquad +C_{\delta} C(M)\ T^{\frac{1}{4}}\ [N(u_0,f)^2+ T\int_0^T \|\nabla \tilde w_{t}\|^{2}_{H^1 (\Omega_0^f;{\mathbb R}^9)} + \int_0^T \|\nabla \tilde w\|^{2}_{H^2 (\Omega_0^f;{\mathbb R}^9)} ]\ .
\label{K4} 
\end{align}

\noindent{\bf Step 5.} For $\displaystyle K_5=\sum_{p=0}^1 \int_0^T (D_{(-1)^p h} (\tilde b_k^r \tilde b_k^s) D_{(-1)^{p}h}{W},_r,\ D_{-h}D_h {W},_s)_{L^2({\mathbb R}^3_+;{\mathbb R}^3)}$ ,
\begin{align}
|K_5|&\le \ \delta \int_0^T \|D_{-h}D_h \nabla W\|^2_{L^2 ({\mathbb R}^3_+;{\mathbb R}^9)} +C_{\delta}\int_0^T \|D_h(\tilde b\tilde b^T)\|^2_{L^4({\mathbb R}^3_+;{\mathbb R}^9)}\|D_h \nabla W\|^2_{L^4({\mathbb R}^3_+;{\mathbb R}^9)}\nonumber\\
%& \le\  \delta \int_0^T \|D_{-h} D_h \nabla W\|^2_{L^2 ({\mathbb R}^3_+;{\mathbb R}^9)} +C_{\delta} C(M)\ \int_0^T \|W\|^{0.5}_{H^2 ({\mathbb R}^3_+;{\mathbb R}^3)} \| W\|^{1.5}_{H^3 ({\mathbb R}^3_+;{\mathbb R}^3)}\ ,
& \le \delta \int_0^t \|\tilde w\|^2_{H^3(\Omega_0^f;{\mathbb R}^3)} \nonumber\\
&\qquad +C_{\delta} C(M)\ T^{\frac{1}{4}}\ [N(u_0,f)^2+ T\int_0^T \| \tilde w_{t}\|^{2}_{H^2 (\Omega_0^f;{\mathbb R}^3)} + \int_0^T \| \tilde w\|^{2}_{H^3 (\Omega_0^f;{\mathbb R}^3)} ]\ .
\label{K5} 
\end{align}

\noindent{\bf Step 6.} For $\displaystyle K_6=\int_0^T (D_{-h} D_h [C^{irks}]\ \int_0^{\cdot} {W^k},_r,\ D_{-h}D_h {W^i},_s)_{L^2({\mathbb R}^3_-;{\mathbb R})}$,
an integration by parts in time gives
\begin{align*}
K_6=&\ - \int_0^T (D_{-h}D_h [C^{irks}]\ {W^k},_r,\ D_{-h} D_h \int_0^{\cdot} {W^i},_s)_{L^2({\mathbb R}^3_-;{\mathbb R})} \\
& +  (\ D_{-h}D_h [C^{irks}]\ \int_0^T {W^k},_r ,\ D_{-h} D_h \int_0^T {W^i},_s)_{L^2({\mathbb R}^3_-;{\mathbb R})}\ , 
\end{align*}
from which we infer from the $H^4$ regularity of $\Omega_0^s$,
\begin{align*}
K_6\le &\ C\ \int_0^T  \|\nabla {W}\|_{L^4({\mathbb R}^3_-;{\mathbb R}^9)}\| D_{-h} D_h \int_0^{\cdot} {\nabla W}\|_{L^2({\mathbb R}^3_-;{\mathbb R}^9)} \\
& +  \|\int_0^T \nabla{W}\|_{L^4({\mathbb R}^3_-;{\mathbb R}^9)} \| D_{-h} D_h \int_0^T \nabla {W}\|_{L^2({\mathbb R}^3_-;{\mathbb R}^9)}\ , 
\end{align*}
leading us to
\begin{align}
|K_6|\le &\ C T\ [\ \sup_{[0,T]} \| \tilde w\|^2_{H^2(\Omega_0^s;{\mathbb R}^3)} + \sup_{[0,T]} \|\int_0^{\cdot} {\tilde w}\|^2_{H^3(\Omega_0^s;{\mathbb R}^3)}\ ] \ . 
\label{K6}
\end{align}

\noindent{\bf Step 7.} Let $$\displaystyle K_7=\sum_{p=0}^1\int_0^T ( D_{(-1)^p h} [C^{irks}]\ \int_0^{\cdot} D_{(-1)^{p} h} {W^k},_r,\ D_{-h}D_h {W^i},_s)_{L^2({\mathbb R}^3_-;{\mathbb R})}\ .$$ 

An integration by parts in time gives
\begin{align*}
K_7=&\ - \sum_{p=0}^1 \int_0^T (D_{(-1)^p h} [C^{irks}]\ D_{(-1)^{p+1} h} {W^k},_r,\ D_{-h} D_h \int_0^{\cdot} {W^i},_s)_{L^2({\mathbb R}^3_-;{\mathbb R})} \\
& +  (\ D_{(-1)^{p} h} [C^{irks}]\ D_{(-1)^{p} h} \int_0^T {W^k},_r ,\ D_{-h} D_h \int_0^T {W^i},_s)_{L^2({\mathbb R}^3_-;{\mathbb R})}\ , 
\end{align*}
and thus from the $H^4$ regularity of $\Omega_0^s$,
\begin{align*}
K_7\le &\ C\ \int_0^T  \|D_h \nabla {W}\|_{L^2({\mathbb R}^3_-;{\mathbb R}^9)}\| D_{-h} D_h \int_0^{\cdot} {\nabla W}\|_{L^2({\mathbb R}^3_-;{\mathbb R}^9)} \\
& +  \|D_h \int_0^T \nabla{W}\|_{L^2({\mathbb R}^3_-;{\mathbb R}^9)} \| D_{-h} D_h \int_0^T \nabla {W}\|_{L^2({\mathbb R}^3_-;{\mathbb R}^9)}\ . 
\end{align*}
Therefore, 
\begin{align}
|K_7|\le &\ C T\ [\ \sup_{[0,T]} \| \tilde w\|^2_{H^2(\Omega_0^s;{\mathbb R}^3)} + \sup_{[0,T]} \|\int_0^{\cdot} {\tilde w}\|^2_{H^3(\Omega_0^s;{\mathbb R}^3)}\ ] \ . 
\label{K7}
\end{align}

\begin{remark}
The $H^4$ regularity of $\Omega_0^s$ is used only for proving (\ref{K6}) and (\ref{K7}). As a matter of fact, $W^{3,p}$ for $p>3$ would have been sufficient.
\end{remark}

\noindent{\bf Step 8.} Let $\displaystyle K_8=\int_0^T ( D_{-h} F_1,\  D_h D_{-h} D_h {W})_{L^2({\mathbb R}^3_+;{\mathbb R}^3)}\ .$ Then
\begin{align}
|K_8|&\le \delta \int_0^T \| W\|^2_{H^3 ({\mathbb R}^3_+;{\mathbb R}^3)} +C_{\delta}\int_0^T \|D_h(\tilde b\tilde b^T)\|^2_{L^4({\mathbb R}^3_+;{\mathbb R}^9)}\|\nabla \tilde w\|^2_{L^4 (\Omega_0^f;{\mathbb R}^9)}\nonumber\\
&\qquad + C_{\delta}\int_0^T \|\tilde b\tilde b^T\|^2_{L^{\infty}({\mathbb R}^3_+;{\mathbb R}^9)}\|D_h \nabla \tilde w\|^2_{L^2 (\Omega_0^f;{\mathbb R}^9)}+C\ N(u_0,f)^2\nonumber\\
&\qquad + C_{\delta}\int_0^T \|\nabla\tilde b\|^2_{L^4({\mathbb R}^3_+;{\mathbb R}^{27})}\|q\|^2_{L^4 (\Omega_0^f;{\mathbb R})} +
C_{\delta}\int_0^T \|\tilde b\|^2_{L^{\infty}({\mathbb R}^3_+;{\mathbb R}^9)}\| \nabla q\|^2_{L^2 (\Omega_0^f;{\mathbb R}^3)}\nonumber\\ 
&\le  \delta \int_0^T \| \tilde w\|^2_{H^3 (\Omega_0^f;{\mathbb R}^3)} +C_{\delta} C(M)\ T\ [\ N(u_0,f)^2+T\ \int_0^T \|\tilde w_{t}\|^{2}_{H^2 (\Omega_0^f;{\mathbb R}^3)}\ ]\nonumber\\
&\qquad +C_{\delta} C(M)\ T\ [N(u_0,f)^2+ T\ \int_0^T \|q_t\|^{2}_{H^1 (\Omega_0^f;{\mathbb R})}\ ] +C\ N(u_0,f)^2\ .
\label{K_8} 
\end{align}

\noindent{\bf Step 9.} For $\displaystyle K_9=\int_0^T ( D_{-h}D_h H_i,\  D_{-h}D_h {W},_i)_{L^2({\mathbb R}^3_+;{\mathbb R}^3)}$, we have
\begin{align}
|K_9|& \le  \delta \int_0^T \| W\|^2_{H^3 ({\mathbb R}^3_+;{\mathbb R}^3)} +C_{\delta}\int_0^T \|D_h (\tilde b\tilde b^T)\|^2_{L^4({\mathbb R}^3_+;{\mathbb R}^9)}\|\nabla \tilde w\|^2_{L^4 (\Omega_0^f;{\mathbb R}^9)}\nonumber\\
&\qquad +C_{\delta}\int_0^T \|D_{-h} D_h(\tilde b\tilde b^T)\|^2_{L^2({\mathbb R}^3_+;{\mathbb R}^9)}\|\tilde w\|^2_{W^{1,4} (\Omega_0^f;{\mathbb R}^3)}\nonumber\\
&\qquad+ C_{\delta}\int_0^T \|\tilde b\tilde b^T\|^2_{L^{\infty}({\mathbb R}^3_+;{\mathbb R}^9)}\| \tilde w \|^2_{H^2 (\Omega_0^f;{\mathbb R}^3)}\nonumber\\ & \le \delta \int_0^T \| \tilde w\|^2_{H^3 (\Omega_0^f;{\mathbb R}^3)} +C_{\delta} C(M)\ T\ [N(u_0,f)^2+ T\int_0^T \|\tilde w_t\|^{2}_{H^2 (\Omega_0^f;{\mathbb R}^3)}\ ]\ .
\label{K9} 
\end{align}

\noindent{\bf Step 10.} Let \begin{align*}
\displaystyle K_{10}=\int_0^T (D_{-h} D_h F_2,&\ D_{-h} D_h {W})_{L^2({\mathbb R}^3_-;{\mathbb R}^3)}\\
& + \int_0^T ( D_{-h} D_{h} D_{-h} K_i,\ D_h {W},_i)_{L^2({\mathbb R}^3_-;{\mathbb R}^3)}\ .
\end{align*}
 Then
\begin{align}
|K_{10}|\le &\ C {T} \ [\ \sup_{[0,T]}\|\tilde w\|^{2}_{H^2(\Omega_0^s;{\mathbb R}^3)}+ \sup_{[0,T]} \|\int_0^{\cdot}\tilde w\|^{2}_{H^3 (\Omega_0^s;{\mathbb R}^3)}\ ]+ C\ N(u_0,f)^2\ .
\label{K10} 
\end{align}

\section*{Acknowledgments}
SS was partially supported by  National Science Foundation under grants DMS-0105004 and 
NSF ITR-0313370.

\end{document}